%% file: main.tex
\documentclass{article} 
\usepackage[accepted]{icml2026}

\usepackage{natbib}

\input{math_commands.tex}

\newcommand*{\PP}{\text{P}}
\usepackage{xurl}
\usepackage[backref=page]{hyperref}
\PassOptionsToPackage{linktocpage}{hyperref}
\usepackage{url}
\usepackage{sidecap}
\usepackage{cancel}
\sidecaptionvpos{figure}{t}
\usepackage{overpic}
\usepackage{tcolorbox}
\newcommand{\de}{\text{d}}
\renewcommand*{\backref}[1]{}
\renewcommand*{\backrefalt}[4]{%
  \ifcase#1
    (not cited)%
  \or
    \textcolor{gray}{\emph{(cited on page:~#2)}}%
  \else
    \textcolor{gray}{\emph{(cited on pages:~#2)}}%
  \fi
}
\usepackage[table]{xcolor}
\usepackage{xr}
\usepackage{mathrsfs}
\usepackage{titletoc}
\usepackage{multicol}
\usepackage{mathtools, cuted}
\usepackage{caption}

\DeclareCaptionLabelFormat{center}{\hspace*{\fill}\textbf{#1 #2}\hspace*{\fill}}
\PassOptionsToPackage{linktocpage}{hyperref}
\usepackage{booktabs}
\newcommand{\rad}[1]{\renewcommand{\arraystretch}{#1}}
\usepackage{xcolor}
\definecolor{darkpurple}{rgb}{0.5,0,0.5}
\definecolor{mydarkgrey}{rgb}{0.27, 0.27, 0.27}
\definecolor{antiquefuchsia}{rgb}{0.57, 0.36, 0.51}
\definecolor{cadetblue}{rgb}{0., 0.70, 0.70}
\definecolor{brightmaroon}{rgb}{0.76, 0.13, 0.28}
\definecolor{mymaroon}{rgb}{0.62, 0.05, 0.19}
\definecolor{mymaroon2}{cmyk}{0.2, 0.32, 0.02, 0.068}
\definecolor{darkpurple}{rgb}{0.5,0,0.5}
\definecolor{mydarkgrey}{rgb}{0.47, 0.47, 0.47}
\definecolor{mymagenta}{rgb}{0.8125, 0, 0.8125}

\usepackage{nicefrac}   
\usepackage{graphicx} 
\usepackage{bibunits}
\usepackage{indentfirst}
\usepackage{makecell}
\usepackage{graphicx} 
\usepackage[linesnumbered,ruled,vlined]{algorithm2e}

\usepackage{xparse}

\newif\iffirstitem
\firstitemtrue

\SetCommentSty{mycommfont}
\newcommand{\nosemic}{\renewcommand{\@endalgocfline}{\relax}}
\newcommand{\dosemic}{\renewcommand{\@endalgocfline}{\algocf@endline}}
\let\oldnl\nl
\newcommand{\nonl}{\renewcommand{\nl}{\let\nl\oldnl}}

\definecolor{orcidlogocol}{HTML}{A6CE39}
\definecolor{darkred}{RGB}{135, 33, 9}
\definecolor{mygreen}{HTML}{f0f1ff}   
\usepackage[revision]{revdiff}
\usepackage{url}
\tcbuselibrary{breakable, skins}

\hypersetup{colorlinks,urlcolor=darkpurple, citecolor=mymaroon2, linkcolor=cadetblue}

\icmltitlerunning{From Geometry to Dynamics}

\usepackage{wrapfig}
\usepackage{overpic}

\begin{document}

\twocolumn[
  \icmltitle{From Geometry to Dynamics: Learning Overdamped Langevin Dynamics from Sparse Observations with Geometric Constraints }

  \begin{icmlauthorlist}
    \icmlauthor{Dimitra Maoutsa}{yyy}
  \end{icmlauthorlist}

  \icmlaffiliation{yyy}{Technical University of Berlin, Berlin, Germany}

  \icmlcorrespondingauthor{Dimitra Maoutsa}{dimitra.maoutsa@gmail.com}
  
\icmlkeywords{ geometry, dynamical systems, inference}

  \vskip 0.3in
]
\printAffiliationsAndNotice{}


\begin{abstract}
How can we learn the laws underlying the dynamics of stochastic systems when their trajectories are sampled sparsely in time? Existing methods either require temporally resolved high-frequency observations, or rely on geometric arguments that apply only to conservative systems, limiting the range of dynamics they can recover. Here, we present a new framework that reconciles these two perspectives by reformulating inference as a stochastic control problem. Our method uses geometry-driven path augmentation, guided by the geometry in the system’s invariant density to reconstruct likely trajectories and infer the underlying dynamics without assuming specific parametric models. Applied to overdamped Langevin systems, our approach accurately recovers stochastic dynamics even from extremely undersampled data, outperforming existing methods in synthetic benchmarks. This work 
demonstrates the effectiveness of incorporating geometric inductive biases into stochastic system identification methods.

\end{abstract}

\section{Introduction}\vspace{-2.5pt}
How can we discover the underlying driving forces that govern the behaviour of complex, stochastic systems when we only measure their state at discrete time points? From pollen motion in a liquid medium~\citep{einstein1905molekularkinetischen} and  chemical reactions~\citep{li2020chemical} to population dynamics~\citep{silva2018spontaneous, fisher2014transition} and cell growth~\citep{alonso2014modeling}, many natural processes evolve following stochastic dynamics, best described by Langevin or stochastic differential equations (SDEs) of the form  
\vspace{-2.5pt}
\begin{equation} \label{eq:system}
    \text{d}\mathbf{X}_t = {\mathbf{f}(\mathbf{X}_t)} \,\text{d}t + {\boldsymbol{\sigma}} \,\text{d}\mathbf{W}_t.
\end{equation}
\vspace{-1.5pt}Under this formalism, the deterministic part of the equation $\mathbf{f}(\cdot):\mathbb{R}^d \rightarrow \mathbb{R}^d$, the \emph{drift} function, captures the long-term evolution of the state variables $\mathbf{X}_t \in \mathbb{R}^d$, while the stochastic part $\boldsymbol{\sigma}:\mathbb{R}^{d \times d}$, the \emph{diffusion}, accounts for the contribution of unresolved degrees of freedom.

A common challenge for most inference methods is the assumption of fine temporal resolution of observations.  Most approaches estimate $\mathbf{f}$ from short-time increments or by approximating transition densities over  a small time interval (Sec~\ref{appsec:related_work}).
However, when the sampling interval $\tau$ is large, these approximations break down, and the inverse problem becomes poorly constrained: multiple drifts can induce similar transition statistics between sparse observations.
Thus, in such a setting, to accurately recover $\mathbf{f}$ requires additional assumptions or inductive biases consistent with the data.

\begin{figure*}[t!]
  \hspace{-23pt}
   \centering
    \includegraphics[width=0.87\linewidth]{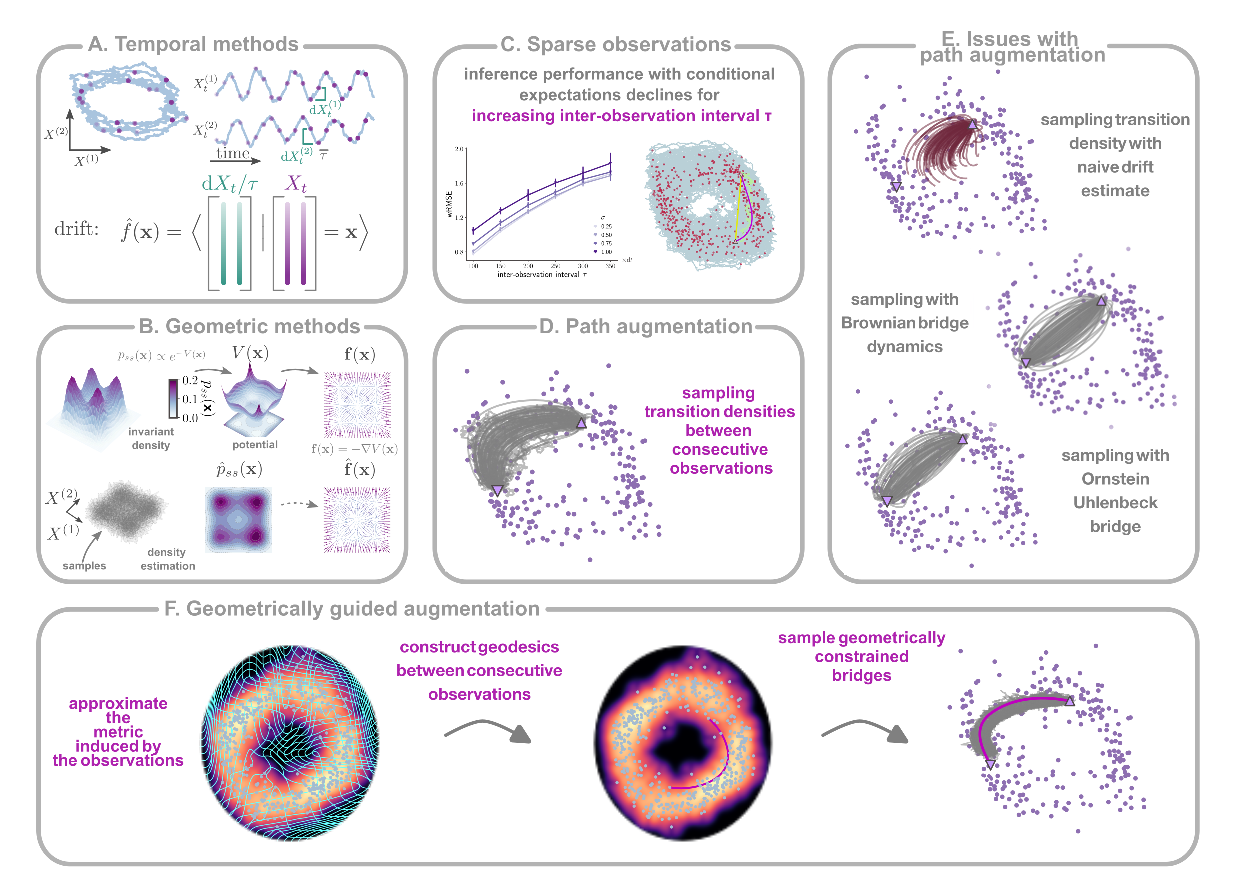} 
    \captionsetup{width=1.00\linewidth}
    \caption{{\textbf{Temporal and geometric perspectives for discovering stochastic dynamics and proposed inference with geometrically guided augmentation.}} 
    (\textbf{A.}) Temporal methods consider the time-ordering of observations $\{\boldsymbol{\mathcal{O}}_k\}^K_{k=1}$ (\emph{purple dots}) to approximate the drift $\mathbf{f}(\mathbf{x})$ with conditional rescaled state increments $\hat{\mathbf{f}}(\mathbf{x}) = \langle \frac{\text{d} \mathbf{X}_t}{\tau} | \mathbf{X}_t=\mathbf{x} \rangle$. (\textbf{B.}) Geometric methods assume a conservative drift $\mathbf{f}(\mathbf{x})=-\nabla V(\mathbf{x})$
    as 
    the gradient of a potential. (\textbf{C.}) 
    With increasing inter-observation interval $\tau$ performance of temporal methods degrades because Euclidean distances ignore the curvature of the latent continuous path between consecutive observations.
    (\textbf{D.}) Path augmentation alternates between state estimation - by sampling diffusion bridges for each inter-observation interval - and drift inference. (\textbf{E.}) Commonly used path augmentation methods employ Brownian or Ornstein-Uhlenbeck bridges that 
    increasingly deviate from the unobserved path as $\tau$ grows. (lower) Illustration of the ground truth (\emph{neon green}) and geodesic (\emph{magenta}) continuous path between two observations and of that assumed during inference with Gaussian likelihood (\emph{yellow line}).  (\textbf{F.}) Geometrically guided augmentation approximates first the metric induced by the invariant density, constructs geodesics connecting consecutive observations, and samples geometrically constrained diffusion bridges. }
    \label{fig1}
  \vspace{-5pt}
  \centering
\end{figure*}

\vspace{-5pt}
 \textbf{Two dominant perspectives for stochastic inference.}  
 Existing data-driven approaches for learning dynamical systems provide powerful tools in the deterministic setting and for densely sampled stochastic dynamics (Sec~\ref{appsec:related_work}).
Broadly, data-driven system identification for stochastic systems largely follows two tracks: \textcolor{mymaroon}{\emph{\textbf{Temporal methods}}} (Fig.~\ref{fig1}\textbf{A.}) rely on the \textbf{temporal ordering} of measurements, regressing state increments against states to estimate the drift. This works well when the inter-observation interval 
 ($\tau$) is small~\citep{ruttor2013approximate, friedrich1997description, friedrich2011approaching, ragwitz2001indispensable, jacobs2023hypersindy}. 
\textcolor{antiquefuchsia}{\emph{\textbf{Geometric methods}}} on the other hand, approximate the \textbf{invariant density}~\citep{batz2016variational,guharlim2021stationary} or eigenstructure of the infinitesimal generator of the diffusion process~\citep{singer2008non,nuske2021spectral,ionides2006inference,talmon2015intrinsic,dsilva2016data,berry2018iterated}) (Fig.~\ref{fig1}\textbf{B.}). 
Each perspective has limitations: temporal approaches deteriorate with increasing inter-observation intervals (Fig.~\ref{fig1}\textbf{C.}), whereas geometric methods are restricted to systems with conservative forces~\citep{berry2015nonparametric,batz2016variational} or decoupled state variables~\citep{singer2008non}.


 
 


\vspace{-5.5pt}
\begin{tcolorbox}[breakable, enhanced jigsaw,width=1.05\linewidth]
\paragraph{A unifying perspective: reconcile temporal and geometric methods by constraining with most probable paths extracted from the invariant density.}  
Here, we recast inference into a stochastic control problem and introduce \textbf{geometry-aware path augmentation}. Our method follows a simple premise that incorporates \textbf{geometric inductive biases} informed by the system's \emph{invariant density} into dynamical inference: we postulate that the augmented paths should lie \textbf{in the vicinity} of \textbf{geodesic curves} (Fig.~\ref{fig1}\textbf{F.}~middle, magenta line) that connect consecutive measurements on the \textbf{empirical manifold} induced by the observations. To achieve this, \textbf{(i)} we approximate the Riemannian metric induced by the observations (Fig.~\ref{fig1}\textbf{F.}) without the need to predefine the dimensionality of the empirical manifold, \textbf{(ii)} compute geodesics between consecutive observations through nonparametric approximation of shortest path distances between consecutive observations according to the approximated metric, and \textbf{(iii)} estimate the unobserved path between consecutive observations by generating \textbf{geometrically constrained diffusion bridges}
that both respect temporal order and are guided toward identified geodesics (Fig.~\ref{fig1} \textbf{F.}). Nonparametric estimation of the drift function based on the augmented paths within an Expectation Maximisation framework (\textbf{E.M.})~\citep{dempster1977maximum} results in accurate approximations of the underlying stochastic dynamics.
Extensive numerical experiments demonstrate the effectiveness of our proposed method in recovering the ground truth stochastic dynamics, even in challenging scenarios where existing approaches fail.
\end{tcolorbox}

\hspace{-23pt}

\vspace{-20pt}

  \begin{figure*}[t]
  \begin{overpic}[width=1.0\textwidth]{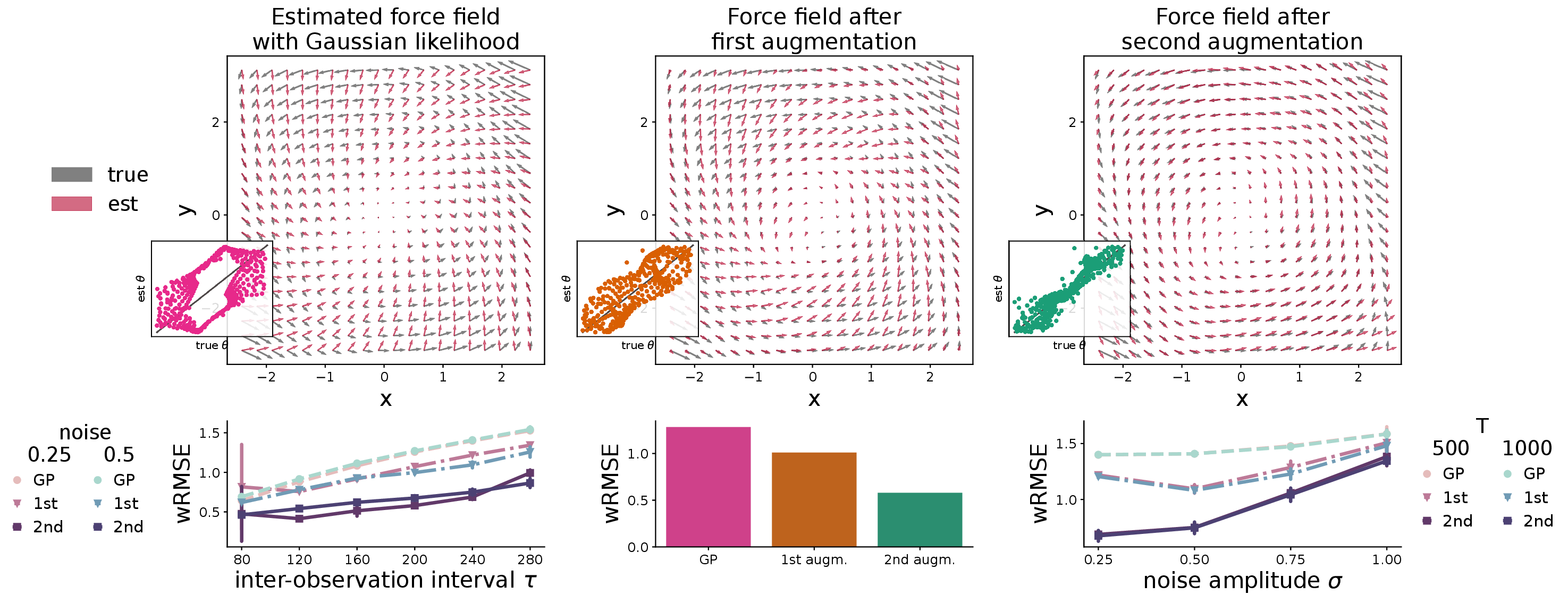}
  \put(6,35){a.}
  \put(36,35){b.}
  \put(66,35){c.}
  \put(6,12){d.}
  \put(36,12){e.}
  \put(66,12){f.}
  \put(35,2.5){\fontsize{6pt}{7pt}{$\times \de t$}}
  \end{overpic}
  \captionsetup{width=0.925\textwidth}
\caption{
\textbf{Geometry-aware path augmentation improves drift inference after two iterations.} 
Estimated (\emph{red}) vs. ground truth (\emph{grey}) force field with \textbf{a.)} Gaussian likelihood, 
\textbf{b.)} after one, and \textbf{c.)} after two augmentations. (\textbf{Insets}) Ground truth vs. estimated angles at grid points. 
\textbf{d.)} Weighted (by observation density) root mean square error (wRMSE) vs. inter-observation interval $\tau$ for different noise levels $\sigma=\{0.25,0.5\}$ for drift estimated with a Gaussian likelihood (\emph{gaus}-circles), after first augmentation (\emph{1st}-triangles), and after second augmentation (\emph{2nd}-squares) for $T=500$ (time units). 
\textbf{e.)} wRMSE across iterations for the presented example. 
\textbf{f.)} wRMSE vs. noise amplitude $\sigma$ for different trajectory durations ${T=\{500,1000\}}$ (time units) for inter-observation interval $\tau=240$ ($dt$).
Markers in \textbf{d.)} and \textbf{f.)} indicate augmentation steps. Error bars: one standard deviation over five independent runs. ($\de t = 0.01)$
}
  \label{fig:res} \vspace{-0pt}
\end{figure*}
Our contributions in this work are the following:\\
\textbf{1.}~We introduce a geometry-guided augmentation framework for inference of sparsely observed SDEs that reconciles temporal and geometric methods.\\
\textbf{2.}~We formulate the geometry-constrained augmentation in terms of a stochastic control problem that can be solved independently across inter-observation intervals.\\
\textbf{3.}~We provide theoretical support showing that for long inter-observation intervals, the bias of short-time approximations is controlled by higher-order terms involving the curvature of the vector field, explaining why purely temporal methods degrade as $\tau$ increases (Sec.~\ref{app:curvature}).\\
\textbf{4.}~We integrate this augmentation into a model-independent drift estimation framework (using a nonparametric function approximation), and demonstrate improved recovery of dynamics in challenging large-$\tau$ settings.

\vspace{-10pt}

\section{Setup and background}
\vspace{-8pt}
\textbf{Setting.}~We consider a system whose state evolves according to Eq.~\ref{eq:system}.
We observe the system state at discrete time points $t_k = k\tau$ at \textbf{inter-observation intervals} of $\tau$ time units, obtaining a time-ordered set of observations $\{\boldsymbol{\mathcal{O}}_k \,\dot= \,\mathbf{X}_{t_k}\}^K_{k=1}$.\vspace{0.75pt}\\
\textbf{Goal.}~Given $\{\boldsymbol{\mathcal{O}}_k\}_{k=1}^K$, estimate the drift $\mathbf{f}$.\\
\textbf{Background.}
Many drift-learning approaches~\citep{batz2018approximate, friedrich1997description} rely on {short-time} approximations: for a small $\tau$ they assume \vspace{-5pt}
\begin{equation}\label{eq:short_time_gaussian}
\mathbf{X}_{t+\tau}\mid \mathbf{X}_t \approx \mathcal{N}\!\big(\mathbf{X}_t + \mathbf{f}(\mathbf{X}_t)\tau,\; \boldsymbol{\sigma} \boldsymbol{\sigma}^\top\tau\big),
\end{equation}
which leads to a quadratic objective in $\mathbf{f}$ and motivates nonparametric estimators such as Gaussian-process drift models~\citep{ruttor2013approximate, hostettler2018modeling, zhao2020state}.
However, in our considered setting the sampling interval $\tau$ is \textbf{not} small: the Gaussian increment model of Eq.~\ref{eq:short_time_gaussian}
(and Euler--Maruyama-based likelihoods) becomes inaccurate because the underlying transition density is generally non-Gaussian (see more details in Sec~\ref{appsec:b}).

Common drift inferences approaches rely on \emph{small} inter-observation intervals $\tau$~\citep{batz2018approximate}. As $\tau$ increases, the EuM approximation becomes inaccurate: transition densities are not Gaussian, and higher-order remainder terms related to the curvature of the flow field become important (see further theoretical analysis in Sec.~\ref{app:curvature} and c.f. Fig.~\ref{fig:increase}).
Attempts to mitigate this problem by introducing bridge sampling to infer the unobserved path between observations~\citep{batz2018approximate, sermaidis2013markov} provide small improvements, because these methods rely on {linearised} or otherwise simplified bridge dynamics that do not match the true transition densities (c.f. Sec.~\ref{appsec:related_work}).

\textbf{Our approach.} Here, we target 
 this large inter-observation interval setting
 by merging insights from both temporal and geometric perspectives. Specifically, our approach combines \textbf{nonlinear} bridge sampling with \textbf{a geometric approximation of the system’s invariant density} as detailed in the following.
This introduces an implicit bias in our approach: among the alternative potential hypotheses regarding the underlying drift dynamics compatible with the data, we favour those whose likely paths follow
high-density regions on the invariant density.


\vspace{-12pt}
\section{Methodology}

\textbf{Core idea.}~The invariant density of the observed system imposes a low-dimensional structure on the state space, within which the observations are confined. We propose that this low-dimensional structure is well approximated by a Riemannian manifold $\mathcal{M}_{\infty} \in \mathbb{R}^{m \leq d}$ in the ambient space (see Sec.~\ref{appsec:theory_low_dim}), and that the ensemble of observations $\{\boldsymbol{\mathcal{O}}_k\}_{k=1}^{K}$ offers a reliable discrete approximation to $\mathcal{M}_{\infty}$. 
We term this observation-based approximation the \textit{empirical manifold} $\mathcal{M}$. The central premise of our approach is that \textbf{unobserved paths between successive observations will be lying either \emph{on} or \emph{in the vicinity} of the empirical manifold} $\mathcal{M}$. In particular, we postulate that unobserved paths should lie \textbf{in the vicinity of geodesics that connect consecutive observations} on $\mathcal{M}$.

However, while this view of a lower dimensional manifold embedded in a higher dimensional ambient space helps to build intuition, for practical purposes we adopt a complementary view of the low dimensional manifold inspired by~\citep{frohlich2021bayesian}. According to this view, we consider the entire observation space $\mathbb{R}^d$ as a smooth Riemannian manifold, $\mathcal{M}\dot{=}\mathbb{R}^d$, characterised by a Riemannian metric $\boldsymbol{\mathfrak{h}}$. The effect of the nonlinear geometry of the observations is then captured by the metric $\boldsymbol{\mathfrak{h}}$. Thus to approximate the geometric structure of the system's invariant density, we learn the Riemannian metric tensor $H:\mathbb{R}^d \rightarrow \mathbb{R}^{d \times d} $ and compute the geodesics between consecutive observations according to the learned metric. Intuitively, according to this view the observations $\{\boldsymbol{\mathcal{O}}_k\}^K_{k=1}$ introduce distortions in the way we compute distances in the state space. The advantage of this approach is that we do not need to estimate the dimensionality of the empirical manifold, which would have been difficult due to the presence of fluctuations in the system's dynamics. Instead, we still operate in the original space, while the empirical manifold introduces distortions in the estimated metric (Fig.~\ref{fig1}\textbf{F.i.}).

  \textbf{Inference framework.}~Our approach comprises three steps: 
 \textbf{\textcolor{mymaroon}{($\boldsymbol{\alpha}$.)}}~Approximation of the geometric structure of the system's invariant density with metric learning, \textbf{\textcolor{mymaroon}{($\boldsymbol{\beta}$.)}}~estimation of the (latent) system state between consecutive observations guided by the invariant density (\textbf{path augmentation}), and \textbf{\textcolor{mymaroon}{($\boldsymbol{\gamma}$.)}}~data-driven estimation of the drift function (Fig.~\ref{fig1}). We perform the two final steps in an iterative manner within an Expectation Maximisation (\textbf{E.M.}) framework~\citep{dempster1977maximum}.
 \vspace{-8.5pt}
\paragraph{\textbf{\textcolor{mymaroon}{($\boldsymbol{\alpha}$.)}}~Approximating the Riemannian geometry induced by the observations.}

Although there are many methods for approximating Riemannian manifolds~\citep{tenenbaum2000global, balasubramanian2002isomap,mead1992review,roweis2000nonlinear}, our objective is to obtain a representation that acts as a \textit{local} constraint for subsequent state estimation between consecutive observations. We achieve this in two steps:
\textbf{\textcolor{cadetblue}{(i.)}}~We approximate in the ambient space $\mathbb{R}^d$ the metric $\boldsymbol{\mathfrak{h}}$ induced by the observations (see Fig.~\ref{fig1}\textbf{F.i.}). This identifies regions of the state space with high observation density (represented with small metric values). 
\textbf{\textcolor{cadetblue}{(ii.)}}~We construct geodesics between consecutive observations on the {empirical manifold} $(\mathcal{M}\dot{=} \mathbb{R}^d, \boldsymbol{\mathfrak{h}})$ (see Fig.~\ref{fig1}\textbf{F.ii.}). The geodesics provide geometry-informed reference curves between consecutive observations that identify high-density paths through the empirical geometry. Each such path subsequently functions as a soft proximity constraint during the latent state estimation.
 \vspace{-6.5pt}
\paragraph{\textbf{\textcolor{cadetblue}{(i.)}}~Approximation of the invariant metric.}
To approximate the (local) metric $\boldsymbol{\mathfrak{h}}$ in a nonparametric form at locations $\mathbf{x}$ of the state space, we follow~\cite{arvanitidis2019fast}, and consider the inverse of the weighted local diagonal covariance computed on the $K$ observations as  \vspace{-6.5pt}
\begin{equation}\label{eq:metric_approx}
    H_{dd}(\mathbf{x}) = \left(  \sum\limits^K_{k=1} w_k(\mathbf{x}) \left( \mathcal{O}^{(d)}_k - x^{(d)}\right)^2 + \epsilon   \right)^{-1},
\end{equation}\vspace{-2.5pt}
with weights $w_k(\mathbf{x}) = \exp \left(- \frac{\|  \boldsymbol{\mathcal{O}}_k - \mathbf{x} \|^2_2}{2 \sigma^2_{\mathcal{M}}}  \right)$, and $A^{(d)}$ denoting the $d$-th dimensional component of the vector $\mathbf{A}$ for $\mathbf{A}\in\{\mathbf{x}, \boldsymbol{\mathcal{O}}_k\}$. The parameter $\epsilon > 0$ is a small value ensuring non-zero diagonals of the weighted covariance matrix, while $\sigma_{\mathcal{M}}$ is a hyper-parameter characterising the curvature of the approximated manifold.
\vspace{-10.5pt}
\paragraph{\textbf{\textcolor{cadetblue}{(ii.)}}~Constructing geodesics between consecutive observations.} To compute the geodesic curves connecting consecutive observations on the empirical manifold, we employ the approximated metric tensor $\mathbf{H}(\mathbf{x})$. We identify the geodesic curve $\boldsymbol{\gamma}^k_{t'}$ between $\boldsymbol{\mathcal{O}}_k$ and $\boldsymbol{\mathcal{O}}_{k+1}$ as the curve with minimum energy that connects these two points, i.e., as the minimiser of the kinetic energy functional ${\mathcal{E}(\boldsymbol{\gamma}^k_{t'}) =\int^1_0 L_{{\mathcal{M}}}(\boldsymbol{\gamma}^k_{t'}, \dot{\boldsymbol{\gamma}}^k_{t'})\, \text{d}t'}$
\begin{align} \label{eq:geodesic}
  \boldsymbol{\gamma}^{k*}_{t'} &=  \underset{  \begin{array}{c} \boldsymbol{\gamma}^k_{t'},\\ \boldsymbol{\gamma}^k_0 = \boldsymbol{\mathcal{O}}_k, \boldsymbol{\gamma}^k_1=\boldsymbol{\mathcal{O}}_{k+1} \end{array}  }{\arg\min} \int^1_0 L_{{\mathcal{M}}}(\boldsymbol{\gamma}^k_{t'}, \dot{\boldsymbol{\gamma}}^k_{t'}) \text{d}t', \\\quad 
\text{with}& \;\;\;\;  \int^1_0 L_{{\mathcal{M}}}(\boldsymbol{\gamma}^k_{t'}, \dot{\boldsymbol{\gamma}}^k_{t'})\, \text{d}t'= \frac{1}{2}  \int^1_0 \|\dot{\boldsymbol{\gamma}}^k_{t'} \|^2_{\mathfrak{h}} \text{d}t'\nonumber  ,
\end{align} 
where $L_{{\mathcal{M}}}(\boldsymbol{\gamma}^k_{t'}, \dot{\boldsymbol{\gamma}}^k_{t'})$ is an appropriately constructed Lagrangian.
The minimising curve of this functional is the same as the minimiser of the curve length functional $\ell(\boldsymbol{\gamma}_{t'})$ (c.f. Eq.~\ref{eq:ell}), i.e., the geodesic~\citep{do1992riemannian}.
This results in a system of second order differential equations (Eq.~\ref{eq:geode})~\citep{arvanitidis2017latent,do1992riemannian} (Sec.~\ref{appsec:with})
with boundary conditions $\boldsymbol{\gamma}^k_0 = \boldsymbol{\mathcal{O}}_k $ and $ \boldsymbol{\gamma}^k_1=\boldsymbol{\mathcal{O}}_{k+1}$ that we solve with a probabilistic differential equation solver as in~\citep{arvanitidis2019fast}.

\begin{figure}[t]
 
  \begin{center}
  \begin{overpic}[width=1.05\linewidth]{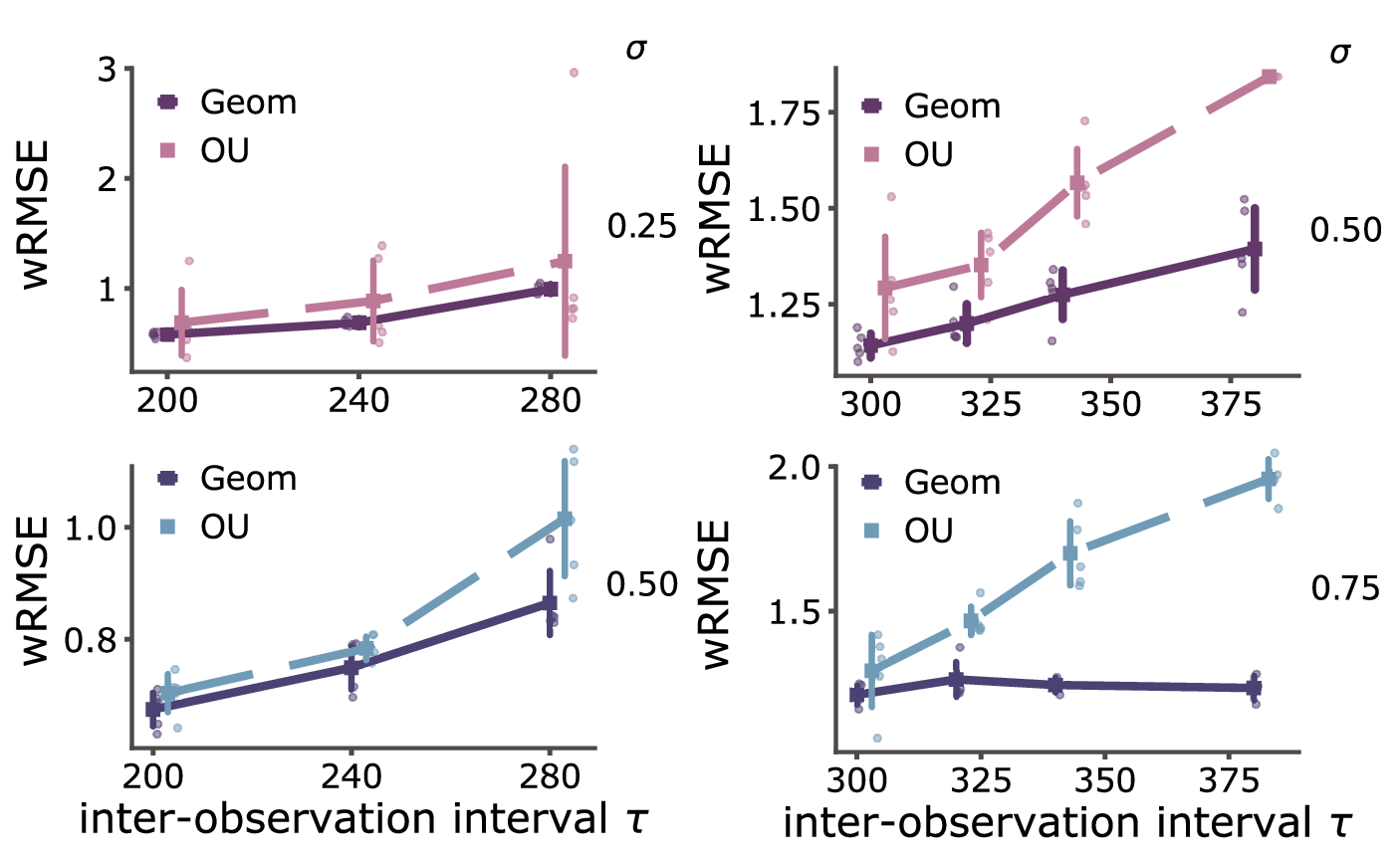}
  \put(0,55){a.}
  \put(0,30){b.}  
  \put(50,55){c.}
  \put(50,30){d.}
  \put(44,5.5){\fontsize{6pt}{7pt}{$\times \de t$}}
  \put(94,5.5){\fontsize{6pt}{7pt}{$\times \de t$}}

  \put(44,32.5){\fontsize{6pt}{7pt}{$\times \de t$}}
  \put(94,32.5){\fontsize{6pt}{7pt}{$\times \de t$}}
  \end{overpic}
  \end{center}
  \captionsetup{width=0.52\textwidth}
  \vspace{-10pt}
  \caption{\textbf{Comparison of geometry-aware inference against inference with Ornstein-Uhlenbeck augmentation.} 
Weighted root mean square error (wRMSE) vs. different inter-observation intervals $\tau$ for different noise amplitudes for moderate inter-observation intervals with \textbf{a.)} $\sigma=0.25$ and \textbf{b.)} $\sigma=0.50$, and for large inter-observation intervals with \textbf{c.)} $\sigma=0.50$ and \textbf{d.)} $\sigma=0.75$, where only one observation per oscillation period is available. 
Error bars indicate one standard deviation over five independent runs. ($\de t = 0.01$)}
  \label{fig:ou} 
 \vspace{-10pt} 
\end{figure}

\begin{table*}[t]
\rad{1.3}
\setlength{\tabcolsep}{6pt}
\rowcolors{2}{white}{mygreen}
\resizebox{\textwidth}{!}{
\begin{tabular}{rrrrcrrrcrrr}\toprule
& \multicolumn{6}{c}{\textbf{Van der Pol}} & \phantom{abc}&\\
\cmidrule{2-5} \cmidrule{6-8} \cmidrule{10-11}
\textbf{wRMSE} $\downarrow$& \makecell{total duration \\T} &  $\tau=80\times \de t$ & $\tau=120\times \de t$ & $\tau=160\times \de t$  &  $\tau=200\times \de t$ & $\tau=240\times \de t$ & $\tau=280$ $\times \de t$\\ \midrule
$\boldsymbol{\sigma}\mathbf{=0.25}$ \\
GP&500 & 0.642 $\pm$ 0.006&0.879 $\pm$ 0.005&1.083 $\pm$ 0.015& 1.258 $\pm$ 0.011&1.399 $\pm$ 0.003& 1.528 $\pm$ 0.0153&\\
SVISE&500 & 1.465 $\pm$ 0.009& 0.857 $\pm$ 0.021&    0.740  $\pm$ 0.072&    0.592 $\pm$ 0.026& $\mathbf{0.587 \pm 0.112}$&   $\mathbf{0.824 \pm  0.003}$&\\
KM-basis&500 & $\mathbf{0.368 \pm 0.054}$ & 0.452 $\pm$ 0.011& 0.671 $\pm$ 0.023 &1.588 $\pm$ 0.021&1.751 $\pm$ 0.008 & 1.735 $\pm$ 0.020&\\
LatentSDE&500 & 1.091 $\pm$0.316&1.091 $\pm$ 0.039 & 1.098 $\pm$ 0.023 & 1.089 $\pm$ 0.036       &1.088 $\pm$ 0.038&  1.091 $\pm$ 0.039 &\\
LatentSDE+GP-pre&500 & 1.095 $\pm$ 0.038 &1.085 $\pm$ 0.039 & 1.101 $\pm$ 0.034 & 1.089 $\pm$ 0.038       &1.106 $\pm$ 0.045&  1.102 $\pm$ 0.039 &\\
{GSBM} &500  & 1.518 $\pm$ 0.033& 1.435 $\pm$ 0.055& -& -& -& -&  \\  
{[SF]2M} &1500  & 1.741 $\pm$ 0.304& 1.801 $\pm$ 0.226& 1.745 $\pm$ 0.322& 1.583 $\pm$ 0.132& 1.816 $\pm$ 0.228& 1.721 $\pm$ 0.094 &  \\    
{MFM}$_{\text{LAND}}$&1500 &1.463 $\pm$ 0.007 & 1.469 $\pm$ 0.005 & 1.469 $\pm$ 0.012  &1.469 $\pm$ 0.008 & 1.469 $\pm$ 0.006 & 1.467 $\pm$ 0.008& \\
\textbf{Geometric (our)}&500 & 0.474 $\pm$ 0.034& $\mathbf{0.413 \pm 0.016}$& $\mathbf{0.514 \pm 0.068}$& $\mathbf{0.578 \pm 0.022 }$&0.687 $\pm$ 0.032& 0.993 $\pm$ 0.037&\\
\cmidrule{1-8}
$\boldsymbol{\sigma}\mathbf{=0.50}$\\
GP&500 & 0.691 $\pm$ 0.029 & 0.916 $\pm$ 0.014 & 1.114 $\pm$ 0.15& 1.272 $\pm$ 0.030& 1.409 $\pm$ 0.019&1.542 $\pm$ 0.044&\\
SVISE&500 & 1.235 $\pm$ 0.083& 0.9935$\pm$ 0.015 & 0.7505 $\pm$ 0.052 &0.736 $\pm$ 0.072 &1.3565 $\pm$ 0.278 &1.425 $\pm$ 0.086 
&\\
KM-basis&500 & 0.495 $\pm$ 0.010&0.727 $\pm$ 0.008& 0.890 $\pm$ 0.024&1.683 $\pm$ 0.020 & 1.744 $\pm$ 0.038&1.732 $\pm$ 0.065&\\
LatentSDE&500 & 1.158 $\pm$ 0.036 & 1.151 $\pm$ 0.045 & 1.160 $\pm$ 0.032 & 1.151$\pm$ 0.036 & 1.146 $\pm$ 0.033 &
               1.176 $\pm$ 0.046&\\
LatentSDE+GP-pre&500 & 1.158 $\pm$ 0.045 & 1.159 $\pm$ 0.034 & 1.159 $\pm$ 0.027 & 1.151$\pm$ 0.034 & 1.150 $\pm$ 0.028 &
               1.191 $\pm$ 0.052&\\
{GSBM} &500  & 6.106 $\pm$ 2.988& 4.818 $\pm$ 3.060& 4.738 $\pm$ 3.304& 4.875 $\pm$ 3.222& 9.076 $\pm$ 1.451& \textcolor{gray}{26.187 $\pm$ 18.804}& \\    
{[SF]2M} &1500  & 1.869 $\pm$ 0.482& 1.813 $\pm$ 0.286& 1.484 $\pm$ 0.096& 1.876 $\pm$ 0.247& 1.753 $\pm$ 0.158& 1.707 $\pm$ 0.233  & \\             
{MFM}$_{\text{LAND}}$&1500 & 1.517 $\pm$ 0.011 & 1.526 $\pm$ 0.006 & 1.536 $\pm$ 0.009 & 1.537 $\pm$ 0.017 & 1.528 $\pm$ 0.015 & 1.545 $\pm$ 0.019&\\        
\textbf{Geometric (our)}&500 &$\mathbf{0.462 \pm 0.019}$&$\mathbf{ 0.541 \pm 0.023}$& $\mathbf{0.621 \pm 0.012} $& $\mathbf{0.675 \pm 0.030}$& $\mathbf{0.750 \pm 0.038}$& $\mathbf{0.865 \pm 0.057}$&\\
\bottomrule
\end{tabular}}
\captionsetup{width=\textwidth}
\caption{Performance comparison in terms of weighted root mean square error (wRMSE) of considered frameworks for different noise conditions $\sigma$ and inter-observation intervals $\tau$ for the Van der Pol system. For all settings $\de t = 0.01$. }
\label{tab:LC}
\end{table*}

\vspace{-6pt}
\paragraph{\textbf{\textcolor{mymaroon}{($\boldsymbol{\beta}$.)}}~Latent state estimation: Geometry-guided augmentation.}
To estimate the unobserved system state between consecutive observations $\boldsymbol{\mathcal{O}}_k$ and $\boldsymbol{\mathcal{O}}_{k+1}$, we perform variational inference~\citep{beal2003variational}(see Sec.~\ref{approx_over_paths}).
Given a prior diffusion process with drift $\hat{\mathbf{f}}(\cdot): \mathbb{R}^{d }\rightarrow \mathbb{R}^{d }$ and diffusion $\sigma$, we construct an \textbf{approximating process} conditioned \textbf{i.)}~to pass through the observations, and \textbf{ii.)}~to respect the local geometry of the invariant density as it is represented by the geodesics. The conditioned process is also a diffusion process with the same diffusion constant and an effective drift function $\mathbf{g}(\mathbf{x},t)$~\citep{chetrite2015variational, majumdar2015effective}. The path probability measure ${Q}_X(\mathbf{X}_{0:T})$ induced by the approximating process
\begin{align}\label{eq:controlled_SDE}
{Q}_X(\mathbf{X}_{0:T}) :\,\,\,\, \text{d}\mathbf{X}_t &= \mathbf{g}\left(\mathbf{X}_t,t\right) \text{d}t + \sigma \text{d}\bar{\mathbf{W}}_t \\&= \left(\widehat{\mathbf{f}}(\mathbf{X}_t) +  \mathbf{u}(\mathbf{X}_t,t) \right) \text{d}t +  \sigma \text{d}\bar{\mathbf{W}}_t,  \nonumber
\end{align}
provides an approximation to the unobserved continuous system state. In Eq.~\ref{eq:controlled_SDE}, ${\mathbf{u}(\cdot,\cdot):\mathbb{R}^{d } \times \mathbb{R}^+ \rightarrow \mathbb{R}^{d }}$  is a time-dependent control term that guides the approximating path distribution through the observations, while staying in the vicinity of the corresponding geodesics between them.

More precisely, we obtain the controlled drift $\mathbf{g}\left(\mathbf{X}_t,t\right)$ from the solution of the variational problem of minimising the functional (see Sec.~\ref{appsec:without})
\begin{align}\label{eq:main_free_energy}
    \mathcal{F}[Q_X] &= \mathcal{KL}\Big(Q_X(\mathbf{X}_{0:T})||\PP(\mathbf{X}_{0:T}\mid \hat{\mathbf{f}}) \Big) \\&\;\;\;\;- \sum \limits_{k=1}^K \Big\langle\ln \PP(\boldsymbol{\mathcal{O}}_{k}\mid \mathbf{X}_{t_k})\Big\rangle_{{Q}} + \Big\langle\| \boldsymbol{\Gamma}_t -\mathbf{X}_{0:T} \|^2 \Big\rangle_{{Q}} \nonumber\\
  &= \scriptstyle \frac{1}{2} \int \limits^T_0 \int \Big[ \|\mathbf{g}(\mathbf{x},t) - \hat{\mathbf{f}}(\mathbf{x})\|_{\mathbf{D}}^2 + U_{\mathcal{O}}(\mathbf{x},t) + \beta U_{\mathcal{G}}(\mathbf{x},t)  \Big]  q_t(\mathbf{x}) \text{d}\mathbf{x} \text{d}t,\nonumber
\end{align}
where $\boldsymbol{\Gamma}_t$ denotes the sequence of $K$ geodesics indexed by time $t$, $\boldsymbol{\Gamma}_t\dot{=} \{\boldsymbol{\gamma}^k_{t'}\}_{t=(k-1)\tau+t' \tau}$, where $\boldsymbol{\gamma}^k_{t'}$ is the geodesic connecting $\boldsymbol{\mathcal{O}}_k$ and $\boldsymbol{\mathcal{O}}_{k+1}$, and $t'\in [0,1]$ denotes a rescaled time variable, and $\beta$ is a weighting term. In Eq.~\ref{eq:main_free_energy}, the term ${U_{\mathcal{O}}(\mathbf{x},t) = - \sum \limits_{t_k} \ln \PP(\boldsymbol{\mathcal{O}}_k \mid \mathbf{x})\, \delta(t- t_k)}$ \textbf{forces the augmentation to pass through the observations at each bridge boundary}, while ${U_{\mathcal{G}}(\mathbf{x},t) \dot{=} \| \boldsymbol{\Gamma}_t -\mathbf{x} \|^2 }$ \textbf{guides the latent path towards the identified geodesics}. This geodesic term should be considered as a soft penalty rather than a hard constraint: the controlled process is allowed to deviate from ($\boldsymbol{\Gamma}_t$), with the strength of the geometric bias determined by the scalar weight $\beta$ (Fig.~\ref{fig:beta}).

This minimisation can be construed as a stochastic control problem~\citep{opper2019variational} with the objective to identify a time-dependent drift adjustment $\mathbf{u}(\mathbf{x},t):=\mathbf{g}(\mathbf{x},t) - \hat{\mathbf{f}}(\mathbf{x})$ for the system with drift $\hat{\mathbf{f}}(\mathbf{x})$ so that the controlled dynamics fulfil the path constraints $U_{\mathcal{O}}(\mathbf{x},t)$ and $U_{\mathcal{G}}(\mathbf{x},t)$. 
 The first term,
$\|g(x,t)-\hat f(x)\|^2_D$, is the quadratic control cost: it penalises changes of the drift away from the current estimate $\hat f$. The potential $U_O(x,t)$ enforces consistency with the observed endpoints, while $U_G(x,t)$ biases the bridge toward the geodesic.

The optimal time-dependent control for the interval between $\boldsymbol{\mathcal{O}}_k$ and $\boldsymbol{\mathcal{O}}_{k+1}$ results from the solution of the backward equation~\citep{kappen2005path,maoutsa2022deterministic}
\begin{equation} \label{main:phipde}
    \frac{\partial \phi_t(\mathbf{\mathbf{x}})}{\partial t} = - \mathcal{L}_{\hat{\mathbf{f}}}^{\dagger} \phi_t(\mathbf{x}) + \beta \,U_{\mathcal{G}}(\mathbf{x},t) \phi_t(\mathbf{x}),
\end{equation}
with terminal condition $\phi_{t_{k+1}}(\mathbf{x}) = \chi(\mathbf{x}) = \delta(\mathbf{x}-\boldsymbol{\mathcal{O}}_{k+1}) $ and with $\mathcal{L}_{\hat{f}}^{\dagger}$ denoting the adjoint Fokker-Planck operator for the process of Eq.~\ref{appeq:prior_sde}.
As shown in~\cite{maoutsa2022deterministic} the optimal drift adjustment $\mathbf{u}(\mathbf{x},t)$ can be expressed in terms of the difference of the logarithmic gradients of two probability flows
\begin{equation}\label{maineq:control}
    \mathbf{u}^*(\mathbf{x},t) = \mathbf{D} \Big( \nabla \ln q_{T-t}(\mathbf{x}) - \nabla \ln \rho_t(\mathbf{x}) \Big),
\end{equation}
where $\rho_t$ fulfils the forward (filtering) partial differential equation (PDE)
 \begin{align}
\frac{\partial \rho_t(\mathbf{x})}{\partial t} = {\cal{L}}_{\hat{\mathbf{f}}} \rho_t(\mathbf{x}) - \beta \,U_{\mathcal{G}}(\mathbf{x},t) \rho_t(\mathbf{x}),
\label{eq:FPE2} 
\end{align}
while $q_t$ is the solution of a time-reversed PDE with initial condition ${q}_{0} (\mathbf{x}) \propto \rho_T(\mathbf{x}) \chi(\mathbf{x})$
\begin{align}\label{Fokker_bridge3}\scriptstyle
\frac{\partial {q}_{t}(\mathbf{x})}{\partial t} = 
\scriptstyle -\nabla\cdot \Big[\big(\mathbf{D}\;\nabla \ln  \rho_{T-t} (\mathbf{x})  - \hat{\mathbf{f}}(\mathbf{x}, T-t)\big)  {q}_{t} (\mathbf{x})\Big]  +  \frac{\mathbf{D}}{2} \nabla^2 {q}_{t} (\mathbf{x}).\nonumber
\end{align}

The density $\rho_t$ is therefore a forward filtering density: it propagates probability mass under the current drift estimate $\hat{\mathbf{f}}$, while down-weighting paths that accumulate large geometric cost. The terminal factor $\chi(\mathbf{x})$ imposes the endpoint condition at $\mathcal{\boldsymbol{O}}_{k+1}$. The density $q_t$ is then evolved backward in time from $q_0(\mathbf{x})\propto \rho_T(\mathbf{x})\chi(\mathbf{x})$, and carries the information required to steer the process toward the endpoint.

 For each interval $[\boldsymbol{\mathcal{O}}_k, \boldsymbol{\mathcal{O}}_{k+1}]$ we identify the posterior path measure (minimiser of Eq.~\ref{apeq:free_energy2}) by solving such a stochastic control problem for the time-varying control $\mathbf{u}(\mathbf{x},t)$ of Eq.~\ref{maineq:control}. This results in a set of $K-1$ \emph{independent} optimal control problems, that are solved in parallel for efficiency.

\vspace{-6pt}
\paragraph{\textbf{\textcolor{mymaroon}{($\boldsymbol{\gamma}$.)}}~Estimating the drift.}
We approximate the drift function in a model independent framework by imposing a Gaussian process prior on the function values ${\mathbf{f} \sim \PP_o({\mathbf{f}}) =\mathcal{GP}(\mathbf{m}^f, k^f)}$, where $\mathbf{m}^f$ and $k^f$ denote the mean and covariance function of the Gaussian process. 
The optimal measure for the drift ${Q}_f$ is a Gaussian process given by~\citep{batz2018approximate}
\begin{equation} \label{appeq:drift_measure}\footnotesize
    {Q}_f \propto \PP_o \exp\left({ -\frac{1}{2} \int  \|\mathbf{f}(\mathbf{x})\|_{D}^2 A(\mathbf{x}) - 2 \langle \mathbf{f}(\mathbf{x}), B(\mathbf{x}) \rangle_{D}  } \text{d}\mathbf{x}\right),
\end{equation}
with 
$\displaystyle A(\mathbf{x})\dot{=} \int^T_{0} {q}_t(\mathbf{x}) \text{d}t 
\;\text{and}\;B(\mathbf{x})\dot{=} \int^T_{0} {q}_t(\mathbf{x}) \mathbf{g}(\mathbf{x},t) \text{d}t,$ where ${q}_t(\mathbf{x})$ denotes the marginal density of the constrained process' state obtained by the state estimation. The function $\mathbf{g}(\mathbf{x},t)$ denotes the effective (time-dependent) drift of the constrained process (Eq.~\ref{eq:controlled_SDE}), resulting from the solution of the individual control problems accounting for the observations and the invariant geometry.

\vspace{-12pt} 

\section{Results}


\textbf{Revealing stochastic dynamics in model systems.}
  To demonstrate the
effectiveness of our approach, we inferred the stochastic dynamics of model
systems, and compared the resulting estimates to those obtained from: \textbf{\textcolor{mymaroon}{(i.)}} Gaussian process regression without state estimation (\textbf{GP}), \textbf{\textcolor{mymaroon}{(ii.)}} path augmentation with Ornstein-Uhlenbeck dynamics (\textbf{OU})~\citep{batz2018approximate}, \textbf{\textcolor{mymaroon}{(iii.)}} sparse variational inference with state estimation (\textbf{SVISE})~\citep{course2023state}, \textbf{\textcolor{mymaroon}{(iv.)}} basis function approximation of Kramers-Moyal coefficients, i.e. the drift function (\textbf{KM-basis})~\citep{nabeel2025discovering}, and \textbf{\textcolor{mymaroon}{(v.)}} latent SDE inference with amortized reparameterization with (\textbf{LatentSDE+GP-pre}) and without pre-training (\textbf{LatentSDE})~\citep{course2023amortized}, \textbf{\textcolor{mymaroon}{(vi.)}} metric flow matching (\textbf{MFM})~\citep{kapusniak2024metric}(with RBF~\citep{arvanitidis2021geometrically} and LAND metric~\citep{arvanitidis2019fast} metric approximations - see Sec~\ref{appsec:ablation_metric}), \textbf{\textcolor{mymaroon}{(vii.)}} generalized Schr{\"o}dinger bridge matching \textbf{(GSBM)}  ~\citep{liu2023generalized}, 
\textbf{\textcolor{mymaroon}{(viii.)}} simulation-free Schr{\"o}dinger bridges via score and flow matching (\textbf{[SF]$^2$M}) ~\citep{tong2023simulation}
(Sec.~\ref{appsec:baselines}). We tested our
method on non-conservative systems inducing diverse types of invariant geometries: \textbf{\textcolor{antiquefuchsia}{(a.)}} a Van der Pol system, \textbf{\textcolor{antiquefuchsia}{(b.)}} 
a process with harmonic trapping and circulation and a Gaussian repulsive obstacle in the centre introduced in~\cite{frishman2020learning},
\textbf{\textcolor{antiquefuchsia}{(c.)}}
a Hopf system, 
\textbf{\textcolor{antiquefuchsia}{(d.)}} a Selkov glycolysis model~\citep{sel1968self} (see Sec.~\ref{sec:details}), \textbf{\textcolor{antiquefuchsia}{(e.)}} a Repressilator system (see Fig.~\ref{fig:3D}), \textbf{\textcolor{antiquefuchsia}{(f.)}} as well as Van der Pol dynamics embedded in higher dimensional spaces (5D and 8D) (see Fig.~\ref{fig:5D} and Fig.~\ref{fig:8D}). For most settings, the proposed framework outperformed existing methods, especially for large inter-observation intervals (Table~\ref{tab:other_systems} and~\ref{tab:LC}).

\begin{figure}
 
  \vspace{-6pt} 
  \begin{overpic}[width=0.99\linewidth]{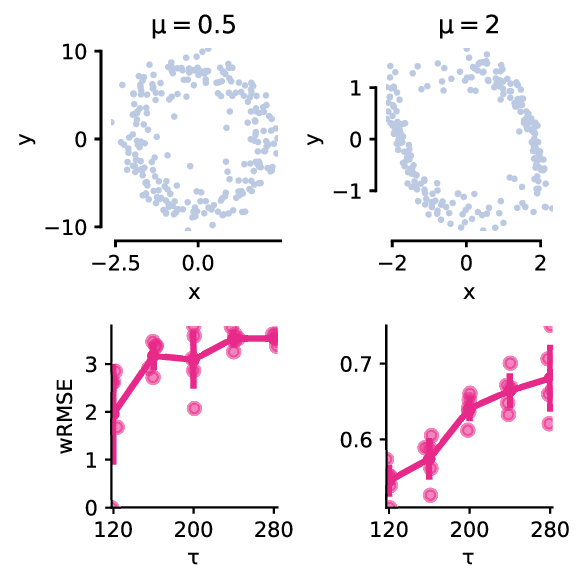}
  \put(5,90){a.}
  \put(5,50){c.}  
  \put(50,90){b.}
  \put(50,50){d.}
  \put(52,8.5){\footnotesize{$\times \de t$}}
  \put(98,8.5){\footnotesize{$\times \de t$}}
  \end{overpic}
  \captionsetup{width=0.49\textwidth}
  \caption{\textbf{Geometry-aware inference provides accurate drift estimation for different empirical manifold geometries resulting from different parameter regimes of the Van der Pol system.} \textbf{a.-b.)} Empirical manifold for the Van der Pol system with different $\mu$ parameters.
Notice the different scales on the axes. \textbf{c.-d.)} Inference performance of the proposed framework against inter-observation interval $\tau$.
Error bars indicate one standard deviation over five independent runs.  ($\de t = 0.01$)}
  \label{fig:mu} 
 \vspace{-12pt} 
\end{figure}

\begin{figure}

  \begin{overpic}[width=0.99\linewidth]{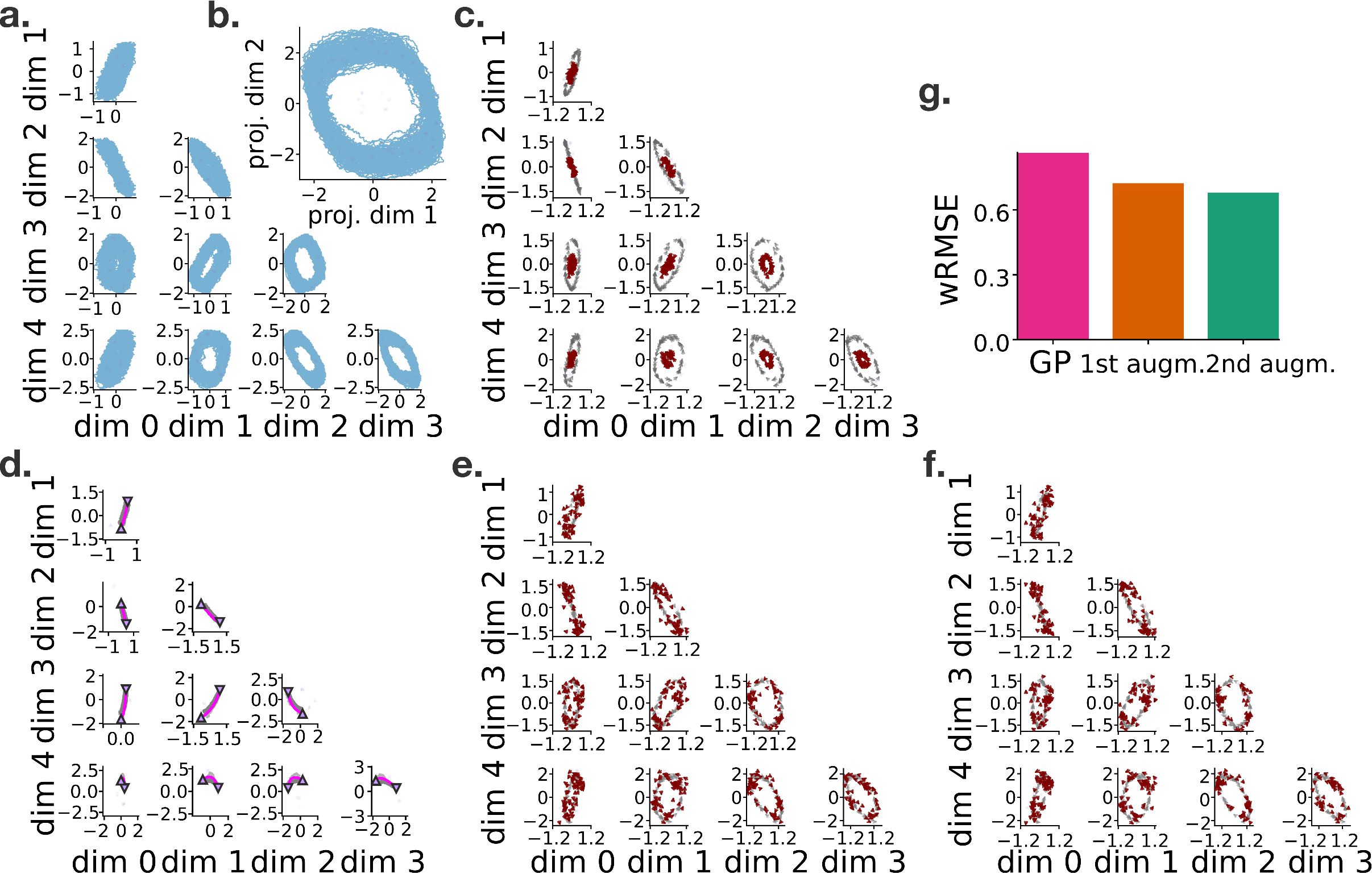}
  \end{overpic}
  \captionsetup{width=0.49\textwidth}
\caption{\textbf{Van der Pol system embedded in a 5D ambient space with isotropic 
5D dynamical noise and a mean restoring force orthogonal to the embedding subspace.} \textbf{a.)} Sparse observations (purple) together with the underlying unobserved path (light blue). \textbf{b.)} Trajectory of the Van der Pol system in its intrinsic 2D coordinates. \textbf{c.)} Observations together with the flow field estimated from sparse data by Gaussian process regression (maroon), compared to the ground-truth flow field (grey). \textbf{d.)} Example of a geometry-guided augmentation step: constructed geodesic (magenta) and sampled bridge trajectory (grey). \textbf{e., f.)} Inferred flow field after the first and second augmentation steps, respectively. \textbf{g.)} Weighted root mean squared error (wRMSE) 
of the inferred flow field at each stage of inference: initial GP estimate, 
after the first augmentation, and after the second augmentation.}
 \label{fig:5D} 
 \vspace{-12pt} 
\end{figure}

 We quantified the quality of the inference in terms of weighted root mean square error (\textbf{wRMSE}) between the estimated and ground truth drift functions evaluated on a $d-$dimensional grid spanning the state space volume of the observations. The weights for each grid point were obtained from a kernel density estimation of the observations. Thus misalignment of ground truth and estimated dynamics were penalised stronger for regions of the state space visited more frequently by the observed process.

For a system with a drift function following Van der Pol dynamics, we found that only after two E.M. iterations, the estimated force field (red arrows) was well aligned to the ground truth force field that generated the observations (grey arrows) (Fig.~\ref{fig:res}a.). For comparison we demonstrate also the result of a naive GP estimation, which resulted in a flow fields orthogonal to the ground truth ones. 
\vspace{-6pt}

We performed systematic estimations for this system under different noise conditions $\sigma$, observed at different inter-observation intervals $\tau$, for different trajectory lengths $T$ (Sec.~\ref{sec:details}). For the examined noise amplitudes (Fig.~\ref{fig:res} f.), the proposed algorithm improved the naive estimation with Gaussian assumptions within two iterations (Fig.~\ref{fig:res}). For increasing noise the improvement contributed by our method decreased (Fig.~\ref{fig:res}f.), as the invariant geometry is less well defined for large noise, but was still considerable.
\begin{table*}\centering
\small
\setlength{\tabcolsep}{4pt}
\rowcolors{2}{white}{mygreen}
\resizebox{\textwidth}{!}{
\begin{tabular}{rrrrcrrrcrrr}\toprule
& \multicolumn{3}{c}{\textbf{Harmonic trapping system}} & \phantom{abc}& 
  \multicolumn{3}{c}{\textbf{Hopf}} & \phantom{abc} & 
  \multicolumn{3}{c}{\textbf{Selkov}}\\
\cmidrule{2-4} \cmidrule{6-8} \cmidrule{10-12}
\textbf{wRMSE} $\downarrow$& $\tau=150$ & $\tau=200$ & $\tau=250$ $\times \de t$&& 
  $\tau=200$ & $\tau=300$ & $\tau=400$ $\times \de t$&& 
  $\tau=100$ & $\tau=200$ $\times \de t$&\\ \midrule
GP & 2.632 $\pm$ 0.007 & 3.387 $\pm$ 0.012& 3.733 $\pm$ 0.011 &&0.781 $\pm$ 0.006   & 0.969 $\pm$ 0.015 & 1.069 $\pm$ 0.006 && 0.550 $\pm$ 0.021 &  0.682 $\pm$ 0.040 &\\
SVISE & 35.204 $\pm$ 39.888 & 3.462 $\pm$ 0.129 &  7.540 $\pm$ 7.602 &&2.113 $\pm$ 0.658  &4.960 $\pm$ 2.687 &3.936 $\pm$ 1.063&&5.793 $\pm$ 0.028& 2.028 $\pm$ 0.045&\\
LatentSDE & \textbf{2.348} $\pm$ 0.032  &   \textbf{2.340}   $\pm$ 0.047  & \textbf{2.356} $\pm$ 0.042  && 1.168 $\pm$  0.052 &   1.161   $\pm$  0.053 &  1.173 $\pm$ 0.046 && 0.742 $\pm$ 0.022  &   0.747   $\pm$  0.021 &\\
\textbf{Geometric (ours)} &2.762 $\pm$ 0.132 & 3.034  $\pm$ 0.143 & 2.693 $\pm$ 0.992 && \textbf{0.210} $\pm$ 0.013 &\textbf{0.237} $\pm$ 0.010  &\textbf{0.255} $\pm$ 0.028 &&\textbf{0.414} $\pm$ 0.245 &  \textbf{0.682} $\pm$ 0.071 &\\
\bottomrule  
\end{tabular}
}
\captionsetup{width=\textwidth}
\caption{Performance comparison in terms of wRMSE for the considered frameworks for three different nonlinear dynamical systems and for increasing inter-observation interval $\tau$. Numbers indicate mean wRMSE and standard deviation of five independent runs for each setting. For all settings $\de t = 0.01$.}
\label{tab:other_systems}
\vspace{-10pt}
\end{table*}

\textbf{Impact of the geometry of empirical manifold.}
We performed inference for different parameter values of the Van der Pol system ($\mu=1$ (as above) and $\mu=0.5$ and $\mu=2$), that resulted in asymmetries of the invariant density (Fig.~\ref{fig:mu}). We observed that the performance of all inference frameworks deteriorated for increasing asymmetry (larger dynamic range along one dimension), yet our method still delivered more accurate predictions compared to the other considered frameworks. Moreover, employing an alternative metric-learning algorithm to approximate the invariant geometry resulted in comparable performance. Specifically, Geometric$_{\text{RBF}}$ in Table~\ref{tab:LCmetric} uses the metric of~\cite{arvanitidis2021geometrically}, further developed in~\cite{kapusniak2024metric}, in which a diagonal metric is represented as a positive linear combination of Gaussian RBFs centred at selected cluster centres.

\vspace{-5pt}
\textbf{Impact of noise amplitude.}
For systems with small dynamical noise (small $\sigma$), geodesics approximate the manifold structure better, however the path integral control is limited by the control costs proportional to inverse noise covariance. Our framework had comparable accuracy for all inter-observation lengths, but improvement was small for small lengths since in that setting the estimation with Gaussian likelihood already provided a good approximation of the ground truth drift.

 We compared our method to the approach proposed in~\cite{batz2018approximate}. In this work, the authors perform  augmentation with Ornstein-Uhlenbeck bridges, i.e. assuming linear underlying dynamics. We found that our approach delivered more accurate estimates for larger inter-observation intervals. For inter-observation intervals with only one observation per oscillation period (Fig.~\ref{fig:ou}c.,d.), our approach delivered better results by considering additionally knowledge of the direction of movement in the state space (Sec.~\ref{sec:details}). The variance of estimates of the proposed method was smaller compared to Batz et al. due to the consistency imposed by conditioning on the invariant geometry of the system. Predictions improved with longer observation intervals $T$, and for decreasing noise amplitude $\sigma$. In both settings the invariant geometry was more well approximated by the empirical manifold. 
Overall, state estimation with Ornstein-Uhlenbeck dynamics~\cite{batz2018approximate}, was increasingly less capable of correctly estimating the latent system state and subsequently correctly approximating the unknown drift function especially as the length of the inter-observation interval $\tau$ increased.

\vspace{-5pt}
  \textbf{Effects of noise miss-estimation.}We further investigated the impact of noise misestimation on the accuracy of drift inference (S.I. Fig.~\ref{fig:mis}). Our findings indicate that after two augmentations conditioned on the invariant geometry, small inaccuracies in the employed dynamical noise during the simulation of augmented paths had a negligible effect on the overall accuracy of the inferred drift. In particular, for small inter-observation intervals, our approach remained robust to misestimated noise amplitudes. As the inter-observation intervals increased, the effect of noise deviations from ground truth on the performance remained small, provided the noise used in the augmentation deviated by at most $\pm 0.1$ from the ground truth noise amplitude. Thus, stochastic dynamics may still be identified with our method even with inaccurate or moderately misestimated diffusion constants.
  
  Additional experiments with three- and eight-dimensional systems are provided in Fig.~\ref{fig:3D} and Fig.~\ref{fig:8D}. We further provide in the Supplement more ablations with respect to \textbf{(i)} the metric-learning algorithm (Sec.~\ref{appsec:ablation_metric}), \textbf{(ii)} and noise misspecification (Sec.~\ref{app:mise}).\footnote{Code is available at: \url{https://github.com/dimitra-maoutsa/From_geometry_to_dynamics}.}


\begin{figure}
 
  \vspace{-6pt} 
  \begin{overpic}[width=0.99\linewidth]{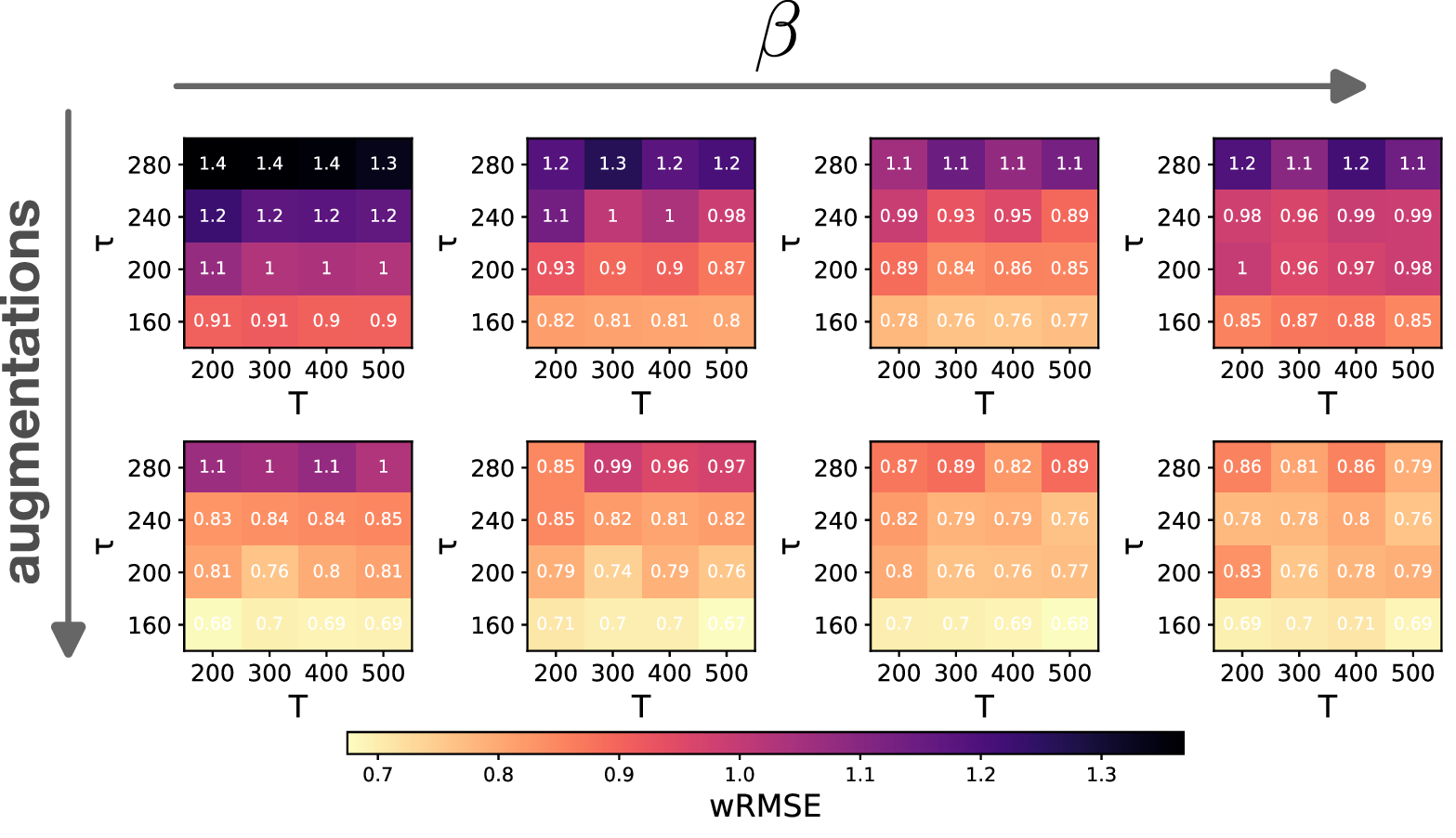}
  \end{overpic}
  \captionsetup{width=0.52\textwidth}
  \caption{\textbf{Sensitivity analysis with respect to the parameter $\beta$, which controls
 the relative weight of the geodesic control term.} Weighted root mean square error values between ground truth and estimated drift, averaged over $3$ independent runs, for different total trajectory lengths $T$ and inter-observation intervals $\tau$ with noise amplitude $\sigma=0.5$. Upper row shows results after the first iteration, while second row shows results after the second one. Values of $\beta$ from left to right $[0.1, 0.5, 1.0, 2.0]$. }
  \label{fig:beta} 
 \vspace{-10pt} 
\end{figure}

\vspace{-5pt}

\section{Discussion}\label{s:discussion}

Discovering unknown driving forces from sparse observations of stochastic systems remains still challenging: when the sampling interval is large, short-time transition approximations break down and the inverse problem of identifying the drift function becomes challenging.
Our results show that combining the \textbf{temporal ordering} of the observations while accounting for the \textbf{geometry of the invariant density} can substantially improve drift recovery in nonlinear stochastic systems in non-conservative settings with steady-state probability currents.

\textbf{Geometric inductive bias.}
We introduced \textbf{geometric inductive biases} into inference of stochastic systems by treating the deterministic flow field as a scaffold upon which system states fluctuate. We approximated this scaffold in terms of \textbf{distortions of a metric induced by the system's measurements}. This approach effectively approximates the low-dimensional invariant density (empirical manifold) without the need to project to a lower dimensional space, whose dimensionality would be hard to estimate due to the presence of fluctuations.

The key insight is that \textbf{geodesics} computed on the empirical manifold with respect to the approximated metric 
provide geometry-informed reference curves for latent path augmentation. These curves are motivated by the \textbf{most probable paths} of the unknown system between consecutive observations in the Onsager-Machlup sense, but are not explicitly assumed to be the true unobserved paths. They rather serve as soft constraints induced by the approximated invariant geometry.
Using these \textbf{geodesics as control constraints}, we formulated a path-augmentation scheme that bridges sparse observations with trajectories consistent with both the temporal order and the geometry of the data.

\textbf{Why standard path augmentation fails at large $\tau$.}
Widely used inference methods, predominantly developed within the statistics community, often employ path (\emph{data}) augmentation to approximate transition densities between successive observations. However, this approach suffers from several challenges: \textbf{1.)} The unobserved information between successive observations is an infinite-dimensional object, requiring the solution of a complex and computationally intensive problem (bridge sampling)~\citep{beskos2006exact}. We addressed this challenging problem using the computationally efficient framework developed in~\cite{maoutsa2022deterministic}. \textbf{2.)} Direct drift estimation from sparse observations results in estimated dynamics that significantly deviate from the ground truth. Thereby consecutive observations of the system have small probability under the law of the estimated SDE. This discrepancy, in turn, leads to several computational difficulties: \textbf{i)} Most bridge sampling schemes become computationally demanding, or fail, when attempting to generate transition densities between atypical states for the considered stochastic dynamics. For instance, the method of~\citep{maoutsa2021deterministica} successfully generates transition densities between atypical states only for conservative systems through a reweighting with Brownian bridge dynamics. Alternatively, an exceedingly large number of samples would be required for accurate numerical approximation.
\textbf{ii)} Iterative algorithms, such as Expectation Maximisation, which exhibit only {local} convergence~\citep{romero2019convergence}, may converge to inaccurate solutions, when the initial estimation significantly deviates from the ground truth.

\textbf{Relation to Schr\"odinger bridges and controlled transport.}
The framework we employed for the augmentation relies on a deterministic particle formulation of the path integral control formalism~\citep{kappen2005linear}. This framework can be connected to the dynamic Schr\"odinger bridge problem, if we consider transferring probability mass between two Dirac measures or very narrow Gaussians that sit on each observation, considering additionally a potential that constraints the intermittent dynamics similar to the one considered in~\cite{neklyudov2023computational}. Thus, in principle, one can employ one of the recently developed alternative frameworks that solve the dynamic Schrödinger bridge problem for path augmentation. 
Recent bridge- and flow-matching methods~\citep{lipman2022flow, albergo2023stochastic, shi2023diffusion, liu2023generalized} solve related controlled-transport problems, but typically without the path constraint we employed here.
 In contrast, the Generalised Schrödinger Bridge Matching (GSBM)
framework~\citep{liu2023generalized} uses a cost functional that is equivalent to the controlled cost we employ to construct our augmentations (SI Eq.~\ref{appeq:free_energy2}). In this setting, the penalty term
corresponds to the geodesic proximity constraint used in our framework. The GSBM could, in
principle, replace the particle-based framework we use.
We nevertheless adopted a particle-based formulation, because the resulting particles provide direct samples from the controlled path measures, and can therefore be reused to construct Monte Carlo approximations of the integrals entering the Gaussian-process drift inference (Eq.~\ref{appeq:drift_measure}).
 Yet, the Gaussian variant of the GSBM framework that incorporates time-dependent penalty constraints (analogous to our geodesic constraints), might be an interesting avenue to explore for potential incorporation in our framework.

\textbf{Limitations and outlook.}
The proposed approach relies on the geometric characterisation of the invariant density of the system's dynamics. This requires sufficiently long observation windows to accurately characterise said density and correctly approximate the unobserved paths with geodesic curves. Thus, our approach is limited to systems where the invariant density can be approximated by a manifold where we can identify geodesics. An alternative method worth exploring would consider the learned invariant metric directly in the dynamics of the augmented process. Our evaluation was restricted to synthetic systems with known ground-truth dynamics. While this allows us to isolate the effect of sparse sampling and directly quantify drift recovery, validating the method on real datasets with naturally sparse observations remains an important direction for future work.
 Moreover, we have considered here inference of stochastic differential equations with known state independent diffusion. While this approach might seem limited, several processes with state dependent diffusions can be transformed into processes with state independent diffusions~\citep{beskos2006exact, roberts2001inference} through the Lamperti transform if they fulfil the appropriate conditions for the drift function.

\newpage



\newpage

\section*{Acknowledgements}

The author is indebted to Manfred Opper for prompting them to work on this problem,  Nina Miolane for early guidance on extracting geodesics from variational autoencoder based approximation of data manifolds, Stefan Sommer for answering questions on diffusions on manifolds, and to Georgios Arvanitidis for maintaining an open repository with his previous work on approximating geodesics from data manifolds. The author gratefully acknowledges the anonymous reviewers for their careful reading of the manuscript and for their constructive feedback, which contributed to improving the final version of this manuscript. An earlier account of this work has been presented at the NeurIPS 2022 workshop \emph{Machine Learning and the Physical Sciences}~\citep{maoutsa2023geometric}, in ICLR 2023 workshop \emph{Physics for Machine Learning}~\citep{maoutsa2023geometricICLR}, and in ICLR 2026 workshop \emph{Geometry-grounded Representation Learning and Generative Modeling}~\citep{maoutsa2026flow}. The author further acknowledges that previous work from the Python~\citep{van1995python}, numpy~\citep{harris2020array}, scipy~\citep{2020SciPy-NMeth}, matplotlib~\citep{Hunter:2007}, seaborn~\citep{waskom2021seaborn}, GPflow~\citep{GPflow2017}, pyEMD~\citep{pele2009}, and pytorch~\citep{paszke2017automatic} communities facilitated the implementation of the computational part of this work.

\section*{Impact Statement}

The aim of this work is to advance the field of dynamical inference for stochastic systems. While we do not foresee any direct societal consequences directly impinging from our work, we recognize that stochastic systems could be applied in military contexts, financial engineering, or more recently in machine learning for data (such as image, audio, video) generation. Still, the proposed method does not propose interventions that might lead to unfavourable societal outcomes. Overdamped Langevin systems are widespread in areas such as physics, biology, neuroscience, and ecology. We anticipate that our contributions will thus help these disciplines by offering a tool to identify and further study relevant systems.

Our contribution emphasises the importance of incorporating concepts from the evolving field of geometric statistics into system identification methods for stochastic systems. Although geometric and topological properties of invariant densities have been extensively studied in the context of deterministic systems, comparable attention is lacking for their stochastic counterparts.
Our work further highlights that in settings where the amount of augmented data exceeds the number of observations, data augmentation frameworks can enhance inference accuracy by incorporating domain knowledge or other relevant information, such as the geometry of the system’s invariant density we consider here. Many algorithms used for data augmentation, including the expectation maximisation algorithm employed in our work~\citep{romero2019convergence}, show only \textbf{local convergence}. As a result, when the initial estimate deviates significantly from the ground truth value, naive data augmentation methods may converge to suboptimal solutions that fail to accurately identify the underlying system.

\bibliography{references.bib}
\bibliographystyle{unsrtnat}

\newpage

\onecolumn
\appendix
\section*{Supplementary Information}

\startcontents[sections]
\printcontents[sections]{l}{1}{\setcounter{tocdepth}{4}}


\section{Drift inference for high and low frequency observations}\label{appsec:b}

Effective dynamics of systems with many degrees of freedom or inherently stochastic are often described in terms of a stochastic differential equation (SDE)
\begin{equation} \label{app_eq:system}
    \text{d}\mathbf{X}_t = {\mathbf{f}(\mathbf{X}_t)} \text{d}t + \boldsymbol{\eta}(t) \text{d}t = {\mathbf{f}(\mathbf{X}_t)} \text{d}t + \boldsymbol{{\sigma}} \text{d}\textbf{W}_t,
\end{equation}
where the drift $\mathbf{f}(\cdot): \mathbb{R}^d \rightarrow \mathbb{R}^d$ describes the deterministic forces acting on the system, while the delta-correlated Gaussian white noise term $\boldsymbol{\eta}(t)$, $\langle \boldsymbol{\eta}(t) \boldsymbol{\eta}(t')  \rangle = \boldsymbol{\sigma} \boldsymbol{\sigma}^\top \delta(t-t')$ describes the effect of stochastic forces as a product of a
diffusion matrix (or constant) ${\boldsymbol{\sigma}:  \mathbb{R}^{d \times d}}$ that accounts for the magnitude of the stochastic forces acting on the system, and a $d$-dimensional Wiener process $\mathbf{W}_t$ that contributes random influences.

Often the detailed equation that governs the evolution of the state of the system is unknown. Therefore, understanding a system of interest often requires identification from time series observations of its state. In more practical terms, given some \textbf{prior probability} for the drift function, we want to compute the \textbf{posterior probability} $\PP(\mathbf{f} | \{\boldsymbol{\mathcal{O}}_k\}^K_{k=1})$ that identifies the unknown drift function of Eq.~\ref{app_eq:system} that most likely gave rise to the observations of the system state $\{\boldsymbol{\mathcal{O}}_k\}^K_{k=1}$. The exact relationship between the observations and the system state will be defined more precisely in the following.

When a system is observed nearly continuously (inter-observation interval length $\tau$ much smaller than the characteristic time scale of the system  ${\tau \ll \tau_{\text{char}}}$), temporal methods regress the system state $\mathbf{X}_t$ against the state increments $\mathbf{Y}_t \dot{=} \frac{\mathbf{X}_{t+\tau}-\mathbf{X}_t}{\tau}$ to identify the drift function~\citep{friedrich1997description, ragwitz2001indispensable}. In a Bayesian framework, this corresponds to Gaussian process regression with a Gaussian likelihood (SI~\ref{app:cont_time}). However, for large inter-observation intervals $\tau$, these methods fail~\citep{batz2018approximate}, as the Gaussian likelihood assumption is invalid for general nonlinear systems with sparse observations (Fig.\ref{fig1}\textbf{C.}). In such cases, the likelihood is a \emph{path integral} over continuous trajectories of the unobserved process (SI~\ref{app:low_freq}), making Gaussian-based estimates inaccurate (Fig.~\ref{fig1}\textbf{C.}).

This underwhelming performance has motivated the development of methods that combine state estimation (or \textbf{path augmentation}) and dynamical inference. These methods reconstruct continuous paths to approximate transition densities between observations, enabling inference by estimating the system’s state between observations. However, for large time intervals, transition densities are usually analytically intractable, except in a few trivial cases of scalar or linear processes. As a result, the prevailing strategy is to approximate transition densities by sampling marginal distributions of \textbf{diffusion bridges}, which are diffusion processes constrained by their initial and terminal states~\citep{golightly2008bayesian,papaspiliopoulos2012nonparametric,sermaidis2013markov,beskos2006retrospective,chib2006likelihood}. Yet, existing methods employ path augmentation with simplified bridge dynamics (e.g., Brownian~\citep{chib2006likelihood,golightly2008bayesian} or Ornstein-Uhlenbeck bridges~\citep{batz2018approximate}) that do not accurately reflect the underlying transition densities for nonlinear systems (Fig.~\ref{fig1}\textbf{E.}).

An alternative path augmentation strategy would obtain a coarse drift estimate, typically achieved by assuming a Gaussian likelihood between observations (see SI Eq.~\ref{apeq:SDE_likelihood}), and would subsequently employ a stochastic bridge sampler~\citep{de2021simulating,maoutsa2022deterministic,maoutsa2021deterministica} to construct stochastic bridges using the coarsely estimated nonlinear drift. However, for large inter-observation intervals, the coarsely estimated drift function often deviates significantly from the true function that generated the observations. Consequently, the observations frequently fall into low-probability regions of the estimated diffusion dynamics (Fig.~\ref{fig1} \textbf{E.}), rendering the construction of diffusion bridges either too computationally demanding or impossible~\citep{liu2020study}.

\subsection{High frequency observations} ~\label{app:cont_time}
In an optimal but rather practically unrealistic scenario, we would observe the system (path) $\mathbf{X}_{0:T}$ in (nearly) continuous time, and thus we would try to identify the drift from $\PP(\mathbf{f} | X_{0:T}
)$. In such a case, the infinitesimal transition probabilities of the diffusion process between consecutive time-points are Gaussian, i.e.,
\begin{equation} \label{appeq:continous_transitions}
    \PP_{{f}}(\mathbf{X}_{0:T} \mid \mathbf{f}) \propto \exp \left( -\frac{1}{2 \text{d}t} \sum_t  \| \mathbf{X}_{t+\text{d}t} -\mathbf{X}_t - \mathbf{f}(\mathbf{X}_t)\text{d}t\|_{D}^2     \right).
\end{equation}
Here we have introduced the weighted norm $\|\textbf{u} \|_{D}  \dot{=} \textbf{u}^{\top} \cdot \mathbf{D}^{-1} \cdot \textbf{u},$ with $\mathbf{D}\dot{=}\boldsymbol{\sigma\sigma}^\top$ indicating the noise covariance.

In turn, the transition probabilities of a discretised drift-less process (a Wiener path) $P_{\mathcal{W}}(\mathbf{X}_{0:T}) $ with same diffusion $\sigma$ is
\begin{equation}\label{appeq:wiener}
    \PP_{\mathcal{W}}(\mathbf{X}_{0:T})=\exp \left( -\frac{1}{2 \text{d}t} \sum_t  \| \mathbf{X}_{t+\text{d}t} -\mathbf{X}_t\|_{D}^2     \right).
\end{equation}

We can thus express the likelihood for the drift $f$ as the likelihood ratio between the transition probabilities of Eq.~\ref{appeq:continous_transitions} and Eq.~\ref{appeq:wiener}, which for diffusion processes is expressed by the Radon-Nykodym derivative between $\PP_f(\mathbf{X}_{0:T}|f)$ and $\PP_{\mathcal{W}}(\mathbf{X}_{0:T})$
 for paths $\mathbf{X}_{0:T}$ within the time interval $[0,\,T]$
~\citep{liptser2013statistics}
\begin{equation} \label{apeq:SDE_likelihood}
    \mathcal{L}(\mathbf{X}_{0:T} \mid \mathbf{f}) = \exp \left( -\frac{1}{2} \sum_t \| \mathbf{f}(\mathbf{X}_t)\|^2_D  \text{d}t + \sum_t \langle \mathbf{f}(\mathbf{X}_t), \mathbf{X}_{t+\text{d}t}-\mathbf{X}_t \rangle_{D}  \right),\end{equation}
where for brevity we have introduced the notation $\langle \textbf{u}, \textbf{v}\rangle_{D} \dot{=} \textbf{u}^{\top} \cdot \mathbf{D}^{-1} \cdot \textbf{v}$ for the weighted inner product with respect to the inverse noise covariance $\mathbf{D}^{-1}$. This expression results from applying the Girsanov theorem on the path measures induced by a process with drift $\mathbf{f}$ and a Wiener process, with same diffusion $\sigma$, and employing an Euler-Maruyama discretisation on the continuous path $\mathbf{X}_{0:T}$.

The likelihood of a continuously observed path of the SDE (Eq.~\ref{apeq:SDE_likelihood}) has a quadratic form in terms of the drift function. Therefore a Gaussian measure over function values (Gaussian process) is a natural conjugate prior for this likelihood. Thus, to identify the drift in a non-parametric form, we assume a Gaussian process prior for the function values ${\mathbf{f} \sim \PP_0(\mathbf{f}) =\mathcal{GP}(\mathbf{m}^f, k^f)}$, where $\mathbf{m}^f$ and $k^f$ denote the mean and covariance function of the Gaussian process ~\citep{ruttor2013approximate}. The prior measure can be written as
\begin{equation}
    \PP_0(\mathbf{f})  = \exp\left(-\frac{1}{2} \int \int \mathbf{f}(\mathbf{x}) \left(k^f(\mathbf{X},\mathbf{X}') \right)^{-1}\mathbf{f}(\mathbf{X}') \text{d}\mathbf{X} \text{d}\mathbf{X}'\right),
\end{equation}
if we consider a zero mean Gaussian process $\mathbf{m}^f=\textbf{0}$.

Bayesian inference for the drift function $\mathbf{f}$ requires the computation of a probability distribution in the function space, the posterior probability distribution $\PP_f(\mathbf{f} \mid \mathbf{X}_{0:T})$. From the Bayes' rule the posterior can be written as 
\begin{equation}
    \PP_f(\mathbf{f}\mid \mathbf{X}_{0:T}) =\frac{ \PP_0(\mathbf{f}) \mathcal{L}(\mathbf{X}_{0:T} \mid \mathbf{f})}{Z} \propto \PP_0(\mathbf{f}) \mathcal{L}(\mathbf{X}_{0:T} \mid \mathbf{f}),
\end{equation}
where $Z$ denotes a normalising factor defined as 
\begin{equation}
    Z = \int \PP_0(\mathbf{f}) \mathcal{L}(\mathbf{X}_{0:T} \mid \mathbf{f}) \mathcal{D}\mathbf{f} ,
\end{equation}
where $ \mathcal{D}\mathbf{f}$ denotes integration over the Hilbert space $\mathbf{f}: H_0[\mathbf{f}] < \infty$ . Here we have expressed the prior probability over functions as $\PP_0(\mathbf{f}) = e^{-H_0[\mathbf{f}]}$.
In~\cite{ruttor2013approximate} the authors show that in this continuous-time setting, nonparametric estimation of the drift can be attained through a  Gaussian process regression~\citep{rasmussen2003gaussian} with the objective to identify the mapping from the system state $\mathbf{X}_t$ to state increments $\text{d}\mathbf{X}_t$. More precisely, we consider as the regressor the $N$ observations of the system state $\mathbf{X}_t$ and as the associated response variables the state increments
\begin{equation}\label{eq:increments}
    \mathbf{Y}_t = \frac{\mathbf{X}_{t+\text{d}t}-\mathbf{X}_t}{\text{d}t},
\end{equation} and select the kernel function of the Gaussian process as $k^f(\mathbf{X},\mathbf{X}')$. 

If we denote with $\mathcal{X}= \{\mathbf{X}_t\}^{T-\text{d}t}_{t=0}$ and $\mathcal{Y}= \{\mathbf{Y}_t\}^{T-\text{d}t}_{t=0}$ the set of state observations and observation increments, the mean of the posterior process over drift functions $\mathbf{f}$ can be expressed as
\begin{equation}\label{eq7:full_gp_mean}
    \bar{\mathbf{f}}(\mathbf{x}) = k^f(\mathbf{x},\mathcal{X})^{\top} \left( \mathcal{K} + \frac{\mathbf{D}}{\text{d}t} I_N\right)^{-1} \mathcal{Y},
\end{equation}
where we abused the notation and denoted with $k^f(\mathbf{x}, \mathcal{X})$ the vector resulting from evaluating the kernel $k^f$ at points $\mathbf{x}$ and $\{\mathcal{O}_k\}^{K-1}_{k=1}$. Similarly $\mathcal{K} = k^f(\mathcal{X}, \mathcal{X})$ stands for the $(K-1)\times (K-1)$ matrix resulting from evaluation of the kernel on all observation pairs. 
In a similar vein, the posterior variance can be written as
\begin{equation} \label{eq7:full_gp_cov}
    \Sigma^2(\mathbf{x}) = k^f(\mathbf{x},\mathbf{x}) - k^f(\mathbf{x},\mathcal{X})^{\top} \left(\mathcal{K} + \frac{\mathbf{D}}{\text{d}t} \right)^{-1} k^f(\mathbf{x},\mathcal{X}),
\end{equation}
where the term $\mathbf{D}/\text{d}t$ plays the role of observation noise.

\subsection{Low frequency observations} ~\label{app:low_freq} 
As the inter-observation interval increases (\emph{low frequency observations}), the validity of the Gaussian likelihood used in drift estimation diminishes as the transition density is no longer Gaussian. Consequently, methods for drift estimation with Gaussian assumptions~\citep{friedrich1997description,ruttor2013approximate} become increasingly inaccurate. To discount the effects of low frequency sampling, Lade~\citep{lade2009finite} proposed a method to compute finite-time corrections for drift estimates, which has been mainly applied to one-dimensional problems~\citep{honisch2011estimation}. In parallel, the statistics community has proposed path augmentation techniques that involve sampling with a simplified system's dynamics between time-consecutive observations to augment the observed trajectory to a nearly continuous-time path~\citep{golightly2008bayesian,papaspiliopoulos2012nonparametric,sermaidis2013markov,beskos2006retrospective,chib2006likelihood}. However, for large inter-observation intervals and nonlinear systems, the augmented trajectories match poorly the underlying path statistics and these methods often exhibit poor convergence rates or fail to identify the correct dynamics (Figure~\ref{fig1} c. and d.). We note that path augmentation using Ornstein-Uhlenbeck bridges and local linearisation of the \textbf{ground truth} dynamics provides a reasonable approximation of the underlying transition density up to a certain inter-observation interval. Nevertheless, during inference, the ground truth dynamics is unknown, and the proposed local linearisations based on inaccurate drift estimates~\citep{batz2018approximate} perform poorly in this sparsely sampled regime.

As the inter-observation interval $\tau$ increases, if the system is nonlinear, the likelihood assumed between two consecutive observations is no longer Gaussian, but is rather expressed as a \emph{path integral}
\begin{equation}\label{apeq7:path_integral_likelihood}
    \PP(\boldsymbol{\mathcal{O}}_{1:K}\mid \mathbf{f}) = \int \PP(\boldsymbol{\mathcal{O}}_{1:K}\mid \mathbf{X}_{0:T}) \PP(\mathbf{X}_{0:T}\mid \mathbf{f}) \mathcal{D}(\mathbf{X}_{0:T}),
\end{equation}
where $\boldsymbol{\mathcal{O}}_{1:K}\dot{=}\{\boldsymbol{\mathcal{O}}_k\}_{k=1}^{K}$ identifies the set of $K$ observations collected within the interval $[0,\,T]$, 
$\PP(\mathbf{X}_{0:T}\mid \mathbf{f})$ the prior path probability resulting from a diffusion process with drift $\mathbf{f}(\mathbf{x})$, $\mathcal{D}(\mathbf{X}_{0:T}) $ identifies the formal volume element on the path space, and $\PP(\boldsymbol{\mathcal{O}}_{1:K}\mid \mathbf{X}_{0:T})$ stands for the likelihood of observations given the latent path $\mathbf{X}_{0:T}$.

However, the path integral of Eq.~\ref{apeq7:path_integral_likelihood} is in general intractable for nonlinear systems. thus we need to simultaneously estimate the drift and latent state of the diffusion process, i.e., to approximate the joint posterior measure of latent paths and drift functions ${\PP(\mathbf{X}_{0:T},\mathbf{f} \mid \boldsymbol{\mathcal{O}}_{1:K})}$. Therefore we consider the unobserved continuous path $\mathbf{X}_{0:T}$ as latent random variables and employ an Expectation Maximisation (E.M.) algorithm to identify a maximum a posteriori estimate for the drift function. 
More precisely, we follow an iterative algorithm, where at each iteration $n$ we alternate between the two following steps:
\vspace{5pt}
\begin{minipage}[c]{0.95\linewidth}
     An \textbf{Expectation} step, where given a drift estimate $\hat{\mathbf{f}}^n(\mathbf{x})$ we construct an approximate posterior over the latent variables ${Q(\mathbf{X}_{0:T}) \approx \PP(\mathbf{X}_{0:T}\mid \boldsymbol{\mathcal{O}}_{1:K}, \hat{\mathbf{f}}^n(\mathbf{x}))}$, and compute the expected log-likelihood of the augmented path 
    \begin{equation} \label{eq:estep}
    \mathfrak{L}\big(\hat{\mathbf{f}}^n(\mathbf{x}), Q\big) = \mathbb{E}_Q\Big[\ln \mathscr{L}\big(\mathbf{X}_{0:T},\boldsymbol{\mathcal{O}}_{1:K}\mid \hat{\mathbf{f}}^n(\mathbf{x})\big) \Big].
    \end{equation}
    \end{minipage}
    
     A \textbf{Maximisation} step, where we update the drift estimation by maximising the expected log likelihood
    \begin{equation} \label{eq:mstep}
        \mathbf{f}^{n+1}(\mathbf{x}) = \arg \max_f \Big[\mathfrak{L}\big(\mathbf{f}^n(\mathbf{x}),Q\big)-\ln \PP_0\big(\mathbf{f}^n(\mathbf{x})\big) \Big].
    \end{equation}
In Eq.~\ref{eq:mstep}, $\PP_0$ denotes the Gaussian process prior over function values.

\subsection{Approximate posterior over paths.}\label{approx_over_paths} 
To obtain an approximate posterior over the latent paths we perform \textbf{variational inference}~\citep{beal2003variational}. 
In this section, we first formulate the approximate posterior over paths (conditional distribution for the path given the observations) by considering only individual observations as constraints (Section~\ref{appsec:without}). However, this approach results computationally taxing calculations during path augmentation, since the observations are atypical states of the initially estimated drift. To overcome this issue, we subsequently extend the formalism (Section~\ref{appsec:with}) to incorporate constraints that consider also the local geometry of the observations.

\subsubsection{Approximate posterior over paths \underline{without} geometric constraints}\label{appsec:without}
Given a drift function (or a drift estimate) $\hat{\mathbf{f}}(\mathbf{x})$ we can apply variational techniques to approximate the posterior measure over the latent path conditioned on the observations $\boldsymbol{\mathcal{O}}_{1:K}$. We consider that the \textbf{prior process} (the process without considering the observations $\boldsymbol{\mathcal{O}}_{1:K}$) is described by the equation
\begin{equation} \label{appeq:prior_sde}
  \PP(\mathbf{X}_{0:T}\mid \hat{\mathbf{f}}): \qquad  \text{d}\mathbf{X}_t = \hat{\mathbf{f}}(\mathbf{X}_t) \text{d}t + \sigma \text{d}\mathbf{W}_t.
\end{equation}
We will define an approximating (posterior) process that is conditioned on the observations. The conditioned process is also a diffusion process with the same diffusion as Eq.~\ref{appeq:prior_sde} but with a modified, time-dependent drift $g(x,t)$ that accounts for the observations~\citep{chetrite2015variational,majumdar2015effective}. 
We identify the approximate posterior measure $Q$ with the posterior measure induced by an approximating process that is conditioned by the observations $\boldsymbol{\mathcal{O}}_{1:K}$~\citep{opper2019variational},
with governing equation
\begin{equation} \label{appeq:sde_q}
    Q(\mathbf{X}_{0:T}): \qquad \text{d}\mathbf{X}_t = \mathbf{g}(\mathbf{X}_t,t) \text{d}t + \sigma \,\text{d}{W}_t=\left(\hat{\mathbf{f}}(\mathbf{X}_t) +  \mathbf{u}(\mathbf{X}_t,t) \right) \text{d}t +  \sigma\, \text{d}\mathbf{W}_t.
\end{equation}

The effective drift $\mathbf{g}(\mathbf{X}_t,t)$ of Eq.~\ref{appeq:sde_q} may be obtained from the solution of the variational problem of minimising the free energy
\begin{equation}\label{eq:free_energy}
    \mathcal{F}[Q] = \mathcal{KL}\Big(Q(\mathbf{X}_{0:T})||\PP(\mathbf{X}_{0:T}\mid \hat{\mathbf{f}}) \Big)- \sum \limits_{k=1}^K \Big\langle\ln \PP(\boldsymbol{\mathcal{O}}_{k}\mid \mathbf{X}_{t_k})\Big\rangle_{{Q}}.
\end{equation}

By applying the Cameron-Girsanov-Martin theorem we can express the Kullback-Leibler divergence between the two path measures induced by the diffusions with drift $\hat{\mathbf{f}}(\mathbf{x})$ and $\mathbf{g}(\mathbf{x},t)$ as
\begin{align}
    \mathcal{KL}\Big(Q(\mathbf{X}_{0:T})||\PP(\mathbf{X}_{0:T}|\hat{\mathbf{f}}) \Big) &= \Bigg\langle\text{ln}\left(\frac{\text{d} {Q}(\mathbf{X}_{0:T})}{\text{d} {P} \left( \mathbf{X}_{0:T} \vert \hat{\mathbf{f}}  \right)} \right) \Bigg\rangle_{{Q}} \\ 
    &=\Bigg\langle \left( - \frac{1}{2}  \int_0^T { { \|\hat{\mathbf{f}}(\mathbf{X}_t)-\mathbf{g}(\mathbf{X}_t,t) \|_{\mathbf{D}}^2}\text{d}t} + \int_0^T {\frac{ \hat{\mathbf{f}}(\mathbf{X}_t)-\mathbf{g}(\mathbf{X}_t,t)  }{{\mathbf{D}}} \text{d}\mathbf{W}_t} \right)\Bigg\rangle_{{Q}} \nonumber\\
    &=\Bigg\langle  \left( - \frac{1}{2}  \int_0^T { { \|\hat{\mathbf{f}}(\mathbf{X}_t)-\mathbf{g}(\mathbf{X}_t,t) \|_{\mathbf{D}}^2}\text{d}t} +V_T \right)\Bigg\rangle_{{Q}} \\
    &= \frac{1}{2} \int \limits^T_0 \int \| \mathbf{g}(\mathbf{x},t) -\hat{\mathbf{f}}(\mathbf{x})   \|_{\mathbf{D}}^2 \, q_t(\mathbf{x})\, \text{d}\mathbf{x} \, \text{d}t + \mathfrak{C} \label{eq:KL_1},
\end{align}
where $q_t(\mathbf{x})$ stands for the marginal density for $\mathbf{X}_t$ of the approximate process. In the third line we have introduced the random variable $V_T = \int_0^T {\frac{ \hat{\mathbf{f}}(\mathbf{X}_t)-\mathbf{g}(\mathbf{X}_t,t)  }{{\mathbf{D}}} \text{d}\mathbf{W}_t}$. Under the assumption that the function ${\ell(\mathbf{X}_t) = \hat{\mathbf{f}}(\mathbf{X}_t)-\mathbf{g}(\mathbf{X}_t,t)}$ is bounded, piece-wise continuous, and in $L^2[0,\infty)$ , $V_T$ follows the distribution $\mathcal{N}\left(V_T \mid 0, \int_0^T \ell^2(s) \text{d}s\right)$, which for a given $T$ will result into a constant $\mathfrak{C}$. Thus the second term in Eq.~\ref{eq:KL_1} is not relevant for the minimisation of the free energy and will be omitted. 

We can thus express the free energy of Eq.~\ref{eq:free_energy} as~\citep{opper2019variational}
\begin{equation} \label{appeq:free_energy2}
    \mathcal{F}[Q] = \frac{1}{2} \int \limits^T_0 \int \Big[ \|\mathbf{g}(\mathbf{x},t) - \hat{\mathbf{f}}(\mathbf{x})\|_{\mathbf{D}}^2 + U(\mathbf{x},t)   \Big] \, q_t(\mathbf{x})\, \text{d}\mathbf{x} \,\text{d}t,
\end{equation}
where the term $U(\mathbf{x},t)$ accounts for the observations $U(\mathbf{x},t) = - \sum \limits_{t_k} \ln \PP(\boldsymbol{\mathcal{O}}_k \mid \mathbf{x})\, \delta(t- t_k)$.

The minimisation of the functional of the free energy can be construed as a stochastic control problem~\citep{opper2019variational} with the objective to identify a time-dependent drift adjustment $\mathbf{u}(\mathbf{x},t):=\mathbf{g}(\mathbf{x},t) - \hat{\mathbf{f}}(\mathbf{x})$ for the system with drift $\hat{\mathbf{f}}(\mathbf{x})$ so that the controlled dynamics fulfil the constraints imposed by the observations.

\subsubsection{Approximate posterior over paths {with} geometric constraints}\label{appsec:with}

The previously described construction of the approximate measure in terms of stochastic bridges is relevant when the observations have non vanishing probability under the law of the prior diffusion process of Eq.~\ref{appeq:prior_sde}. However, when the prior process (with the estimated drift $\hat{f}$) differs considerably from the process that generated the observations, such a construction might either provide a bad approximation of the underlying path measure, or show slow numerical convergence in the construction of the diffusion bridges.
To overcome this issue, we consider here additional constraints for the posterior process that force the paths of the posterior measure to respect the local geometry of the observations. In the following we provide a brief introduction on the basics of Riemannian geometry and consequently continue with the geometric considerations of the proposed method.

\paragraph{Riemannian geometry.}
A $d$-dimensional \textbf{Riemannian manifold}~\citep{do1992riemannian,lee2018introduction} $\left(\mathcal{M}, \boldsymbol{\mathfrak{h}} \right)$ embedded in a $d$-dimensional ambient space $\mathcal{X} = \mathbb{R}^d$ is a smooth curved $d$-dimensional surface 
endowed with a smoothly varying inner product (Riemannian) \textbf{metric} $\boldsymbol{\mathfrak{h}}: \mathbf{x} \rightarrow \langle \cdot | \cdot \rangle_{\mathbf{x}}$ on $\mathcal{T}_{\mathbf{x}}\mathcal{M}$. A tangent space $\mathcal{T}_{\mathbf{x}} \mathcal{M}$ is defined at each point $\mathbf{x} \in \mathcal{M}$. The Riemannian metric $\boldsymbol{\mathfrak{h}}$ defines a canonical volume measure on the manifold $\mathcal{M}$. Intuitively this characterises how to compute inner products locally between points on the tangent space of the manifold $\mathcal{M}$, and therefore determines also how to compute norms and thus distances between points on $\mathcal{M}$.

A \textbf{coordinate chart} $(G,\phi)$ provides the mapping from an open set $G$ on $\mathcal{M}$ to an open set $V$ in the Euclidean space. The dimensionality of the manifold is $d$ if for each point $\mathbf{x}\in \mathcal{M}$ there exists a local neighborhood $G \subset  \mathbb{R}^d$.
We can represent the metric $\mathfrak{h}$ on the local chart $(G,\phi)$ by the positive definite matrix (\textbf{metric tensor}) ${H(\mathbf{x}) = (\mathfrak{h}_{i,j})_{\mathbf{x}, 0 \leq i,j,\leq d} = \left( \langle  \frac{\partial}{\partial x_i}|  \frac{\partial}{\partial x_j}\rangle_{\mathbf{x}}  \right)_{0 \leq i,j,\leq d} }$ at each point $\mathbf{x} \in G$.

For $\mathbf{v},\mathbf{w} \in \mathcal{T}_{\mathbf{x}}\mathcal{M}$ and $\mathbf{x} \in G$, their inner product can be expressed in terms of the matrix representation of the metric  $\boldsymbol{\mathfrak{h}}$ on the tangent space $\mathcal{T}_{\mathbf{x}}\mathcal{M}$ as $\langle \mathbf{v}|\mathbf{w} \rangle_{\mathbf{x}} = \mathbf{v}^{\top} H(\mathbf{x})\mathbf{w}$, where $H(\mathbf{x})\in \mathbb{R}^{d \times d}$ .

The \textbf{length of a curve} $\gamma:[0,1]\rightarrow \mathcal{M}$ on the manifold is defined as the integral of the norm of the tangent vector 
\begin{equation}\label{eq:ell}
\ell(\gamma_{t'}) = \int^1_0\| \dot{\gamma}_{t'}\|_{\boldsymbol{\mathfrak{g}}} \text{d}t' = \int^1_0 \sqrt{ \dot{\gamma}_{t'}^{\top} H(\gamma_{t'}) \dot{\gamma}_{t'}    } \text{d}t',\end{equation}
where the dotted letter indicates the velocity of the curve $\dot{\gamma}_{t'}=\partial_{t'} \gamma_{t'}$. A \textbf{geodesic curve} is a locally length minimising smooth curve that connects two given points on the manifold.

\paragraph{Riemannian geometry of observations.} 
For approximating the posterior over paths we take into account the geometry of the invariant density as it is represented by  the observations.  
To that end, we consider systems whose dynamics induce invariant (inertial) manifolds that contain the global attractor of the system and on which system trajectories concentrate~\citep{wiggins1994normally,mohammed1999stable, girya1995inertial, fenichel1971persistence, arnold1990stochastic, carverhill1985flows}. We assume thus that the continuous-time trajectories $\mathbf{X}_{0:T} \in \mathbb{R}^d$ of the underlying system concentrates on an invariant manifold $\mathcal{M} \in \mathbb{R}^{m \leq d}$ of dimensionality $m$ (possibly) smaller than $d$.
The discrete-time observations $\boldsymbol{\mathcal{O}}_k$ are thus samples of the manifold $\mathcal{M}$.
The central premise of our approach is that \textbf{unobserved paths between successive observations will be lying either \emph{on} or \emph{in the vicinity} of the manifold} $\mathcal{M}$. In particular, we postulate that unobserved paths should lie \textbf{in the vicinity of geodesics that connect consecutive observations} on $\mathcal{M}$. To that end we propose a path augmentation framework that constraints the augmented paths to lie in the vicinity of identified geodesics between consecutive observations.

However, while this view of a lower dimensional manifold embedded in a higher dimensional ambient space helps to build our intuition for the proposed method, for computational purposes we adopt a complementary view inspired by the discussion in~\citep{frohlich2021bayesian}. According to this view, we consider the entire observation space $\mathbb{R}^d$ as a smooth Riemannian manifold, $\mathcal{M}\dot{=}\mathbb{R}^d$, characterised by a Riemannian metric $\boldsymbol{\mathfrak{h}}$. The effect of the nonlinear geometry of the observations is then captured by the metric $\boldsymbol{\mathfrak{h}}$. Thus to approximate the geometric structure of the system's invariant density, we learn the Riemannian metric tensor $H:\mathbb{R}^d \rightarrow \mathbb{R}^{d \times d} $ and compute the geodesics between consecutive observations according to the learned metric. Intuitively according to this view the observations $\{\boldsymbol{\mathcal{O}}_k\}^K_{k=1}$ introduce distortions in the way we compute distances on the state space.

In effect this approach does not reduce the dimensionality of the space we operate, but changes the way we compute inner products and thus distances, lengths, and geodesic curves on $\mathcal{M}$. The alternative perspective of working on a lower dimensional manifold would strongly depend on the correct assessment of the dimensionality of said manifold. For example, one could use a Variational Autoencoder to approximate the observation manifold and subsequently obtain the Riemannian metric from the embedding of the manifold mediated by the decoder. 
However, our preliminary results of such an approach revealed that such a method requires considerable fine tuning to adapt to the characteristics of each dynamical system and is sensitive to the estimation of the dimensionality of the approximated manifold. 

To learn the Riemannian metric and compute the geodesics we follow the framework proposed by Arvanitidis et al. in~\citep{arvanitidis2019fast}.
In particular, we approximate the local metric induced by the observations at location $\mathbf{x}$ of the state space, in a non-parametric form by the inverse of the weighted local diagonal covariance computed on the observations as~\citep{arvanitidis2019fast}
\begin{equation}\label{appeq:metric}
    H_{dd}(\mathbf{x}) = \left(  \sum\limits^K_{i=1} w_i(\mathbf{x}) \left( x^{(d)}_i - x^{(d)}\right)^2 + \epsilon   \right)^{-1},
\end{equation}
with weights $w_i(\mathbf{x}) = \exp \left(- \frac{\|  \mathbf{x}_i - \mathbf{x} \|^2_2}{2 \sigma^2_{\mathcal{M}}}  \right)$, and $x^{(d)}$ denoting the $d$-th dimensional component of the vector $\mathbf{x}$. The parameter $\epsilon > 0$ ensures non-zero diagonals of the weighted covariance matrix, while $\sigma_{\mathcal{M}}$ characterises the curvature of the manifold.

Between consecutive observations for each interval $[\boldsymbol{\mathcal{O}}_k, \boldsymbol{\mathcal{O}}_{k+1}]$, we identify the geodesic $\boldsymbol{\gamma}^k_{t'}$ as the energy minimising curve, i.e., as the minimiser of the kinetic energy functional ${\mathcal{E}(\boldsymbol{\gamma}^k_{t'}) =\int^1_0 L_{{\mathcal{M}}}(\boldsymbol{\gamma}^k_{t'}, \dot{\boldsymbol{\gamma}}^k_{t'})\, \text{d}t'}$
\begin{equation} \nonumber
  \boldsymbol{\gamma}^{k*}_{t'} =  \underset{\boldsymbol{\gamma}^k_{t'}, \boldsymbol{\gamma}^k_0 = \boldsymbol{\mathcal{O}}_k, \boldsymbol{\gamma}^k_1=\boldsymbol{\mathcal{O}}_{k+1}}{\arg\min} \int^1_0 L_{{\mathcal{M}}}(\boldsymbol{\gamma}^k_{t'}, \dot{\boldsymbol{\gamma}}^k_{t'}) \,\text{d}t',
\end{equation}

\begin{equation} \label{eq:machlup}
\text{with} \;\;\;\;  \int^1_0 L_{{\mathcal{M}}}(\gamma^k_{t'}, \dot{\gamma}^k_{t'}) \text{d}t'= \frac{1}{2}  \int^1_0 \|\dot{\gamma}^k_{t'} \|^2_{\boldsymbol{\mathfrak{h}}}  ,
\end{equation} \label{eq:geodesic_lagrangian}
where $L_{{\mathcal{M}}}(\boldsymbol{\gamma}^k_{t'}, \dot{\boldsymbol{\gamma}}^k_{t'})$ denotes the Lagrangian.
The minimising curve of this functional is the same as the minimiser of the curve length functional $\ell(\boldsymbol{\gamma}_{t'})$ (Eq.~\ref{eq:ell}), i.e., the geodesic~\citep{do1992riemannian}.

By applying calculus of variations, the minimising curve of the functional $\mathcal{E}(\boldsymbol{\gamma}^k_{t'}) $ can be obtained from the Euler-Lagrange equations, resulting in the following system of second order differential equations~\citep{arvanitidis2017latent,do1992riemannian}
\begin{equation}\label{eq:geode}
    \ddot{\boldsymbol{\gamma}_t}^k = -\frac{1}{2} {H(\boldsymbol{\gamma}^k_t)}^{-1} \Bigg( 2 \left( I \otimes (\dot{\boldsymbol{\gamma}_t}^k)^{\top} \right) \frac{\partial \text{vec}[H(\boldsymbol{\gamma}^k_t)]}{\partial \boldsymbol{\gamma}^k_t}  \dot{\boldsymbol{\gamma}_t}^k  - \frac{\partial \text{vec}[H(\boldsymbol{\gamma}^k_t)]^{\top}}{\partial \boldsymbol{\gamma}^k_t} \left(\dot{\boldsymbol{\gamma}_t}^k \otimes \dot{\boldsymbol{\gamma}_t}^k  \right)\Bigg),
\end{equation}
with boundary conditions $\boldsymbol{\gamma}^k_0 = \boldsymbol{\mathcal{O}}_k $ and $ \boldsymbol{\gamma}^k_1=\boldsymbol{\mathcal{O}}_{k+1}$,
where $\otimes$ stands for the Kroenecker product, and $\text{vec}[A]$ denotes the vectorisation operation of matrix $A$ through stacking the columns of $A$ into a vector.
We follow~\cite{arvanitidis2019fast} and obtain the geodesics by approximating the solution of the boundary value problem of Eq.~\ref{eq:geode} with a probabilistic differential equation solver.

\paragraph{Extended free energy functional.} We denote the collection of individual geodesics by $\boldsymbol{\Gamma}_t\dot{=} \{\boldsymbol{\gamma}^k_{t'}\}_{t=(k-1)\tau+t' \tau}$, where $\boldsymbol{\gamma}^k_{t'}$ is the geodesic connecting $\boldsymbol{\mathcal{O}}_k$ and $\boldsymbol{\mathcal{O}}_{k+1}$, and $t'\in [0,1]$ denotes a rescaled time variable. Additional to the constraints imposed in the previously explained setting (Sec~\ref{appsec:without}), here we add an extra term in the free energy ${U_{\mathcal{G}}(\mathbf{x},t) \dot{=} \| \boldsymbol{\Gamma}_t -\mathbf{x} \|^2 }$ that accounts for the local geometry of the invariant density, and guides the latent path towards the geodesic curves $\gamma^k_{t'}$ that connect consecutive observations 
\begin{equation} \label{apeq:free_energy2}
 \small   \mathcal{F}[Q] = \frac{1}{2} \int \limits^T_0 \int \Big[ \|g(\mathbf{x},t) - \hat{f}(\mathbf{x})\|^2_D + U_{\mathcal{O}}(\mathbf{x},t) + \beta U_{\mathcal{G}}(\mathbf{x},t)   \Big] \, q_t(\mathbf{x})\, d\mathbf{x} \,\text{d}t.
\end{equation}
Here we denote the observation term by $U_{\mathcal{O}}(\mathbf{x},t) \dot{=} - \sum_{t_k} \ln \PP(\boldsymbol{\mathcal{O}}_k|\mathbf{x}) \delta(t- t_k)$, while $\beta$ stands for a weighting constant that determines the relative weight of the geometric term in the control objective.

Eq.~\ref{apeq:free_energy2} can be construed as a least-action principle for the augmented path measure. The first term is the energetic cost of changing the drift from ($\hat f$) to ($g=\hat f+u$). The observation potential ($U_{\mathcal{O}}$) imposes the bridge boundary conditions at the sparse measurement times. The geometric potential ($U_{\mathcal{G}}$) does not force the path to equal the geodesic, but penalises deviations from it, with strength controlled by ($\beta$). Therefore, the optimal-control formulation provides a principled way to combine three requirements: dynamical consistency with the current drift estimate, interpolation through the observations, and regularisation by the invariant-density geometry.

 Following~\citep{opper2019variational}, for each inter-observation interval $[\boldsymbol{\mathcal{O}}_k, \boldsymbol{\mathcal{O}}_{k+1}]$ we identify the posterior path measure (minimiser of Eq.~\ref{apeq:free_energy2}) by the solution of a stochastic optimal control problem~\citep{maoutsa2022deterministic} with the objective to obtain a time-dependent drift adjustment ${\mathbf{u}(\mathbf{x},t):=\mathbf{g}(\mathbf{x},t) - \hat{\mathbf{f}}(\mathbf{x})}$ for the system with drift $\hat{\mathbf{f}}(\mathbf{x})$ with initial and terminal constraints defined by $U_{\mathcal{O}}(\mathbf{x},t)$, and additional path constraints $U_{\mathcal{G}}(\mathbf{x},t)$.

 For the case of exact observations, i.e., for an observation process $\boldsymbol{\psi}(\mathbf{x}) = \mathbf{x}$, we can compute the drift adjustment for each of the $K-1$ inter-observation intervals independently. Thus for each interval between consecutive observations, we identify the optimal control $\mathbf{u}(\mathbf{x},t)$ required to construct a stochastic bridge following the dynamics of Eq.~\ref{appeq:prior_sde} with initial and terminal states the respective observations $\boldsymbol{\mathcal{O}}_k$ and $\boldsymbol{\mathcal{O}}_{k+1}$. 

The optimal drift adjustment for such a stochastic control problem for the inter-observation interval between $\boldsymbol{\mathcal{O}}_k$ and $\boldsymbol{\mathcal{O}}_{k+1}$ can be obtained from the solution of the backward equation (see~\citep{maoutsa2022deterministic, maoutsa2021deterministica})
\begin{equation} \label{app:phipde}
    \frac{\partial \phi_t(\mathbf{\mathbf{x}})}{\partial t} = - \mathcal{L}_{\hat{f}}^{\dagger} \phi_t(\mathbf{x}) + U_{\mathcal{G}}(\mathbf{x},t) \phi_t(\mathbf{x}),
\end{equation}
with terminal condition $\phi_T(\mathbf{x}) = \chi(\mathbf{x}) = \delta(\mathbf{x}-\boldsymbol{\mathcal{O}}_{k+1}) $ and with $\mathcal{L}_{\hat{f}}^{\dagger}$ denoting the adjoint Fokker-Planck operator for the process of Eq.~\ref{appeq:prior_sde}.
As shown in~\citep{maoutsa2022deterministic} the optimal drift adjustment $\mathbf{u}(\mathbf{x},t)$ can be expressed in terms of the difference of the logarithmic gradients of two probability flows
\begin{equation}
    \mathbf{u}^*(\mathbf{x},t) = D \Big( \nabla \ln q_{T-t}(\mathbf{x}) - \nabla \ln \rho_t(\mathbf{x}) \Big),
\end{equation}
where $\rho_t$ fulfils the forward (filtering) partial differential equation (PDE)
 \begin{equation} 
\frac{\partial \rho_t(\mathbf{x})}{\partial t} = {\cal{L}}_{\hat{f}} \rho_t(\mathbf{x}) - U_{\mathcal{G}}(\mathbf{x},t) \rho_t(\mathbf{x}),
\label{eq:FPE2} 
\end{equation}
while $q_t$ is the solution of a time-reversed PDE that depends on the logarithmic gradient of $\rho_t(\mathbf{x})$
\begin{align}\label{Fokker_bridge3}
\frac{\partial {q}_{t}(\mathbf{x})}{\partial t} &= 
-\nabla\cdot \Bigg[\Big(\sigma^2\nabla \ln  \rho_{T-t} (\mathbf{x})  - \mathbf{f}(\mathbf{x}, T-t)\Big)  {q}_{t} (\mathbf{x})\Bigg] +  \frac{\sigma^2}{2} \nabla^2 {q}_{t} (\mathbf{x}) , 
\end{align}
with initial condition ${q}_{0} (\mathbf{x}) \propto \rho_T(\mathbf{x}) \chi(\mathbf{x})$
.

We obtain the optimal control by solving the backward equation associated with this controlled diffusion problem. The solution defines the change of path measure required to condition the prior process on the endpoint and path potentials. The resulting control can be expressed as a logarithmic gradient of the corresponding desirability or path-density functions, resulting in a time-dependent correction ($u^\ast(x,t)$). In the implementation, we apply this construction independently to each interval ($[\mathcal{\boldsymbol{O}}_k,\mathcal{\boldsymbol{O}}_{k+1}]$), obtaining one controlled bridge per pair of consecutive observations.

 For the numerical solution of the control problem we use the numerical framework accompanying ~\cite{maoutsa2022deterministic}, where the path constraints associated with the geodesic curves are imposed through the two staged process for particle propagation described in the paper for path constraints, with the particle reweighting being performed through optimal transport implemented using the PyEMD python toolbox~\citep{pele2009}.

More precisely, according to this framework we propagate a particle representation of the probability density $\rho_t(\mathbf{x})$ according to the filtering equation of Eq.~\ref{eq:FPE2}. This follows the dynamics of the uncontrolled process with drift $\hat{\mathbf{f}}$ and particle reweighting at each time step as determined by the path constrained (potential) $U_{\mathcal{G}}(\mathbf{x},t)$, that quantifies the proximity to the geodesic at each time point. In the particle representation we apply this reweighting in the form of a deterministic optimal transportation of the particles
~\citep{reich2013nonparametric}.

\subsection{Approximate posterior over drift functions.}\label{appsec:drift_inference}

For a fixed path measure ${Q}$, the optimal measure for the drift ${Q}_f$ is a Gaussian process given by
\begin{equation} \label{appeq:drift_measure}
    {Q}_f \propto \PP_f \exp\left({ -\frac{1}{2} \int  \|\mathbf{f}(\mathbf{x})\|_{D}^2 A(\mathbf{x}) - 2 \langle \mathbf{f}(\mathbf{x}), B(\mathbf{x}) \rangle_{D}  } \text{d}\mathbf{x}\right),
\end{equation}
with $$A(\mathbf{x})\dot{=} \int^T_{0} q_t(\mathbf{x}) \text{d}t,$$ and $$B(\mathbf{x})\dot{=} \int^T_{0} q_t(\mathbf{x}) g(\mathbf{x},t) \text{d}t, $$ where $q_t(\mathbf{x})$ denotes the marginal constrained density of the state $\mathbf{X}_t$. The function $\mathbf{g}(\mathbf{x},t)$ denotes the effective drift.

 We assume a Gaussian process prior for the unknown function $\mathbf{f}$, i.e., $\mathbf{f} \sim \PP_0(\mathbf{f}) =\mathcal{GP}(\boldsymbol{m}^f, k^f)$ where $m^f$ and $k^f$ denote the mean and covariance function of the Gaussian process. Following Ruttor \emph{et al.}~\citep{ruttor2013approximate}, we employ a sparse kernel approximation for the drift $f$ by optimising the function values over a sparse set of $S$ inducing points $\{Z_i\}^{S}_{i=1}$.
We obtain the resulting drift from
\begin{equation}
    \hat{\mathbf{f}}_S(\mathbf{x}) = k^\mathbf{f}(\mathbf{x},\mathcal{Z}) \left( I + \Lambda \, \mathcal{K}_S  \right)^{-1} \mathbf{d},
\end{equation}
where we have defined introduced the notation $\mathcal{K}_S \dot{=} k^f(\mathcal{Z},\mathcal{Z}) $
\begin{equation}
    \Lambda = \frac{1}{\sigma^2} \mathcal{K}^{-1}_S \left(   \int k^f(\mathcal{Z},\mathbf{x}) A(\mathbf{x}) k^f(\mathbf{x},\mathcal{Z}) \text{d}\mathbf{x} \right)   \mathcal{K}^{-1}_S.
\end{equation}

\begin{equation}
    \mathbf{d} = \frac{1}{\sigma^2} \mathcal{K}^{-1}_S \left(   \int k^f(\mathcal{Z},\mathbf{x}) B(\mathbf{x}) \text{d}\mathbf{x} \right)   \mathcal{K}^{-1}_S,
\end{equation}

The associated variance results similarly from the equation
\begin{equation}
    \Sigma^2_S (\mathbf{x}) = k^f(\mathbf{x},\mathbf{x}) - k^f(\mathbf{x},\mathcal{Z}) \left( I + \Lambda \, \mathcal{K}_S  \right)^{-1} \Lambda\, k^f(\mathcal{Z},\mathbf{x}).
\end{equation}

We employ a sample based approximation of the densities in Eq.~\ref{appeq:drift_measure} resulting from the particle sampling of the path measure $Q$ resulting from the geometric augmentation, i.e. the integrals over $\int q_t(\mathbf{x})$ are over the samples of the augmented paths. 
Thus by representing the densities by samples, we can rewrite the density $p_{t}(x)$ in terms of a sum of Dirac delta functions centered around the particles positions
$$p_{t}(\mathbf{x}) \approx \frac{1}{N} \sum^N_{j=1} \delta(\mathbf{x} - \textbf{X}_j(t)),$$ and replace the Riemannian integrals with summation over particles, i.e. perform a Monte Carlo integration. Here $\textbf{X}_j(t)$ represents the position of the $j$-th particle at time point $t$.

\section{Sparse Gaussian process estimation}
Since the amount of required observations for accurate drift estimation is generally large for systems with nonlinear dynamics, regular Gaussian process regression becomes computationally intensive. Its computational complexity scales as $\mathcal{O}(N^3)$ with the number of observations $N$ due to the $N \times N$ kernel matrix inversions required for inference (c.f. Eq.\ref{eq7:full_gp_cov} and~\cite{rasmussen2003gaussian}). Therefore, \cite{ruttor2013approximate} employ the sparse (low dimensional approximation) counterpart of Gaussian process regression~\citep{titsias2009variational, csato2002sparse} that reduces significantly the computation time by reducing the computational complexity to $\mathcal{O}(NM^2)$, where $M \ll N$ denotes the number of selected sparse (inducing) points. Here we present briefly the derivation. 

For sparse Gaussian process drift inference, we augment the distributions with $M$ inducing points $\mathbf{z} = \left[z_1, \dots, z_M \right] $ with inducing values $\mathbf{u} = \left[ \mathbf{f}(z_m)\right]_{m=1}^M$ that are jointly Gaussian distributed with the latent function values $\{\mathbf{f}(\mathbf{X}_t)\}_{t=0}^T$.

As demonstrated previously the true posterior for function values $\mathbf{f}$ is expressed as a product
\begin{equation} \label{eq7:posterior_GP1}
P_f(\mathbf{f}) = \frac{1}{Z} P_o(\mathbf{f}) e^{-\mathcal{A}(\mathbf{f})},
\end{equation}
where $Z$ a normalisation constant, $\mathcal{A}(\mathbf{f}) = \frac{1}{2}\mathbf{f}^T \Lambda \mathbf{f} - \mathbf{a}^T \mathbf{f}$ a quadratic form of $\mathbf{f}$ (see Eq.~\ref{apeq:SDE_likelihood}), while $P_o(\mathbf{f})$ denotes a prior Gaussian measure. Thus the posterior $P_f(\mathbf{f})$ is also Gaussian.

To employ sparse Gaussian process inference, we approximate $P_f$ with ${Q_f=\mathcal{GP}\left(m^q(\cdot), k^q(\cdot,\cdot) \right)}$, with mean and variance functions to be calculated, depending only on the \emph{smaller} subset ($M \ll N$) of inducing function values $\mathbf{u}$,
\begin{equation}
Q_f(\mathbf{f}) \propto {R(\mathbf{u})} P_o(\mathbf{f}).
\end{equation}
The effective likelihood $R(\mathbf{u})$ is chosen as the minimiser of the Kullback-Leibler divergence $\mathcal{KL}\left(Q_f||P_f  \right)$.

We may now express the prior $P_o(\mathbf{f})$ and the approximate marginal $Q_f(\mathbf{f})$ in terms of the inducing points
\begin{equation} \label{eq7:prior_u}
P_o(\mathbf{f}) = P_o(\mathbf{f}|\mathbf{u}) P_o(\mathbf{u}),
\end{equation}and
\begin{equation}\label{eq7:qfu}
Q_f(\mathbf{f}) = Q_f(\mathbf{f}|\mathbf{u}) Q_f(\mathbf{u})=P_o(\mathbf{f}|\mathbf{u}) Q_f(\mathbf{u}),
\end{equation}
under the assumption that the posterior conditional $Q_f(\mathbf{f}|\mathbf{u})$ matches the prior conditional $ P_o(\mathbf{f}|\mathbf{u})$.

We select the effective likelihood $R(u)$ as the minimiser of the relative entropy between $Q_f$ and $P_f$
\begin{equation} \label{eq7:kl_effective_likelihood}
\begin{split}
\mathcal{KL}\left( Q_f|| P_f  \right) &= \int Q_f(\mathbf{f})\; \text{ln}\frac{Q_f(\mathbf{f})}{P_f(\mathbf{f})} \de \mathbf{f}\\
&= \int P_o(\mathbf{f}|\mathbf{u}) Q_f(\mathbf{u})\; \text{ln}\frac{P_o(\mathbf{f}) R(\mathbf{u})}{\frac{1}{Z} P_o(\mathbf{f}) e^{-\mathcal{A}(\mathbf{f})}} \de \mathbf{f} \de \mathbf{u}\\
&= \int P_o(\mathbf{f}|\mathbf{u}) Q_f(\mathbf{u}) \; \text{ln}\frac{P_o(\mathbf{f}) R(\mathbf{u})}{\frac{1}{Z}P_o(\mathbf{f}|\mathbf{u}) e^{-\mathcal{A}(\mathbf{f}|\mathbf{u})}P_o(\mathbf{u})} \de \mathbf{f} \de \mathbf{u}\\
&= \int P_o(\mathbf{f}|\mathbf{u}) Q_f(\mathbf{u})\; \text{ln}\frac{P_o(\mathbf{u}) R(\mathbf{u})}{\frac{1}{Z} e^{-\mathcal{A}(\mathbf{f}|\mathbf{u})}P_o(\mathbf{u})} \de \mathbf{f} \de \mathbf{u}\\
&= \int P_o(\mathbf{f}|\mathbf{u}) Q_f(\mathbf{u}) \;\text{ln}\frac{ R(\mathbf{u})}{\frac{1}{Z} e^{-\mathcal{A}(\mathbf{f}|\mathbf{u})}} \de \mathbf{f} \de \mathbf{u}\\
&= \text{ln}Z + \int{Q_f(\mathbf{u})} \;\text{ln}\left( \frac{e^{\text{ln}R(\mathbf{u})}}{  e^{-\mathbb{E}_o \left[ \mathcal{A}(\mathbf{f}|\mathbf{u})  \right]  }}   \right)\de \mathbf{u}.
\end{split}
\end{equation}

In Eq.~\ref{eq7:kl_effective_likelihood} in the second line, we have introduced Eq.~\ref{eq7:posterior_GP1}-Eq.~\ref{eq7:qfu}. In the third line we have introduced $\frac{P_0(\mathbf{f})}{P_0(\mathbf{f}|\mathbf{u})} = P_0(\mathbf{u})$ from Eq.~\ref{eq7:prior_u}. In the final line we rearranged the terms that do not depend on $\mathbf{f}$  outside of the integral over $\mathbf{f}$, moved the $\ln Z$ term out of the integration over $\mathbf{u}$, and denoted $\langle\cdot\rangle_0 = \int P_0(\mathbf{f}|\mathbf{u}) \; \de \mathbf{f}$.
 
 To minimise the relative entropy $\mathcal{KL}\left[ Q_f|| P_f  \right]$  we conclude that the optimal choice for the effective likelihood $\label{R} R(\mathbf{u}) $ is 
 \begin{equation} \label{Ruu}
 R(\mathbf{u}) \propto e^{-\langle \mathcal{A}(\mathbf{f}|\mathbf{u})  \rangle_o  } .
 \end{equation}
  
  Given the quadratic form of $A(\mathbf{f})$ we may write the conditional expectation in Eq.~\ref{Ruu} as a quadratic form too
  \begin{equation}\label{app:expect_f}
  \begin{split}
  \langle \mathcal{A}(\mathbf{f}|\mathbf{u}) \rangle_o &= \frac{1}{2} \langle \mathbf{f}|\mathbf{u} \rangle_o^\top \,\Lambda \,\langle\mathbf{f}|\mathbf{u} \rangle_o + \frac{1}{2} \text{Tr}\left(  \text{Cov}_o[\mathbf{f}|\mathbf{u}] \Lambda\right)- a^\top \langle\mathbf{f}|\mathbf{u} \rangle_o \\
  &= \frac{1}{2} \langle\mathbf{f}|\mathbf{u} \rangle_o ^\top \; \Lambda \;\langle\mathbf{f}|\mathbf{u} \rangle_o - a^{\top} \langle\mathbf{f}|\mathbf{u} \rangle_o + \text{const.},
  \end{split}
  \end{equation}
 where in the last line we take into account that the term $\text{Tr}\left(  \text{Cov}_o[\mathbf{f}|\mathbf{u}] \Lambda\right)$ is independent of the sparse function values $\mathbf{u}$ (c.f.~\cite{ruttor2013approximate}). In Eq.~\ref{app:expect_f} $
\Lambda \;\dot{=}\; \mathrm{diag}\!\big[\Delta t\, D^{-1},\ldots,\Delta t\, D^{-1}\big].$
 
 In particular, the conditional expectation of function values $f$ conditioned on the inducing point function values $\mathbf{u}\equiv \mathcal{U}$ at inducing point locations $\mathbf{z}\equiv\mathcal{Z}$ equals 
 \begin{equation}
    \bar{f}^s(\mathbf{x}) =  \langle f|\mathbf{u} \rangle_o = k(\mathbf{x},\mathcal{Z}) k(\mathcal{Z},\mathcal{Z})^{-1}\mathcal{U},
 \end{equation}
 while the covariance equals
 \begin{equation}
     (\Sigma^s)^2(\mathbf{x}) = k(\mathbf{x},\mathbf{x}) - k(\mathbf{x},\mathcal{Z})k(\mathcal{Z}, \mathcal{Z})^{-1}k(\mathcal{Z},\mathbf{x}),
 \end{equation}
 where we have employed similar notation for the kernel functions as in Eqs.~\ref{eq7:full_gp_mean}-\ref{eq7:full_gp_cov}.
 

\section{Theoretical evidence that supports the use of geodesics as geometric constraints representing most probable paths between observations} \label{app:OM}

The Onsager-Machlup functional for diffusion processes has been known in theoretical physics as a characteriser of the most probable path (MPP) between two pre-defined states of the process. In~\citep{onsager1953fluctuations}, Onsager and Machlup used the thermal fluctuations of a diffusion process to show that the probability density of a path $\boldsymbol{\gamma} \in C^1 \left([0,T], \mathbb{R}^d \right)$ in $\mathbb{R}^d$ over finite interval can be expressed as a Boltzmann factor 
\begin{equation}
\PP(\boldsymbol{\gamma}) \sim \exp \left[ - \int^{T}_{0} L(\boldsymbol{\gamma}(t), \dot{\boldsymbol{\gamma}}(t) ) dt \right],
\end{equation}
where \begin{equation}L(\boldsymbol{\gamma}(t), \dot{\boldsymbol{\gamma}}(t) ) = \frac{1}{2} \| {\dot{\boldsymbol{\gamma}}(t) - \mathbf{f}(\boldsymbol{\gamma}(t))}\|_{\mathbf{D}^{-1}}^2 + \frac{1}{2} \nabla \cdot \mathbf{f}(\boldsymbol{\gamma}(t)).\footnote{Onsager and Machlup's initial work concentrated around linear processes and therefore the functional initially introduced by the did not include the second term with the divergence of $\mathbf{f}$ as this is a constant for linear $\mathbf{f}$. It was later added to the OM function to account for trajectory entropy corrections~\citep{taniguchi2007onsager,adib2008stochastic}}\end{equation}
The function $L(\boldsymbol{\gamma}(t), \dot{\boldsymbol{\gamma}}(t) )$ is known as the \textbf{Onsager-Machlup} function (action), while its integral over time is known as Onsager-Machlup action functional. It has been used as Lagrangian in Euler-Lagrange minimisation schemes to identify the \textbf{most probable path (MPP)} of a diffusion process between two given points in the state space~\citep{graham1977path, stratonovich1971probability}.

Stratonovich~\citep{stratonovich1971probability} considered the probability that a sample of a multidimensional diffusion process will lie in the vicinity of (within a tube of infinitesimal thickness around) an idealised smooth path in the state space. To compute this probability he constructed a probability functional which is identical to the Onsager-Machlup functional considered as Lagrangian for the diffusion process.
 Duerr et al.~\citep{durr1978onsager} considered scalar diffusion processes and constructed the Onsager-Machlup function from the asymptotic limit of the transition probability between the starting and end state of the path using a Girsanov transformation. 


\if False
Fujita et al.~\citep{fujita1982onsager} have shown that the log probability of a tube of width $\epsilon \rightarrow 0$ in the state space has limiting value value 
\begin{equation}
    \log \PP(\gamma) \rightarrow \mathfrak{C}_1 + \mathfrak{C}_2/\epsilon^2 - \int^T_0 \frac{1}{2} \| \dot{\gamma}(t) \|^2 dt,
\end{equation}
where the term $\| \dot{\gamma}(t) \|^2$ can be identified as the energy of the path, and $\mathfrak{C}_1$,  $\mathfrak{C}_2$ denote constants. This already bares some similarities with our approach, since the curve that minimises the energy functional along a path in path space can be identified as a geodesic curve.
\fi

Considering Brownian motions defined on  a Riemannian manifold $(\mathcal{M}, \boldsymbol{\mathfrak{g}})$ with associated Riemannian metric $\boldsymbol{\mathfrak{g}}$, the Onsager-Machlup functional can be expressed as the integral over the Lagrangian~\citep{takahashi1981probability,graham1980onsager,grong2022most}
\begin{equation} \label{eq:OM_manifold}
    L(\boldsymbol{\gamma}, \dot{\boldsymbol{\gamma}}) = \frac{1}{2} \| \dot{\boldsymbol{\gamma}}(t) \|_{\boldsymbol{\mathfrak{g}}}^2 - \frac{1}{12} S(\boldsymbol{\gamma}(t)),
\end{equation}
where $\| \cdot \|_{\boldsymbol{\mathfrak{g}}}$ denotes the Riemannian norm on the tangent space $\mathcal{T}_X \mathcal{M}$ of the manifold with respect to the metric $\boldsymbol{\mathfrak{g}}$, and $S(\cdot)$ stands for the scalar curvature of the manifold at each point. The first term is the Lagrangian used to identify geodesic curves on manifolds (c.f. ~\ref{eq:geodesic_lagrangian})

In our proposed formalism, for computational purposes we have assumed the entire $\mathbb{R}^d$ as smooth manifold. We can identify the first term of Eq.~\ref{eq:OM_manifold} with the Lagrangian we optimised for computing the geodesics on the manifold $(\mathbb{R}^d, \boldsymbol{\mathfrak{g}})$, where $\boldsymbol{\mathfrak{g}}$ is the metric learned from the observations.

However the system we observed was a diffusion process defined in $\mathbb{R}^d$ with an Euclidean metric. Constructing a path augmentation scheme that guides the augmented paths towards the geodesics of a diffusion defined with respect to a different metric raises questions about the validity of our approach. Here we should note that diffusions with a general state dependent diffusion coefficient $\boldsymbol{\sigma} \in \mathbb{R}^{d \times m}$, and $m$-dimensional Brownian motion, can be considered as evolving on the manifold $\left(\mathbb{R}^{d}, \boldsymbol{\mathfrak{g}} \right)$, with the associated metric $\boldsymbol{\mathfrak{g}} = \left(\boldsymbol{\sigma \sigma}^{\top} \right)^{-1}$~\citep{capitaine2000onsager}. Thus it may be possible to associate the metric learned from the data with the metric arising from a state dependent diffusion by applying a transformation akin to an inverse Lamperti transform~\citep{oksendal2003stochastic} to transform our learned SDE to one that would have induced the learned metric due to the state dependent diffusion. The existence of such a transformation would justify the proposed method. Our empirical results demonstrate that such a transformation may be possible.

\section{Does the proposed approach invalidate the Markovian property of the diffusion process?}
The proposed path augmentation seemingly invalidates the Markovian property of the diffusion process.
According to the Markov property of the diffusion of Eq.~\ref{eq:system}, the system state $\mathbf{X}_{k \tau+\delta t}$ should depend only the state $\mathbf{X}_{k \tau}$, i.e., the observation $\boldsymbol{\mathcal{O}}_k$. The proposed augmentation makes the state $\mathbf{X}_{k \tau+\delta t}$ depending not only on the next observation $\boldsymbol{\mathcal{O}}_{k+1}=\mathbf{X}_{(k+1) \tau}$, but also on past and future states that lie in the vicinity of these observations. 

We effectively construct the augmented paths to compute the likelihood of a drift estimate. To compute this likelihood we require to evaluate the transition probabilities between consecutive observations. Since for general nonlinear systems the transition probabilities are in general intractable, we have to resort to numerical approximations. Ideally we would approximate the transition density with a bridge sampler that would consider the nonlinear estimated SDE conditioned to pass though consecutive observations. However for coarse drift estimates, the observations have zero probability under the law of the estimated SDE, and construction of those bridges would result either in very taxing computations or would fail altogether. Instead, here, we compute the likelihood of a "corrected" estimate (the correction resulting from the invariant density) under which the observations have non-zero probability, and subsequently re-estimate the drift on the augmented path with this "corrected" estimate.
 By taking into account the local geometry of the observations, we provide systematic corrections for the misestimated drift function to generate the augmented paths. This effectively nudges the augmentation process towards the second observation of each inter-observation interval through the path constraint that forces the augmented paths towards the geodesics.

\section{Related work and positioning of the present work}\label{appsec:related_work}

Here, we briefly review further related work on inference or modelling of SDEs and position our work further with respect to the existing literature.

Recent advances in dynamical system inference have delivered valuable tools for identifying continuous-time \emph{deterministic} systems from observations~\citep{cremers1987construction, brunton2016discovering, daniels2015automated, casadiego2018inferring, mcgoff2015statistical, kantz2004nonlinear, schmidt2009distilling, gerardos2025principled, packard1980geometry, crutchfield1987equations}. \textbf{Data-driven} (or \textbf{nonparametric}, or \textbf{equation-free}) approaches  
 seek to reconstruct the governing equations of observed systems directly from state observations, without imposing explicit assumptions or inductive biases about the underlying dynamical models. They rely on function approximation to infer the system’s structure from observations, such as basis functions~\citep{acosta1995radial,small2002minimum,judd1995selecting,small1998comparisons,bruckner2020inferring, frishman2020learning}, symbolic regression~\citep{kaiser2018sparse,brunton2016discovering, bongard2007automated,daniels2015automated}, spectral approximations~\citep{kevrekidis2003equation,theodoropoulos2000coarse}, Gaussian processes~\citep{alvarez2009latent,sanguinetti2006probabilistic,sarkka2019use}, or neural networks~\citep{teng2018machine, bhattoo2022learning,jung2019reco}.
 However, extending these methods to \emph{stochastic} systems remains difficult. In this setting, inference must disentangle the influence of underlying deterministic forces from random fluctuations, a task that is particularly difficult for low sampling rates.

\paragraph{{$\triangleright$~Modelling general SDEs from state observations}.}

As already mentioned in the Introduction and in Sec.~\ref{appsec:b} existing inference methods for SDEs can be broadly clustered in temporal and geometric methods, where the former accounts for the temporal order of the observations, while the latter approximate the invariant system density and discard any time information.

\textbf{Temporal methods} rely on the Euler-Maruyama discretisation of the SDE paths approximating conditional expectations of state increments (i.e. the Krammers Moyal coefficients). They model the drift either in terms of Gaussian processes~\citep{ruttor2013approximate, batz2018approximate, hostettler2018modeling, zhao2020state,yildiz2018learning}, basis functions~\citep{nabeel2025discovering, ragwitz2001indispensable,friedrich1997description, peinke1997new,  gao2024learning, friedrich2000extracting,ferretti2020building} or libraries of functions~\citep{boninsegna2018sparse,huang2022sparse},
 kernel regression~\citep{lamouroux2009kernel,jiang1997nonparametric}, dynamic mode decomposition to learn the eigenfunctions of the Koopman operator~\citep{klus2020data}, by approximating the central moments of the transition densities~\citep{stanton1997nonparametric}, Bayesian neural networks~\cite{bae2025inferring}, or by applying generalised methods of moments~\citep{hansen1993back}.

 As explicitly detailed in Sec.~\ref{appsec:b}, most temporal methods do not provide accurate drift estimates when the interval between observations is large. The two prevailing approaches to mitigate this finite-sampling rate effects is to either account for the systematic bias introduced by the finite sampling rate by estimating
 an explicit correction term for the inferred drift~\citep{ragwitz2001indispensable, ragwitz2002ragwitz, kleinhans2005iterative, kleinhans2007maximum}, or by performing state estimation for the unobserved paths (also known as path or data augmentation) and then estimating the drift from the continuous paths.

 The former approach works only for scalar systems, 
while the latter employs simplified bridge dynamics (e.g., Brownian~\citep{chib2006likelihood,eraker2001mcmc, sermaidis2013markov} or Ornstein Uhlenbeck~\citep{batz2018approximate, billio1998simulated} bridges) that are analytically tractable or computationally non-demanding. However, for large $\tau$ and for nonlinear systems, these simplified bridge dynamics match poorly the underlying path statistics.
 (Fig.~\ref{fig1} {\bf{D.}}). It is important to mention here, that path augmentation with Ornstein Uhlenbeck bridges similar to ~\cite{batz2018approximate} provides a good approximation of the underlying transition density, when the underlying linear process employed for each bridge has a drift that comes from the local linearisation of the \textbf{ground truth} drift function. However, during inference the ground truth dynamics are unknown and the local linearisations on inaccurate drift estimates employed in \cite{batz2018approximate} provide imprecise approximations for large $\tau$.

Alternative methods, employ variational inference~\citep{batz2016variational,opper2019variational, duncker2019learning, verma2024variational} and approximate the posterior path measure with a tractable Gaussian process induced by a time–varying linear SDE. This results in ODEs for the posterior mean and covariance matrix and an ELBO that is optimized directly \citep{archambeau2007variational, duncker2019learning}.

Building on the prolific line of work on neural ODEs,
neural SDEs~\citep{li2020scalable} employ
 gradient-based stochastic variational inference
and the
 stochastic adjoint sensitivity method to compute gradients of solutions of stochastic equations with respect to their parameters. Building on these methods, ~\cite{course2023amortized} remove the need for adjoint-based gradient computations by combining amortized inference with a reparametrization of the ELBO by assuming a latent linear process that generates the latent path.

\textbf{Geometric approaches} 
 on the other hand, discard the temporal structure of the observations, and treat them as samples of the invariant density. Thereby these methods either employ density estimation to identify the drift as the gradient of a potential~\cite{kutoyants2004statistical}, or resort to spectral approximations of the generator of the diffusion process through manifold learning.

Manifold learning methods employ often the \emph{diffusion maps} algorithm, introduced by Coiffman and colleagues ~\cite{singer2008non}, to learn the dominant part of the spectrum of the transfer operator of the observed diffusion process~\cite{coifman2005geometric,nadler2006diffusion,giannakis2019data,ferguson2011nonlinear,talmon2015intrinsic}. In essence, these methods, learn from the data the few leading eigenfunctions of
the Laplace–Beltrami operator that captures the Riemannian geometry of the observations, and consider them as a parametrisation of the manifold representing the invariant density.


\paragraph{$\triangleright$~Modelling SDEs from population level snapshots/boundary conditions.}
With recent computational advances in optimal transport, a growing body of work focuses on the implementation of Schr\"odinger bridge sampling methods, including formulations with additional path constraints. These mostly generative methods aim to transport the data distribution from some initial boundary condition to a terminal one, typically by learning the underlying stochastic equation to perform this transport through Schr\"odinger bridge sampling~\citep{lipman2022flow, pooladian2023multisample, albergo2023stochastic, albergo2022building, zhang2026learning}.
Flow matching~\citep{lipman2022flow} identifies the probability flow ODE that pushes forward an initial Gaussian density to a target one by solving a regression problem. The method relies on analytically tractable probability paths that provide closed-form regression targets for learning the velocity field, resulting in simulation-free training of deterministic flows. However, the framework is restricted to Gaussian distributions since the employed objective becomes intractable for general source distributions.
Conditional flow matching (CFM)~\citep{tong2023improving} generalizes flow matching by introducing conditional probability paths between paired samples, allowing the marginal velocity field to be learned with regression without requiring explicit evaluation of the marginal densities or restrictive assumptions on the source distribution.
Generalized Schr\"odinger Bridge Matching (GSBM)~\citep{liu2023generalized} follows an alternating optimisation scheme that learns both drift and marginals. Given prescribed boundary conditions for initial and terminal densities, the framework minimises a kinetic energy term, and formulates the resulting problem in terms of a stochastic optimal control problem conditioned on the boundary conditions and a path cost that accounts for additional constraints. 
In the multi-marginal setting,~\citet{shen2025multimarginal} consider the problem of inferring unobserved population trajectories from marginal samples observed at several time points. Their method alternates between learning piecewise Schr\"odinger bridges between adjacent population snapshots and refining a reference SDE within a selected parametric family using trajectories sampled from the current bridge estimate. This addresses the limitation of standard multi-marginal Schr\"odinger bridges, that although they enforce multiple marginal constraints, they can reduce to independent adjacent-time bridge problems and therefore fail to identify the long-range temporal structure of the underlying dynamics. The proposed iterative approach introduces a global dynamical bias across time, but the method remains a snapshot-based bridge construction whose output is a path measure compatible with prescribed marginals.
Action matching~\citep{neklyudov2023action} introduces a simulation-free variational objective that identifies a time-dependent scalar potential (entropic action) $s_t$, whose gradient $\nabla s_t$ transports the densities from the initial to the boundary condition through the continuity equation. In its entropic formulation the $\nabla s_t$ can be considered as the drift of the underlying SDE, whose marginals match the boundary conditions. However, by construction, the framework can recover only gradient drifts and is therefore not suitable for identifying general stochastic systems with stationary probability currents. In contrast, simulation-free score and flow matching ([SF]$^2M$)~\citep{tong2023simulation} jointly learns the probability-flow ODE and the score function by regressing against closed-form quantities derived from conditional Brownian bridge paths, facilitating simulation-free identification of general Schr\"odinger bridge dynamics with non-gradient drifts. ~\citet{zhang2026learning} extend this simulation-free Schr\"odinger bridge perspective to non-equilibrium reference dynamics, by replacing the Brownian reference process with a multivariate Ornstein–Uhlenbeck process. This provides the added flexibility of employing reference dynamics of non gradient-type with steady state probability currents, extending thereby the simulation-free Schr\"odinger bridge matching framework to non-conservative, non-equilibrium stochastic systems.



\paragraph{Geometry aware generative methods.}

Metric flow matching (MFM) generalizes CFM by learning
interpolants that account for the geometry of the
data. However, MFM does not assume a stochastic underlying process, as our framework does, only a deterministic interpolation (transport) that respects the data manifold. 
However, by assuming a specific noise amplitude for the underlying SDE, one can consider the flow field as generated by the effective drift of a probability flow ODE associated with the considered SDE and make inferences about the underlying drift function. This is the approach we followed when comparing the performance of MFM to our framework in Table~\ref{tab:LC}.

\paragraph{Approximating observation geometry in the ambient space.}
In our work, we approximate the geometry induced by the observations by endowing the ambient space $\mathbb{R}^d$ with an observation-dependent Riemannian metric $H(\mathbf{x})$ (Eq.~\ref{eq:metric_approx}) that encodes the local anisotropy of the data distribution. In our framework this metric acts as a constraint for data-augmentation and as a geometric inductive bias for drift function inference: augmented paths are encouraged to remain in regions where the metric $H(\mathbf{x})$ induces smaller distances, i.e. in the vicinity of geodesics computed with respect to this metric, thereby aligning the augmented paths with the empirical observation geometry.

This perspective connects to a growing body of work that \textbf{approximates Riemannian metrics directly in the ambient space} as a proxy for the unknown curved low-dimensional data manifold, instead of first estimating its intrinsic dimensionality and then constructing explicit low-dimensional embeddings. 

In  parallel, an increasing body of literature focuses on 
endowing generative models with geometric constraints or inductive biases. While most methods
function in an autoencoder-like setting, by learning an embedding function for projecting to a lower dimensional space that respects prescribed or learned geometric constraints~\citep{duque2022geometry, kalatzis2020variational, arvanitidis2017latent} geometry, "Riemannian" methods, similar to our proposed method, operate in the ambient space by directly a Riemannian geometry embedded there
and define normalizing flows or other generative processes directly on the
manifold of interest.
~\cite{mathieu2020riemannian} introduce a framework for continuous normalizing flows  defined in the ambient space, respecting a prescribed Riemannian geometry. Similarly, ~\cite{de2022riemannian} proposed a score-based generative model that models target densities with support on prescribed Riemannian manifolds in terms of a time-reversal of Langevin dynamics.

Metric flow matching~\citep{kapusniak2024metric} interpolates data distributions that respect the geodesic interpolants computed according to the metric induced by the observations. The method employs a data-adapted metric in the ambient space to design interpolants (geodesic curves) with low kinetic energy under the approximated geometry, and constrains the generative paths to respect manifold induced by the data samples. Our construction is conceptually similar with these approaches, in that we also avoid explicit low-dimensional embeddings and instead approximate the observation manifold through a Riemannian metric living in the ambient space. However, in contrast to methods focused on deterministic transport or simulation-free matching, we use the learned metric to regularise continuous-time diffusion bridges and drift inference, through the stochastic controlled geometric augmentation, so that the recovered stochastic dynamics are geometrically consistent with the geometry of the observation-induced invariant measure.

\paragraph{Positioning of the present work.}
Our approach combines the nonparametric flexibility of Gaussian-process–based drift inference from time-series data with recent geometric ideas for population-level SDE modelling. Similar to Metric Flow Matching~\citep{kapusniak2024metric}, we posit that augmented trajectories should remain on the manifold induced by the observations: both frameworks estimate a data-adapted Riemannian metric and construct interpolants (geodesics and bridges) that respect this geometry. MFM learns the underlying ODE necessary to transport an initial distribution to a target one under the data-adapted metric, while our framework assumes underlying stochastic dynamics. Nevertheless, once the diffusion is known or coarsely estimated, one can interpret the inferred ODE as a probability flow ODE and make inferences about the underlying drift function of a stochastic system.  The GSBM framework~\citep{liu2023generalized} employs a stochastic control objective that is similar to the objective we consider for constructing the augmented paths. 
In its standard formulation, GSBM does not impose the geometric constraints on augmented paths that we propose in our framework. For the baseline comparison presented here, we therefore employed a modified variant of GSBM in which we consider the same geometric path constraint term used in our method to the GSBM loss, while keeping all other components of GSBM unchanged.
 Finally, whereas these methods typically learn a drift that transports a single source distribution to a single terminal snapshot, yielding thus a \textbf{locally valid dynamics}, our method, akin to multi-marginal bridge sampling~\citep{shen2024multi}, fits a sequence of bridges across multiple time points to recover a \textbf{single global drift} consistent with the underlying drift dynamics. 


\section{Geometric constraints on inference}
Our method bridges the gap between approaches that rely only on the temporal structure of observations and those that approximate the invariant density, while ignoring temporal order. Motivated by advances in geometric statistics ~\citep{miolane2020geomstats, sommer2020probabilistic}, and the growing interest on the concept of manifold hypothesis~\citep{fefferman2016testing,shnitzer2020manifold}, i.e., the consideration that the state of multi-dimensional dynamical systems often resides in low-dimensional regions of the state space,
several recent methods integrate geometric and temporal constraints in stochastic system identification. In \emph{Langevin regression} framework~\citep{callaham2021nonlinear}, the Kramers-Moyal (KM) coefficients are estimated and low sampling effects are accounted for by solving an adjoint Fokker-Planck equation, with regularisation via moment matching~\citep{lade2009finite}.\cite{tong2020trajectorynet} consider the manifold of the observations for inference of cellular dynamics. Their method employs dynamic optimal transport to interpolate between measured distributions constrained to lie in the vicinity of the observations. While sharing similar intuitions with our method, Tong et al. do not employ SDE modelling for inherently stochastic cellular dynamics and do not consider the underlying geometry of the observations, relying solely on constraints penalizing pairwise distances between them. Shnitzer et al.~\citep{shnitzer2020manifold, shnitzer2016manifold}
employ diffusion maps to approximate the eigenfunctions of the backward Kolmogorov operator (the generator of the stochastic Koopman operator~\citep{giannakis2019data,vcrnjaric2020koopman}). By evolving the dominant operator eigenspectrum with a Kalman filter, they account for the temporal order of observations. However, their approach is limited to conservative systems and requires the presence of a spectral gap in the approximated operator's spectrum.

\section{Theoretical justification for Riemannian manifold approximation of the invariant density}\label{appsec:theory_low_dim}

\input{geometry}

To demonstrate the intuition of behind the low dimensional manifold approximation mentioned in the main text, we embedded the 2 dimensional Van der Pol drift into an D-dimensional ($D \in [3,10]$ space via a linear projection, add a stable linear force in the subspace orthogonal to the embedding, adding also D-dimensional isotropic noise, and simulated the resulting system. The simulated systems were
full D-dimensional stochastic dynamical systems. Subsequently, we estimated the effective dimensionality of the observation manifold by computing the participation ratio of the covariance spectrum of the observations. Let $\lambda_1,\ldots,\lambda_D$ denote the eigenvalues of the empirical covariance matrix of the observed states. The participation ratio is defined as
\begin{equation}
\mathrm{PR}
=
\frac{\left(\sum_{i=1}^{D} \lambda_i\right)^2}
{\sum_{i=1}^{D} \lambda_i^2}.
\end{equation}

\begin{figure}
 \includegraphics[width=0.96\textwidth]{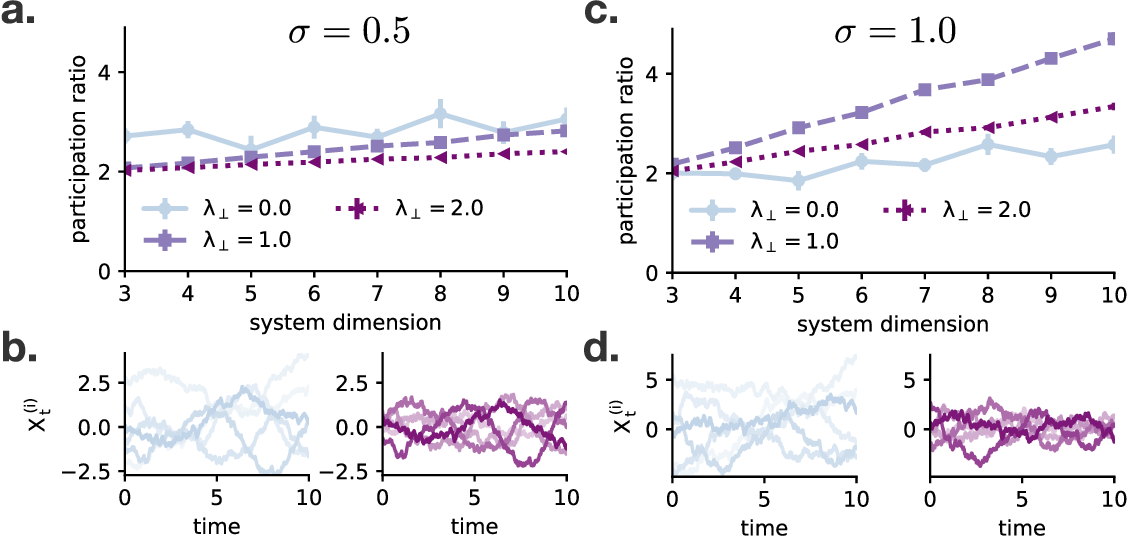}
\caption{\textbf{Effective dimensionality of the invariant density as a function of ambient system dimension for a $D$-dimensional system obtained by linearly embedding a 2-dimensional Van der Pol dynamics into spaces of dimension $D \in \{3,\dots,10\}$, with isotropic $D$-dimensional noise and an additional restoring drift acting in directions perpendicular 
to the embedded subspace.} Effective dimensionality is quantified by the participation ratio for noise amplitudes \textbf{a.)} $\sigma=0.5$ and \textbf{c.)} $\sigma=1.0$. Indicated values are averages over 10 independent runs for each setting. Different curves correspond to different values of $\lambda_{\perp}$, which controls the strength of the restoring drift in directions orthogonal to the embedded 2-dim. subspace: $\lambda_{\perp}=0$ corresponds to the absence of perpendicular confinement, whereas larger values of $\lambda_{\perp}$ increasingly pull trajectories toward the embedded manifold.
\textbf{b., d.)} Example sample trajectories for systems of dimension $D=6$, shown for two values of $\lambda_{\perp}=0.0$ (\emph{left}) and $\lambda_{\perp}=2.0$ (\emph{right}).}
\end{figure}

\newpage

\section{Theoretical justification of geometric augmentation for large inter-observation intervals}

In the following sections we provide a theoretical analysis justifying our choice to employ geometric path augmentation to improve inference in the large inter-observation limit. In particular, in Sec~\ref{app:deteriorates}, we revisit the fact that inference starting from the Euler-Maruyama discretisation deteriorates for increasing inter-observation interval. Then we study the terms in the remainder of the discretisation that become important when the time step (or inter-observation interval) is large, and connect these terms with the geometry of the unknown vector field. We show that for non-linear systems the remainder contains terms related to the curvature of the flow, and that neglecting these terms amounts to assuming a vector field with straight flow-lines in-between observations. This introduces a bias in inference that is linear in the step size. By approximating the curvature by means of controlled path augmentation with reference the geodesic curves of the invariant manifold, our method partially accounts for these remainder terms.

\subsection{Inference performance deteriorates with increasing inter-observation interval for existing frameworks } \label{app:deteriorates}

\begin{SCfigure}[50][tbh]
  \includegraphics[width=0.6\textwidth]{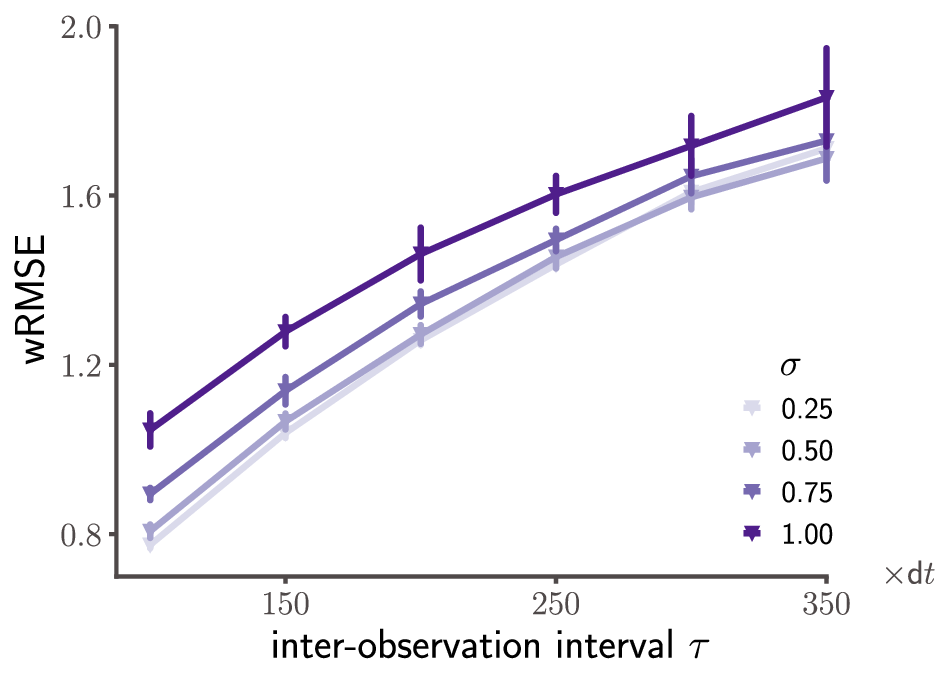}
  \caption{\textbf{Increasing observation interval between successive observations $\tau$ deteriorates performance quantified by increasing weighted root mean squared error (wRMSE) for Gaussian process-based inference.} 
  Weighted root mean square error between estimated and ground truth drift vector fields for increasing observation interval $\tau$ between subsequent observations for different noise conditions (indicated by different hues). Observations were collected from a Van der Pol oscillator system simulated with $dt=0.01$ for $T=500$ time units. 
Error bars indicate one standard deviation over ten independent realizations.}
  \label{fig:increase} 
\end{SCfigure}

We computed the weighted root mean square error (wRMSE) between ground truth flow fields and estimated ones for several commonly applied inference frameworks. We observed that the performance of all of them deteriorates once the inter-observation interval becomes large.

We started with the method that motivated our research, approximating drift functions through Gaussian processes, the method outlined in~\cite{ruttor2013approximate}. The method approximates the drift functions with Gaussian process regression, using the system state $\mathbf{X}_t$ as the regressor and state increments as the response variable $\mathbf{Y}_t \dot{=} \frac{\mathbf{X}_{t+\tau}-\mathbf{X}_t}{\tau}$. This is the Bayesian counterpart of earlier methods encountered in physics literature~\citep{friedrich1997description, ragwitz2001indispensable}, providing additionally uncertainty estimation through the Gaussian process approximation. 

As is evident from Figure~\ref{fig:increase} the discrepancy between ground truth and estimated vector fields  increases for increasing temporal distance between successive observations. This should be understood, under
the consideration that inference of the drift based on regression on state increments results from an approximation relying on a truncated Ito-Taylor expansion. This is also the starting point of the Euler Maruyama discretisation. As the time interval between successive steps of this approximation increases, the truncated approximation does not longer hold, and higher order terms should be considered.

\newpage
\subsection{Inference based on Euler-Maruyama discretisation does not account for the curvature of the trajectories in the state space}\label{app:curvature}
To be more precise, a general SDE of the form
\begin{equation}\label{appeq:generic_sde}
 \text{d}\mathbf{X}_t = \mathbf{f}(\mathbf{X}_t,t)\text{d}t + \boldsymbol{\sigma}(\mathbf{X}_t, t) \text{d}\mathbf{W}_t.
\end{equation}
is a shorthand for the integral equation 
\begin{equation}
    \mathbf{X}_t = \mathbf{X}_{t_0} + \int_{t_0}^t \mathbf{f}(\mathbf{X}_s, s)\,\de s + \int_{t_0}^t \boldsymbol{\sigma}(\mathbf{X}_s, s)\,\de \mathbf{W}_s,
\end{equation}
where as previously in this manuscript, we consider the stochastic integrals in the \textbf{Itô sense}. (Here we start from a more general formulation of the stochastic equation with both diffusion and drift terms being state- and time-dependent to highlight that also for more general SDEs our geometric argument is valid.)

Applying the Itô formula on each integrand, and integrating from $t_0$ to $t$, we obtain the Itô-Taylor expansion of Eq.~\ref{appeq:generic_sde}
\begin{align}
\mathbf{f}(\mathbf{X}_t, t)
=& \mathbf{f}(\mathbf{X}_{t_0}, t_0) 
+ \int_{t_0}^t \frac{\partial \mathbf{f}(\mathbf{X}_{s}, s)}{\partial s} \, \de s 
+ \int_{t_0}^t \sum_u \frac{\partial \mathbf{f}(\mathbf{X}_{s}, s)}{\partial X^{(u)}} f_u(\mathbf{X}_{s}, s) \, \de s 
\nonumber
\\
&\,\,
+ \int_{t_0}^t \sum_u \frac{\partial \mathbf{f}(\mathbf{X}_{s}, s)}{\partial X^{(u)}} 
\left[\boldsymbol{\sigma}(\mathbf{X}_{s}, s)\, \de \mathbf{W}_s\right]_u 
+ \int_{t_0}^t \frac{1}{2} \sum_{u,v} \frac{\partial^2 \mathbf{f}(\mathbf{X}_{s}, s)}{\partial X^{(u)} \partial X^{(v)}} 
\left[\boldsymbol{\sigma}(\mathbf{X}_{s}, s)\,  \boldsymbol{\sigma}^\top(\mathbf{X}_{s}, s)\right]_{uv} \de s
\nonumber \\
 =&  \mathbf{f}(\mathbf{X}_{t_0}, t_0)
+ \int_{t_0}^{t} \mathcal{L}^\dagger_s \mathbf{f}(\mathbf{X}_{s}, s) \, \de s
+ \sum_\nu \int_{t_0}^{t} \mathcal{L}_{W,\nu} \mathbf{f}(\mathbf{X}_{s}, s) \, \de W^{(\nu)}_s,
\end{align}
 and
\begin{align}
\boldsymbol{\sigma}(\mathbf{X}_t, t)
=& \boldsymbol{\sigma}(\mathbf{X}_{t_0}, t_0)
+ \int_{t_0}^t \frac{\partial \boldsymbol{\sigma}(\mathbf{X}_{s}, s)}{\partial s} \, \de s 
+ \int_{t_0}^t \sum_u \frac{\partial \boldsymbol{\sigma}(\mathbf{X}_{s}, s)}{\partial X^{(u)}} f_u(\mathbf{X}_{s}, s) \, \de s \nonumber \\
&\,\,+ \int_{t_0}^t \sum_u \frac{\partial \boldsymbol{\sigma}(\mathbf{X}_{s}, s)}{\partial X^{(u)}} 
\left[\boldsymbol{\sigma}(\mathbf{X}_{s}, s)\, \de \mathbf{W}_s\right]_u  + \int_{t_0}^t \frac{1}{2} \sum_{u,v} \frac{\partial^2 \boldsymbol{\sigma}(\mathbf{X}_{s}, s)}{\partial X^{(u)} \partial X^{(v)}} 
\left[\boldsymbol{\sigma}(\mathbf{X}_{s}, s)\,  \boldsymbol{\sigma}^\top(\mathbf{X}_{s}, s)\right]_{uv} \de s \nonumber \\
=& \boldsymbol{\sigma}(\mathbf{X}_{t_0}, t_0)
+ \int_{t_0}^{t} \mathcal{L}^\dagger_s \boldsymbol{\sigma}(\mathbf{X}_{s}, s) \, \de s
+ \sum_\nu \int_{t_0}^{t} \mathcal{L}_{W,\nu} \boldsymbol{\sigma}(\mathbf{X}_{s}, s) \, \de W^{(\nu)}_s,
\end{align}
where we have used the fact that the product of stochastic differentials due to the Ito isometry and multiplication rules equals the noise covariance times the time step $$dX^{(u)}_t dX^{(v)}_t = \left[ \boldsymbol{\sigma} \boldsymbol{\sigma}^\top \right]_{uv}\,dt,$$ where $${ dX^{(u)}_s = f_u\,ds + \sum_{j=1}^m \sigma_{uj}\,dW^{(j)}_s },$$ while the superscripts/subscripts $u,v$ indicate dimensional components.

In the above equations, we have introduced the operators acting on test-functions $\mathbf{h}:\mathbb{R}^D \rightarrow \mathbb{R}^D$
\begin{equation}
\mathcal{L}^\dagger_t \mathbf{h} = 
\frac{\partial \mathbf{h}}{\partial t}
+ \sum_u \frac{\partial \mathbf{h}}{\partial X^{(u)}} f_u
+ \frac{1}{2} \sum_{u,v} \frac{\partial^2 \mathbf{h}}{\partial X^{(u)} \partial X^{(v)}} 
\left[\boldsymbol{\sigma}(\mathbf{X}_{s}, s) \,  \boldsymbol{\sigma}^\top (\mathbf{X}_{s}, s) \right]_{uv}
\end{equation}
and 
\begin{equation}
\mathcal{L}_{W, v} \mathbf{h} = \sum_u \frac{\partial \mathbf{h}}{\partial X^{(u)}} \, \boldsymbol{\sigma}_{uv}, \qquad \text{for\,\,}v = 1, \ldots, n.
\end{equation}

With these expressions, the original integral equation for $\mathbf{X}_t$ can be written as
\begin{align}
\mathbf{X}_t 
&= \mathbf{X}_{t_0}
+ \mathbf{f}(\mathbf{X}_{t_0}, t_0)(t - t_0)
+ \boldsymbol{\sigma}(\mathbf{X}_{t_0}, t_0)\left( \mathbf{W}_t - \mathbf{W}_{t_0} \right) + \\
\quad \color{orange}{R_1=}& 
\left\{
\begin{aligned}
& + \int_{t_0}^t \int_{t_0}^{s} \mathcal{L}^\dagger_u \mathbf{f}(\mathbf{X}_u, u)\, \de u\, \de s
+ \sum_\nu \int_{t_0}^t \int_{t_0}^{s} \mathcal{L}_{W,\nu} \mathbf{f}(\mathbf{X}_u, u)\, \de W^{(\nu)}_u\, \de s \nonumber \\
&+\int_{t_0}^t \int_{t_0}^{s} \mathcal{L}^\dagger_u \boldsymbol{\sigma}(\mathbf{X}_u, u)\, \de u\, \de \mathbf{W}_s + \sum_\nu \int_{t_0}^t \int_{t_0}^{s} \mathcal{L}_{W,\nu} \boldsymbol{\sigma}(\mathbf{X}_u, u)\, \de W^{(\nu)}_u\, \de \mathbf{W}_s.
\end{aligned}
\right.
\end{align}

In the last equation, dropping the terms in the remainder $R_1$ results in the Euler–Maruyama integration scheme~\citep{jentzen2011taylor}. Introducing the discrete time and noise increments
\begin{equation}
\Delta t_n = t_{n+1} - t_n = \int_{t_n}^{t_{n+1}} \de s, 
\quad 
\Delta \mathbf{W}_n = \mathbf{W}_{t_{n+1}} - \mathbf{W}_{t_n} = \int_{t_n}^{t_{n+1}} \de \mathbf{W}_s,
\end{equation}
we result in the discrete time equation commonly used for numerical integration of SDEs
\begin{equation}
\mathbf{X}_{n+1} = \mathbf{X}_n + \mathbf{f}( \mathbf{X}_n,t_n)\, \Delta t_n + \boldsymbol{\sigma}\, \Delta \mathbf{W}_n.
\end{equation}
 This is also the starting point of most inference methods that employ
the regression scheme mentioned above by approximating the drift as
\begin{equation}
\hat{\mathbf{f}}( \mathbf{X}_n,t_n) \approx \frac{\mathbf{X}_{n+1} - \mathbf{X}_n }{\Delta t}\sim \mathcal{N}\left(\mathbf{\hat{f}}(\mathbf{X}_n,t_n), \frac{ \boldsymbol{\sigma}\,\boldsymbol{\sigma}^{\top} }{\Delta t }\right).
\end{equation}
This discretisation is a zero-order approximation of the ground truth dynamics, and assumes that $\textbf{f}(\cdot)$ remains constant throughout the interval $\Delta t$, i.e. throughout the inter-observation interval $\tau$ in the inference setting.   
However as $\tau$ increases, higher-order terms in the remainder $R_1$ of the Itô-Taylor expansion become significant, since the assumption that the drift is approximately constant over $\tau$ does not hold.

We can glean onto the terms that become important once the inter-observation interval becomes large, by applying the Itô formula onto each one of the integrands in $R_1$ separately \textbf{for the specific setting we consider in this manuscript}, i.e. that of time-independent drift function $\mathbf{f}(\mathbf{x})$ and constant diffusion matrix $\boldsymbol{\sigma}$.
In the following, we demonstrate that the leading-order error of this approximation is governed by the intrinsic geometry of the drift vector field. This provides further insight and a geometric explanation for the deterioration of inference methods for increasing 
inter-observation interval $\tau$.

 In short we show that, inference methods based on the Euler-Maruyama discretisation-based inference effectively assume that the vector field between consecutive observations $\mathbf{X}_n$ and $\mathbf{X}_{n+1}$ does not change. Our analysis shows this is equivalent to assuming trajectories are straight lines ($\mathbf{J}_f \mathbf{f} \parallel \mathbf{f}$) and the It\^o correction is constant. In reality, trajectories curve ($\mathbf{J}_f \mathbf{f}$ has also a perpendicular component), and this curvature itself changes along the vector field. The Euler-Maruyama discretisation-based inference scheme systematically misses these higher-order geometric features, leading to biased drift estimates.

\subsubsection{First remainder term $R_{1,a}$}
We denote the first term of the reminder by $R_{1,a}$
\begin{equation}
R_{1,a} = \int_{t_0}^t \int_{t_0}^s \mathcal{L}^\dagger_u \mathbf{f}(\mathbf{X}_u)\, \de u\, \de s.
\end{equation}

Applying Itô's formula to the integrand  $\mathcal{L}^\dagger_t \mathbf{f}(\mathbf{X}_u, u) $, we get
\begin{align}
\de \mathcal{L}^\dagger_u \mathbf{f}(\mathbf{X}_u)
&= {\frac{\partial}{\partial u} \mathcal{L}^\dagger_u \mathbf{f}(\mathbf{X}_u)\, \de u }
+ \sum_{j=1}^d \frac{\partial \mathcal{L}^\dagger_u \mathbf{f}}{\partial X^{(j)}}(\mathbf{X}_u)\, \de X_u^{(j)} + \frac{1}{2} \sum_{j,k=1}^d \frac{\partial^2 \mathcal{L}^\dagger_u \mathbf{f}}{\partial X^{(j)} \partial X^{(k)}}(\mathbf{X}_u)\,
\left[ \boldsymbol{\sigma} \boldsymbol{\sigma}^\top \right]_{jk}\, \de u.
\end{align}

Plugging in the original equation $ \de X_u^{(j)} = f_j\, \de u + \sum_{\nu=1}^m \sigma_{j\nu}\, \de W_u^{(\nu)} $, and integrating over the time from $t_0$ to $u$
\begin{align}
\mathcal{L}^\dagger_u \mathbf{f}(\mathbf{X}_u) =& \;\mathcal{L}^\dagger_{t_0} \mathbf{f}(\mathbf{X}_{t_0}) + \int_{t_0}^{u} \left( \frac{\partial}{\partial w} (\mathcal{L}^\dagger_w \mathbf{f}(\mathbf{X}_w)) + \sum_{j} \frac{\partial (\mathcal{L}^\dagger_w \mathbf{f})}{\partial X^{(j)}} f_j + \frac{1}{2} \sum_{j,k} \frac{\partial^2 (\mathcal{L}^\dagger_w \mathbf{f})}{\partial X^{(j)} \partial X^{(k)}} [\boldsymbol{\sigma} \boldsymbol{\sigma}^\top]_{jk} \right) \, \de w \nonumber \\ &\,+ \int_{t_0}^{u} \sum_{j} \frac{\partial (\mathcal{L}^\dagger_w \mathbf{f})}{\partial X^{(j)}} [\boldsymbol{\sigma} \de \mathbf{W}_w]_j .
\end{align}

Applying Fubini’s theorem in the original double integral, we change the order of integration
\begin{equation}
\int_{t_0}^t \int_{t_0}^s \phi(u)\, \de u\, \de s = \int_{t_0}^t (t - u)\, \phi(u)\, \de u,
\end{equation}
and applying it twice we obtain
\begin{align}
R_{1,a} = \int_{t_0}^t \int_{t_0}^s \mathcal{L}^\dagger_u \mathbf{f}(\mathbf{X}_u)\, \de u\, \de s
&= \int_{t_0}^t \frac{(t - u)^2}{2}\left[ 
\underbrace{ \sum_j \frac{\partial \mathcal{L}^\dagger_u \mathbf{f}}{\partial X^{(j)}} f_j }_{\textcolor{mymaroon}{R^1_{1,a}}}
+ \underbrace{\frac{1}{2} \sum_{j,k} \frac{\partial^2 \mathcal{L}^\dagger_u \mathbf{f}}{\partial X^{(j)} \partial X^{(k)}} 
[\boldsymbol{\sigma} \boldsymbol{\sigma}^\top]_{jk} }_{\textcolor{mymaroon}{R^2_{1,a}}}
\right] \de u \nonumber \\
&\quad + \int_{t_0}^t \frac{(t - u)^2}{2} \underbrace{\sum_j \frac{\partial \mathcal{L}^\dagger_u \mathbf{f}}{\partial X^{(j)}} \left[ \boldsymbol{\sigma}\, \de \mathbf{W}_u \right]_j }_{\textcolor{mymaroon}{R^3_{1,a}}} + \frac{\tau^2}{2}\mathcal{L}^\dagger_t \mathbf{f}(\mathbf{X}_{t_0}).
\end{align}
In the previous equation we have dropped the term $\frac{\partial}{\partial w}\!\left(\mathcal{L}^\dagger_{w}\,\mathbf{f}(\mathbf{X}_{w})\right) $ that is equal to zero and that would require the drift $\mathbf{f}$ to be time-dependent to be non-negligible.

\paragraph{First component $R^{\,1}_{1,a}$ of remainder term $R_{1,a}$: Flow curvature term.}
The Backward Kolmogorov generator applied to a vector field $\mathbf{f}$ can be written as
\begin{equation}\label{appeq:jacobian_expansion}
\mathcal{L}^\dagger \mathbf{f} \;=\; \mathbf{J}_f\,\mathbf{f} \;+\; \frac12\,\Delta_D \mathbf{f}
 .
\end{equation}
In Eq.~\ref{appeq:jacobian_expansion}, $\mathbf{J}_f \,\dot{=} \,{\nabla} \mathbf{f}$ denotes the Jacobian of $\mathbf{f}$, $\mathbf{D}\dot{=}\,\boldsymbol{\sigma} \boldsymbol{\sigma}^\top $ the noise covariance, and ${\Delta_\mathbf{D} \, \dot{=} \,\sum_{j,k} \mathbf{D}_{jk}\,\partial_{X^{(j)} X^{(k)}}^2 }$ is the noise-covariance weighted Laplacian operator.
Thus each component of $\mathcal{L}^\dagger \mathbf{f}$ comprises the directional derivative of the drift $\mathbf{J}_f \mathbf{f}$ plus an anisotropic/noise-covariance weighted Laplacian of $\mathbf{f}$, which in component-wise form is expressed as
\begin{equation}
\big[\mathcal{L}^\dagger\mathbf{f}\big]_{i}
= \sum_{k}\frac{\partial f_{i}}{\partial X^{(k)}}\,f_{k}
\;+\;\frac12\sum_{k,\ell} \mathbf{D}_{k\ell}\,
\frac{\partial^{2} f_{i}}{\partial X^{(k)}\partial X^{(\ell)}} .
\end{equation}

Differentiating wrt to $X^{(j)}$ yields
\begin{equation}\label{appeq:72}
\frac{\partial}{\partial X^{(j)}}\big[\mathcal{L}^\dagger\mathbf{f}\big]_{i}
= \sum_{k}\frac{\partial^{2} f_{i}}{\partial X^{(j)}\partial X^{(k)}}\,f_{k}
\;+\;\sum_{k}\frac{\partial f_{i}}{\partial X^{(k)}}\,
\frac{\partial f_{k}}{\partial X^{(j)}}
\;+\;\frac12\sum_{k,\ell} \mathbf{D}_{k\ell}\,
\frac{\partial^{3} f_{i}}{\partial X^{(j)}\partial X^{(k)}\partial X^{(\ell)}} ,
\end{equation}
and thus we rewrite the $i$-th component of the term $R_{1,a}^{\,1}$ as
\begin{equation}\label{appeq:73}
\big[R_{1,a}^{\,1}\big]_{i}
= \int_{t_{0}}^{t} \frac{(t - u)^2}{2}\!
\left[
\sum_{j,k}\frac{\partial^{2} f_{i}}{\partial X^{(j)}\partial X^{(k)}}\,f_{k}\,f_{j}
+\sum_{j,k}\frac{\partial f_{i}}{\partial X^{(k)}}
\,\frac{\partial f_{k}}{\partial X^{(j)}}\,f_{j}
+\frac12\sum_{j,k,\ell} \mathbf{D}_{k\ell}\,
\frac{\partial^{3} f_{i}}{\partial X^{(j)}\partial X^{(k)}\partial X^{(\ell)}}\,f_{j}
\right]_{\!i}\,\de u . 
\end{equation}

The third-order state-derivative in the last summand of Eq.~\ref{appeq:73}, indicates that this last term is inactive for linear or quadratic drift functions $\mathbf{f}$.

We re-write again this part of the remainder in a more compact vector notation in terms of the directional derivative of $(\mathbf{J}_{f}\mathbf{f})$ and $\tfrac12\,\Delta_{D}\mathbf{f}$ along the vector field as
\begin{equation}\label{appeq:remainder_maroon}
{\color{mymaroon}
{R_{1,a}^{\,1}
=\int_{t_{0}}^{t} \frac{(t - u)^2}{2}\,
\big[\underbrace{{\nabla}(\mathbf{J}_{f}\mathbf{f})\cdot\mathbf{f}}_{\text{flow curvature}} \;+\; \underbrace{ {\nabla}\!\big(\tfrac12\,\Delta_{D}\mathbf{f}\big)\cdot\mathbf{f} }_{\text{diffusive term along the flow}}\big]\; \de u }}.
\end{equation}

This part of the remainder captures two geometric effects that standard inference methods neglect: the \textbf{intrinsic curvature of deterministic flow trajectories in state space}, and the \textbf{systematic bias introduced by the spatial variation of both drift and diffusion} along these trajectories, when both drift and diffusion are assumed as constant between inter-observation intervals.
\begin{itemize}
\item To understand the {\color{mymaroon}\textbf{first term}}, $\nabla(\mathbf{J}_f\mathbf{f}) \cdot \mathbf{f}$, from a geometric perspective, let us consider a deterministic dynamical system with dynamics $\dot{\mathbf{x}}_t = \mathbf{f}(\mathbf{x}_t)$. A trajectory initiated from an initial condition $\mathbf{x}_0$ traces a streamline in the state space $\mathbb{R}^d$. We express the acceleration of this trajectory in terms of the directional derivative
\begin{equation}
\ddot{\mathbf{x}}_t = \frac{\de}{\de t}\mathbf{f}(\mathbf{x}_t) 
= \mathbf{J}_f(\mathbf{x}_t) \cdot \mathbf{f}(\mathbf{x}_t)
= \mathbf{J}_f\cdot\mathbf{f}.
\end{equation}

The acceleration vector admits a natural orthogonal decomposition
 comprising a component parallel to the vector field $\mathbf{f}$ and an orthogonal component to $\mathbf{f}$
\begin{equation}
\mathbf{J}_f\cdot\mathbf{f} 
= P_{\parallel}(\mathbf{f}) \, \mathbf{J}_f\cdot\mathbf{f}
+ P_{\perp}(\mathbf{f}) \, \mathbf{J}_f \cdot\mathbf{f}.
\end{equation}
Here $P_{\parallel}(\mathbf{f}(\mathbf{x})) = \frac{\mathbf{f}(\mathbf{x})\mathbf{f}^\top(\mathbf{x})}{\|\mathbf{f}(\mathbf{x})\|^2}$ and 
$P_{\perp}(\mathbf{f}(\mathbf{x})) = \mathbb{I} - P_{\parallel}(\mathbf{f}(\mathbf{x}))$ stand for the parallel and orthogonal projectors. The parallel component quantifies the rate of change of speed along the trajectory (tangential acceleration), whilst the perpendicular component defines the \textbf{curvature vector} $\kappa_{\text{flow}}(x)$~\citep{kuhnel2002differential}, which quantifies the bending of the trajectories
\begin{equation}
\label{eq:curvature_vector}
\boldsymbol{\kappa}_{\mathrm{flow}}(\mathbf{x}) = \frac{P_{\perp}(\mathbf{f}(\mathbf{x})) \mathbf{J}_f(\mathbf{x}) \mathbf{f}(\mathbf{x})}{\|\mathbf{f}(\mathbf{x})\|^2}.
\end{equation}
When $\boldsymbol{\kappa}_{\text{flow}} = 0$, the trajectories are straight lines in the state space, while
when ${\|\boldsymbol{\kappa}_{\text{flow}}\| > 0}$ they are curved.

The term $\nabla(\mathbf{J}_f \mathbf{f}) \cdot \mathbf{f}$ quantifies the \textbf{evolution of the trajectory curvature} \footnote{More precisely the directional derivative of the acceleration, $\mathbf{J}_f(\mathbf{x})\cdot \mathbf{f}$ along the flow direction, or the \textbf{rate at which the acceleration changes along the flow, or a measure of how the local curvature of $\mathbf{f}$ as a vector field influences trajectory evolution}.} as the system moves along the flow field. From Eq.~\ref{appeq:73} we have for each dimensional component $i$ of this term
\begin{equation}
\label{eq:nabla_jf_expansion}
\begin{split}
\left[\nabla(\mathbf{J}_f \mathbf{f}) \cdot \mathbf{f}\right]_i &= \sum_{j,k} \frac{\partial^2 f_i}{\partial X^{(j)} \partial X^{(k)}} f_k f_j + \sum_{j,k} \frac{\partial f_i}{\partial X^{(k)}} \frac{\partial f_k}{\partial X^{(j)}} f_j \\
&= [\mathbf{f}^\top (\nabla^2 f_i) \mathbf{f}] + [\mathbf{J}_f^2 \mathbf{f}]_i.
\end{split}
\end{equation}
We observe that this term captures the effects of how both second-order spatial variation of the flow field (the Hessian $\nabla^2 f_i$) and the Jacobian of the acceleration ($\mathbf{J}_f^2\mathbf f$) influence the evolution of trajectories. 
\begin{itemize}

\item In Eq.~\ref{eq:nabla_jf_expansion}, the \textbf{first sub-term}, $\mathbf{f}^\top (\nabla^2 f_i) \mathbf{f}$, represents the \textbf{second directional derivative} (or quadratic variation) of $f_i$ along the flow direction $\mathbf{f}$. It measures the curvature or second-order spatial variation of the $i$-th component of $\mathbf{f}$ in the direction $\mathbf{f}$.  In regions where the Hessian $\nabla^2 \mathbf{f}$ is large (as is for the case of a highly nonlinear drift with curving or bending behaviour), this term becomes significant, and it vanishes for linear or constant drift $\mathbf{f}$. Neglecting this term corresponds to approximating the flow by its linearisation. 

\item The \textbf{second sub-term}, $\mathbf{J}_f^2 \mathbf{f} = \mathbf{J}_f(\mathbf{J}_f \mathbf{f})$, of Eq.~\ref{eq:nabla_jf_expansion} represents the action of the Jacobian operator on the acceleration vector. Geometrically, it describes how the local linearised field acts on the acceleration as we move an infinitesimal step along the flow field, or in other words how the linear approximation changes when following the flow direction $\mathbf{f}$.

\end{itemize}

By temporal integration we have
\begin{equation}
\label{eq:r1a_flow_scaling}
R_{1,a}^{\,1}=\int_{t_0}^t \frac{(t - u)^2}{2} \Big( \nabla(\mathbf{J}_f \mathbf{f}) \cdot \mathbf{f} +{\nabla}\big(\tfrac12\,\Delta_{D}\mathbf{f}\big)\cdot\mathbf{f} \Big) \,\de u \sim \frac{\tau^2}{2} \Big( \nabla(\mathbf{J}_f \mathbf{f}) \cdot \mathbf{f}+ {\nabla}\big(\tfrac12\,\Delta_{D}\mathbf{f}\big)\cdot\mathbf{f} \Big),
\end{equation}
indicating that the evolution of trajectory curvature introduces an $O(\tau^3)$ correction to the transition density.

Drift inference based on Euler–Maruyama–type discretisation 
ignores between others the term $R_{1,a}^{\,1}$ introducing thereby a {mean} bias at each point $\mathbf{x}$ in the state space,
\begin{equation}
\label{eq:bias_main}
\beta^1_{1,a}(\mathbf{x})
= \frac{1}{\tau}\,R_{1,a}^{\,1}
\;\approx\; \frac{\tau}{2}\,\big[\nabla (\mathbf{J}_f \mathbf{f}) \cdot \mathbf{f} + {\nabla}\big(\tfrac12\,\Delta_{D}\mathbf{f}\big)\cdot\mathbf{f}\big](\mathbf{x}).
\end{equation}

This bias induces a mean error in drift estimate, when using Euler–Maruyama-based inference, leading to under- or over-estimation of the ground truth drift at state $\mathbf{x}$. This error scales as $\mathcal{O}(\tau^2)$ with the interval $\tau$.

Let us now consider the temporal rate of change experienced by a particle travelling along the flow field. The instantaneous speed of the particle at location $\mathbf{x}$ is $\|\mathbf{f}(\mathbf{x})\|$. The quantity in the brackets in Eq.~\ref{eq:bias_main},  $\nabla(\mathbf{J}_f \mathbf{f}) \cdot \mathbf{f} + \nabla\left(\frac{1}{2} \Delta_D \mathbf{f}\right) \cdot \mathbf{f}$, is a spatial derivative measuring how quickly the curvature and diffusion variation change as the particle moves in space. The rate of change of this variation per unit of time is expressed as
\begin{equation}
{\frac{\left\| \nabla(\mathbf{J}_f \mathbf{f}) \cdot \mathbf{f} + \nabla\left(\frac{1}{2} \Delta_D \mathbf{f}\right) \cdot \mathbf{f} \right\|(\mathbf{x})}{\|\mathbf{f}(\mathbf{x})\|}} \;\dot{=} \;\tau_{\mathrm{curv}}^{-2}(\mathbf{x}).
\end{equation}
In the last equation we have introduced the time scale of change $\tau_{\mathrm{curv}}$ as the inverse of the rate of change, which captures the characteristic time it takes for the curvature/diffusion variation to change significantly along the particles trajectory.
Then the relative magnitude error in the Euler-Maruyama-based drift estimate satisfies
\begin{equation}
\frac{\|\beta^1_{1,a}(\mathbf{x})\|}{\|\mathbf{f}(\mathbf{x})\|} = \frac{\tau^2}{6\, \tau^2_{\mathrm{curv}}(\mathbf{x})}  ,
\end{equation}
implying that the estimate is reliable only when the inter-observation interval $\tau \ll  \tau_{\mathrm{curv}}(\mathbf{x})$.

\item The {\color{mymaroon}\textbf{second term}} in Eq.\ref{appeq:remainder_maroon}, $\nabla(\frac{1}{2}\Delta_D \mathbf{f}) \cdot \mathbf{f}$, accounts for the diffusion part of the {backward generator} acting on the vector field $\mathbf{f}$. 
The anisotropic Laplacian $\Delta_D \mathbf{f}$ quantifies the \textbf{diffusion–weighted second-order spatial variation of the vector field}
\begin{equation}
[\Delta_D \mathbf{f}]_i = \sum_{j,k} D_{jk} \frac{\partial^2 f_i}{\partial X^{(j)} \partial X^{(k)}} =  \nabla\!\cdot\!\big(\mathbf{D}\,\nabla f_i\big).
\end{equation}

 The directional derivative quantifies how this term evolves along the flow field
\begin{equation}
\label{eq:diffusion_term}
\left[\nabla\!\left(\tfrac{1}{2}\Delta_D \mathbf{f}\right) \cdot \mathbf{f}\right]_i = \frac{1}{2}\sum_{j,k,\ell} D_{k\ell} \frac{\partial^3 f_i}{\partial X^{(j)} \partial X^{(k)} \partial X^{(\ell)}} f_j.
\end{equation}

This term captures how the diffusion-weighted spatial variation of the flow field varies across the state space. As trajectories traverse regions of varying drift curvature, the effective Itô correction itself changes, introducing systematic bias in inference methods that assume that drift is piece-wise constant in-between observations.

\end{itemize}
\if False

 \textbf{Deterministic curvature evolution}: how the bending of flow trajectories and the rate of change of acceleration evolve along the flow;

 \textbf{Stochastic curvature evolution}: how the diffusion-weighted curvature of the drift field varies spatially, including noise-induced drift when $D(\mathbf{x})$ is state-dependent.

Both effects scale as $O(\tau^2)$ and accumulate over the inter-observation interval. Neglecting $R_{1,a}^1$ introduces systematic bias into likelihood-based inference, particularly affecting parameters that govern the nonlinear structure of the drift (e.g., interaction strengths, feedback coefficients, threshold effects). When $\tau$ is large relative to $\tau_{\mathrm{drift}}$, the ratio
\begin{equation}
\frac{\|R_{1,a}^1\|}{\|\mathbf{f}\| \tau} \sim \frac{\tau}{\tau_{\mathrm{drift}}} \times \frac{\|\nabla(\mathbf{J}_f \mathbf{f}) \cdot \mathbf{f}\| + \|\nabla(\frac{1}{2}\Delta_D \mathbf{f}) \cdot \mathbf{f}\|}{\|\mathbf{f}\|/L}
\end{equation}
becomes $O(1)$, signalling that the curvature correction is comparable to the leading-order drift term. Ignoring this geometric structure leads to mis-estimation of drift parameters and underestimation of posterior uncertainty.

The remaining terms in the remainder $R_1$ (cf.\ equation~\eqref{eq:remainder_full}) involve higher-order diffusion corrections $R_{1,a}^2$ and stochastic contributions $R_{1,a}^3$, which we now examine.

\fi

\paragraph{Second component $R^{\,2}_{1,a}$ of remainder term $R_{1,a}$.}
The second component of the remainder term $R_{1,a}$ reads
\begin{equation}
R_{1,a}^{\,2}
= \int_{t_0}^{t} (t-u)\,\frac{1}{2}\sum_{j,k}
\frac{\partial^{2}\big(\mathcal L^\dagger_{u} \mathbf{f}\big)}{\partial X^{(j)}\partial X^{(k)}}
\,\big[\boldsymbol{\sigma\sigma}^{\top}\big]_{jk}\,\de u .
\end{equation}

For the $i$-th dimensional component we have

\begin{equation}
\begin{aligned}
\frac{\partial^{2}}{\partial X^{(h)}\partial X^{(j)}}\big[\mathcal L^\dagger_{u} f\big]_{i}
&= \sum_{k}\frac{\partial^{3} f_{i}}{\partial X^{(h)}\partial X^{(j)}\partial X^{(k)}}\,f_{k}
+ \sum_{k}\frac{\partial^{2} f_{i}}{\partial X^{(j)}\partial X^{(k)}}\,
\frac{\partial f_{k}}{\partial X^{(h)}} \\
&\quad+ \sum_{k}\frac{\partial^{2} f_{i}}{\partial X^{(h)}\partial X^{(k)}}\,
\frac{\partial f_{k}}{\partial X^{(j)}}
+ \sum_{k}\frac{\partial f_{i}}{\partial X^{(k)}}\,
\frac{\partial^{2} f_{k}}{\partial X^{(h)}\partial X^{(j)}} \\
&\quad+ \frac{1}{2}\sum_{k,\ell} \mathbf{D}_{k\ell}\,
\frac{\partial^{4} f_{i}}{\partial X^{(h)}\partial X^{(j)}\partial X^{(k)}\partial X^{(\ell)}} .
\end{aligned}
\end{equation}

For this remainder term, we have for each dimensional component $i$

\begin{equation}
{\ 
\big[R_{1,a}^{\,2}\big]_{i}
= \int_{t_0}^{t} (t-u)\,\frac{1}{2}\sum_{j,k} \mathbf{D}_{jk}\,
\left[
\frac{\partial^{2}}{\partial X^{(k)}\partial X^{(j)}}\big[\mathcal L^\dagger_{u} \mathbf{f}\big]_{i}
\right]\, \de u \ }.
\end{equation}

Geometrically, $R^{2}_{1,a}$ captures the \textbf{diffusion-weighted second-order spatial variation} of the generator $\mathcal{L}_u^\dagger \mathbf{f}$ across the $\sqrt{\tau}$-sized ellipsoid set by $\mathbf{D}$, i.e. the anisotropic Laplacian $\Delta_D(\mathcal{L}_u^\dagger \mathbf{f})$, the diffusion-weighted second spatial variation of the drift along the flow.  
Dropping this term in inference amounts to assuming $\mathcal{L}_u^\dagger \mathbf{f}$ is locally flat and results in an $O(\tau)$ drift bias of size $\beta^2_{1,a}\approx (\tau/4)\,\Delta_D(\mathcal{L}_u^\dagger \mathbf{f})$, underestimating anisotropy and the evolution of curvature of the flow field, so inferred flow-lines appear too straight.

\paragraph{Third component $R^{\,3}_{1,a}$ of remainder term $R_{1,a}$.}

\begin{equation}
R_{1,a}^{\,3}
= \int_{t_0}^{t} (t-u)\,\sum_{j}
\frac{\partial\,\mathcal L^\dagger_{u} \mathbf{f}}{\partial X^{(j)}}\,[\sigma\,\de \mathbf{W}_{u}]_{j},
\end{equation}

\begin{equation}
\big[R_{1,a}^{\,3}\big]_{i}
= \int_{t_0}^{t} (t-u)\,\sum_{j,m}
\frac{\partial}{\partial X^{(j)}}\big[\mathcal L^\dagger_{u} \mathbf{f}\big]_{i}
\,\sigma_{jm}\,\de \mathbf{W}_{u}^{(m)}  ,
\end{equation}

This is a martingale term capturing the stochastic coupling between diffusion and the spatial inhomogeneity of the generator. 
In inference, this term doesn't introduce bias, since $\langle R^{3}_{1,a} \rangle=0$. However, neglecting this term, ignores a second–order variance contribution with $\operatorname{Var}(R^{3}_{1,a}/\tau)=O(\tau)$.

\subsubsection{Second remainder term $R_{1,b}$}
We denote the second term of the reminder by $R_{1,b}$
\begin{equation}
R_{1,b} = \sum_{\nu=1}^n \int_{t_0}^t \int_{t_0}^s \mathcal{L}_{W, \nu}\, \mathbf{f} \, \de W_u^{(\nu)} \, \de s.
\end{equation}

Applying Fubini's theorem again to change the order of integration, we re-write $R_{1,b}$ in the form of a stochastic integral
\begin{equation}
R_{1,b} = \sum_{\nu=1}^n \int_{t_0}^t (t - u) \, \mathcal{L}_{W, \nu}\, \mathbf{f} \, \de W_u^{(\nu)}.
\end{equation}

Substituting the operator results in an expression for each dimensional component $i$
\begin{equation}
[R_{1,b}]_i = \sum_{\nu=1}^n \int_{t_0}^t (t - u) \left( \sum_{j=1}^D \frac{\partial f_i}{\partial X^{(j)}} \sigma_{j \nu} \right) \de W_u^{(\nu)},
\quad \text{for  }\, i=1,\ldots,D.
\end{equation}

In matrix notation, this corresponds to
\begin{equation}
R_{1,b} = \int_{t_0}^t (t - u) \, \mathbf{J}_f \, \boldsymbol{\sigma} \, \de \mathbf{W}_u.
\end{equation}

The remainder $R_{1,b}$ is a stochastic integral with zero mean, but non-zero  covariance, given by
\begin{equation}
\mathrm{Cov}(R_{1,b}) = \langle R_{1,b}\, R_{1,b}^\top \rangle = \int_{t_0}^t (t - u)^2 \, \mathbf{J}_f \, \boldsymbol{\sigma} \boldsymbol{\sigma}^\top \, \mathbf{J}_f^\top \, \de  u.
\end{equation}
For sufficiently smooth $\mathbf{J}_f$ and small time step $\tau = t - t_0$, this covariance scales on the order of $\tau^3$.

 The term $R_{1,b}$ quantifies the contribution to the remainder arising from stochastic fluctuations of the noise acting through the spatial derivatives of the drift $\mathbf{f}$. It does not contribute to additional systematic bias, but introduces variance in the drift estimator, especially when $\boldsymbol{\sigma}$ or $\mathbf{J}_f$ are large.

\subsubsection{Third remainder term $R_{1,c}$ }

We denote the third remainder term by $R_{1,c}$ and re-write here for convenience
\begin{equation}
R_{1,c} = \int_{t_0}^t \int_{t_0}^s \mathcal{L}^\dagger_u \boldsymbol{\sigma}(\mathbf{X}_u, u) \, \de u \, \de \mathbf{W}_s.
\end{equation}

 In the general case of time- and state- dependent diffusion the integrand of this term would be expressed for the $i$-th row and $\ell$-th column component of $\boldsymbol{\sigma} $ as follows
\begin{align}
\left[\mathcal{L}^\dagger_u \boldsymbol{\sigma}(\mathbf{X}_u, u)\right]_{i\ell} 
&= \frac{\partial}{\partial u} \sigma_{i\ell}(\mathbf{X}_u, u) + \sum_{j=1}^D \frac{\partial \sigma_{i\ell}}{\partial X^{(j)}}(\mathbf{X}_u, u) f_j(\mathbf{X}_u, u) \\
&\,\,\,\,\,\,\,\,+ \frac{1}{2} \sum_{j,k=1}^D \frac{\partial^2 \sigma_{i\ell}}{\partial X^{(j)} \partial X^{(k)}}(\mathbf{X}_u, u) [\boldsymbol{\sigma} \boldsymbol{\sigma}^\top]_{jk}(\mathbf{X}_u, u).
\end{align}

However, in our setting we consider state- and time-independent diffusion matrix, and thus $\mathcal{L}^\dagger_u \boldsymbol{\sigma}(\mathbf{X}_u, u) = \mathbf{0}$, and by consequence $R_{1,c} = \mathbf{0}$

\subsubsection{Fourth remainder term $R_{1,d}$ }

The fourth remainder term is 
\begin{equation}
R_{1,d} = \sum_{\nu=1}^n \int_{t_0}^t \int_{t_0}^s \mathcal{L}_{W, \nu}\, \boldsymbol{\sigma} \, \de W_u^{(\nu)} \, \de \mathbf{W}_s.
\end{equation}

For each component $(i,\ell)$ of $\boldsymbol{\sigma}$ 
\begin{equation}
\left[\mathcal{L}_{W, \nu} \boldsymbol{\sigma}\right]_{i \ell} = \sum_{j=1}^D \frac{\partial \sigma_{i \ell}}{\partial X^{(j)}} \, \sigma_{j \nu} \,=\, \mathbf{0}.
\end{equation}

Thus, the omission of this remainder term does not contribute any bias or variance to the EuM-based drift estimator.

\if False
\subsection{Concrete numerical evaluation of the remainder terms for different systems }
To understand better when these remainder terms contribute significantly to the discretisation, we study here closer how each term contributes to the remainder term $R_1,a$ for \textbf{(i)} a linear system under different parameter regimes that influence the curvature of the flow field, \textbf{(ii)} a nonlinear system.
Our analysis focuses mainly on the first remainder term, and in particular on $ R^{1}_{1,a}$ that depends on how quickly the vector field and diffusion operator change when moving along the flow, thus directly controlling the size of the deterministic bias.
\fi

\newpage

 \section{Additional numerical results}\label{sec:additional}

 \subsection{Inference for a 3-dimensional system}

\begin{figure}[h!]
 
  \vspace{-6pt} 
  \begin{center}
  \begin{overpic}[width=0.95\linewidth]{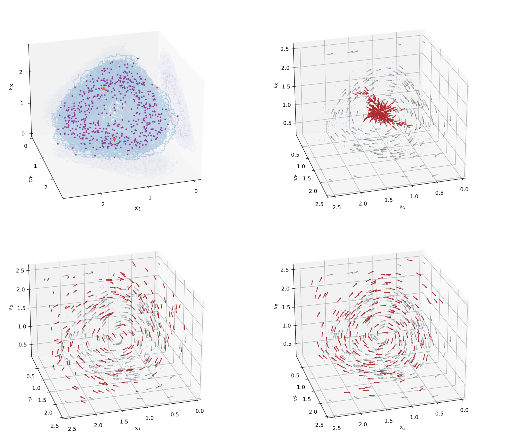}
  \put(0,75){a.}
  \put(0,30){c.}  
  \put(50,75){b.}
  \put(50,30){d.}
  \end{overpic}
  \end{center}
  \captionsetup{width=0.92\textwidth}
  \caption{\textbf{Inference for a three dimensional system following Repressilator dynamics.} \textbf{a.)}Empirical manifold of the Repressilator system. Purple dots indicate the observations and blue line denotes the continuous unobserved trajectory. Orange triangles indicate two consecutive observations. \textbf{b.)}
Ground truth (\emph{grey}) and estimated (\emph{maroon}) flow field with Gaussian process regression without augmentation \textbf{c.-d.)} Ground truth (\emph{grey}) and estimated (\emph{maroon}) flow field with the proposed framework after c.) one and d.) two augmentations.}
  \label{fig:3D} 
\end{figure}
\newpage
\subsection{Inference for 8-dimensional system}

\begin{figure}[h!]
  \begin{center}
  \begin{overpic}[width=0.95\linewidth]{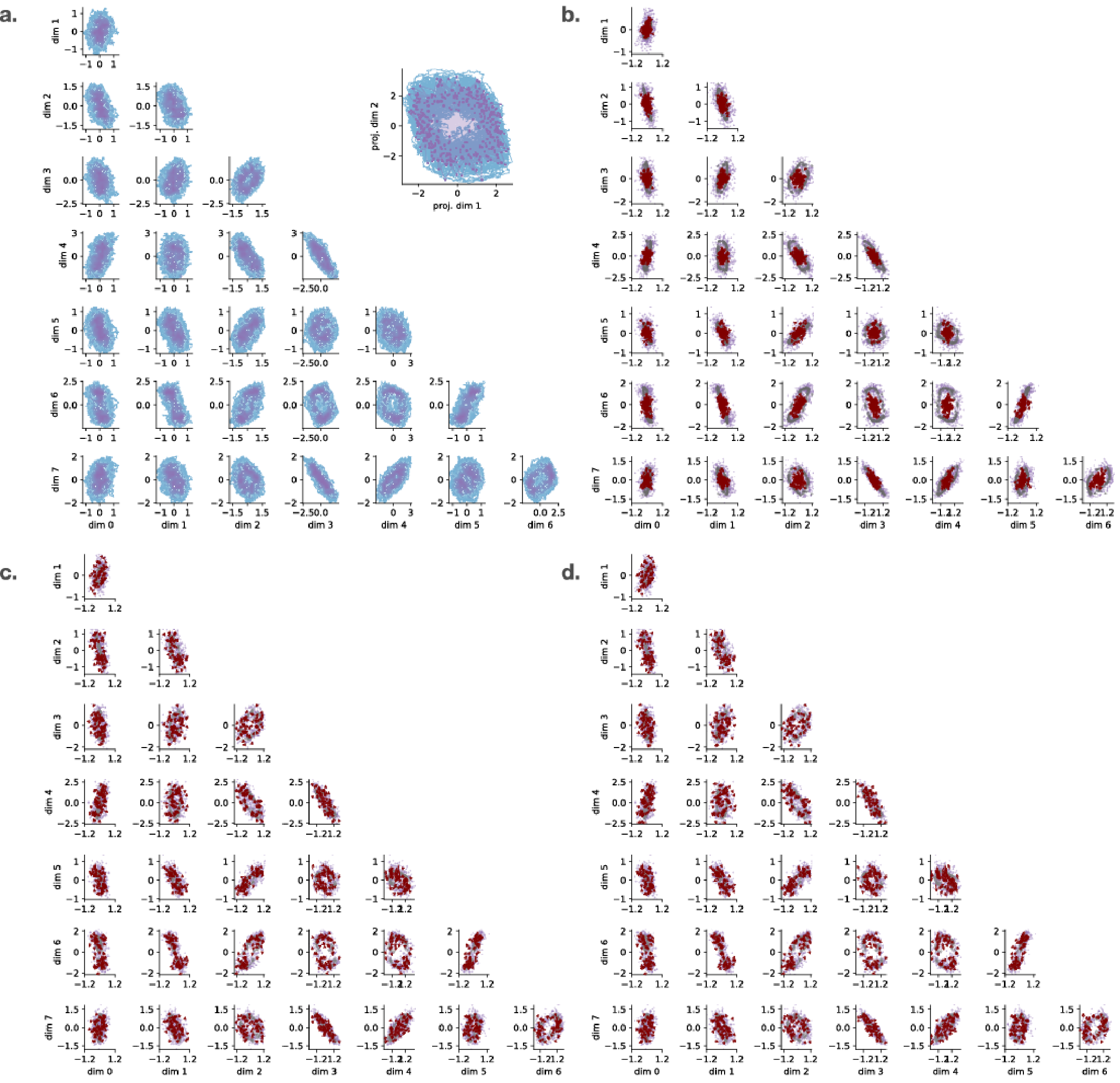}
\end{overpic}
  \end{center}
  \captionsetup{width=0.92\textwidth}
\caption{\textbf{Example inference of an 8-dimensional system. }
\textbf{a.)} Observations (\emph{purple}) and unobserved path (\emph{light blue}) across 
all dimensions, \textbf{b.)} Inferred drift (\emph{maroon}) with Gaussian process regression
 without augmentation against ground truth drift (\emph{grey}). 
The vector fields are plotted with streamlines initiated from $60$ 
random trajectories. Identified vector field after \textbf{c.)} one and 
\textbf{d.)} two augmentations.}

 \label{fig:8D} 
\end{figure}

\newpage

 \subsection{Inference with noise miss-estimation}\label{app:mise}

 \begin{figure*}[ht!] 
  \begin{center}
  
  \begin{overpic}[width=0.98\textwidth]{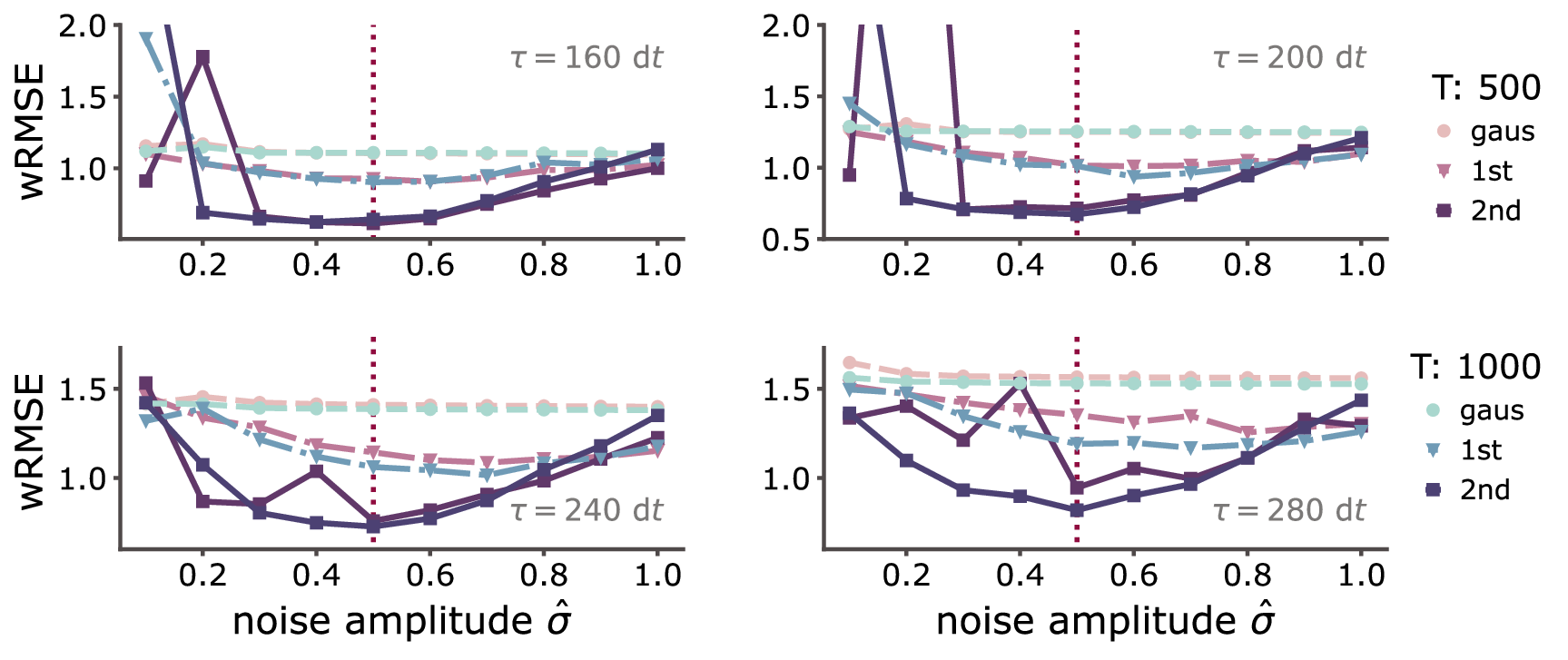}
  \end{overpic}
   \end{center}
  \captionsetup{width=0.95\textwidth}
  \caption{\textbf{Small noise misestimation has small impact on estimation accuracy.} 
Weighted root mean square error (wRMSE) vs. noise amplitude $\sigma$ employed in the augmentation for different inter-observation intervals  with \textbf{a.)} $\tau=160\,dt$ \textbf{b.)} $\tau=200\,dt$,  \textbf{c.)} $\tau=240\,dt$ \textbf{d.)} $\tau=280\,dt$. Pink-purple lines correspond to estimation with total simulation length $T=500$ time units, and blue markers correspond to total simulation length of $T=1000$ time units.
Red dotted line identifies the noise amplitude employed in the simulation of the observations.
}
  \label{fig:mis} 
\end{figure*}

\newpage
\subsection{Convergence of inference over iterations}

In the main text we have reported performance results of the proposed approach for inferring the underlying stochastic dynamics after two augmentations. This was motivated by the observation that already after two augmentations the method outperformed existing approaches, and was a trade-off between compute requirements and performance. Here we provide demonstration that our approach converges both in terms of wRMSE and negative log likelihood over multiple augmentations for a Van der Pol system.

\begin{figure*}[ht!] 
  \begin{center}
  \begin{overpic}[width=0.58\textwidth]{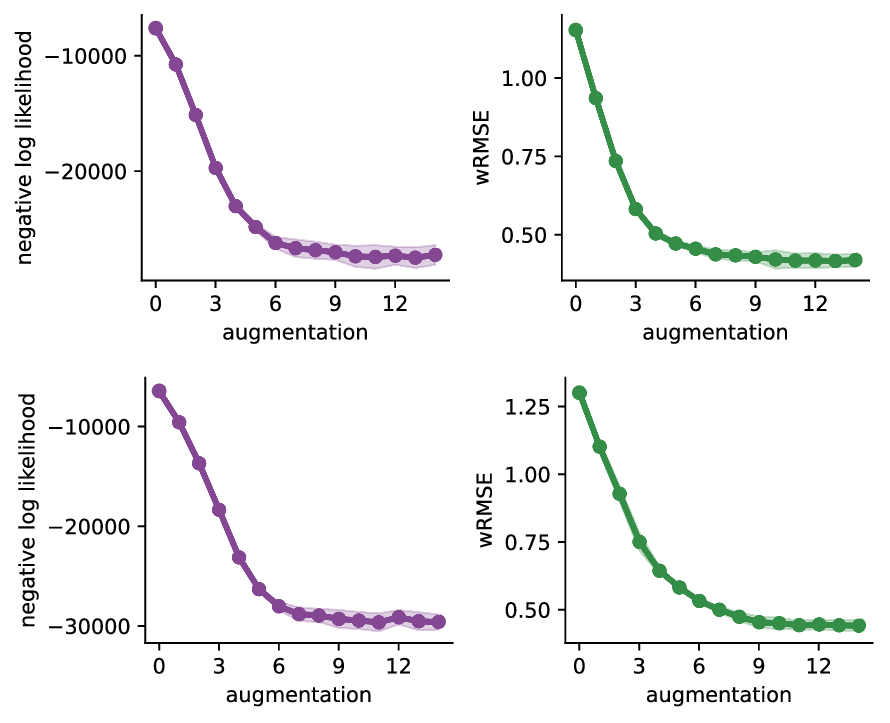}
  \end{overpic}
   \end{center}
  \captionsetup{width=0.95\textwidth}
  \caption{\textbf{Convergence of the E.M. optimisation over augmentations both in terms 
of the negative log likelihood and weighted root mean squared error for 
inter-observation interval $\tau=240\,\mathrm{d}t$ (\emph{upper}) 
and $\tau=280\,\mathrm{d}t$ (\emph{lower}) for noise amplitude $\sigma=0.25$ 
and total trajectory length $T=500$ time units.} Points indicate average 
over $5$ (\emph{upper}) and $3$ (\emph{lower}) independent runs.
}
  \label{fig:convergence} 
\end{figure*}

\newpage

\subsection{Ablations with respect to metric learning algorithm }\label{appsec:ablation_metric}

For approximating the metric induced by the observations, we employed the framework of~\cite{arvanitidis2019fast}, while we could have employed alternative metric learning approaches~\citep{scarvelis2022riemannian, hauberg2012geometric, barua2025geodesic, gruffaz2025riemannian}. However, the framework of~\cite{arvanitidis2019fast} perfectly fitted the purposes of our work, because it employs
 a non-parametric (kernel) estimation for approximating the metric and
computes the geodesics through GP regression. This allowed us to evaluate the geodesic equation at different
increments, that is necessary for imposing the time dependent geodesic constraint. 
A similar metric approximation has been recently employed in~\cite{kapusniak2024metric} for metric flow matching, i.e., for augmentation that respects the geometry of the dataset. While our approach has a similar flavour to this work, our framework additionally requires the augmented data to be temporary ordered and to respect the stochastic flow of the estimated system. This results in learning a global drift that approximates the underlying stochastic dynamics, instead of learning a local drift that transports a snapshot of states from some initial to a final configuration.

To probe the robustness of our framework, when we employ a different approach to estimate the metric, following~\citep{kapusniak2024metric} we tested our method when we employ a radial based function approximation to estimate the diagonal metric, similar to~\cite{arvanitidis2017latent}.
In the table~\ref{tab:LCmetric} we report the performance of our method when we employ the locally adaptive normal distribution framework (Geometric (our)) for approximation the metric~\cite{arvanitidis2019fast} and when we employ the radial basis function variant of the metric approximation (Geometric$_{\text{RBF}}$ (our)) for the Van der Pol system for different inter-observation intervals and noise conditions.
We observe that the resulting drift is accurate also with the RBF method for estimating the metric, yet the method proposed in the main text performs slight better across all inter-observation intervals and noise conditions.

\begin{table*}[t]
\rad{1.3}
\setlength{\tabcolsep}{6pt}
\rowcolors{2}{white}{mygreen}
\resizebox{\textwidth}{!}{
\begin{tabular}{rrrrcrrrcrrr}\toprule
& \multicolumn{6}{c}{\textbf{Van der Pol}} & \phantom{abc}&\\
\cmidrule{2-5} \cmidrule{6-8} \cmidrule{10-11}
\textbf{wRMSE} $\downarrow$& \makecell{total duration \\T} &  $\tau=80\times \de t$ & $\tau=120\times \de t$ & $\tau=160\times \de t$  &  $\tau=200\times \de t$ & $\tau=240\times \de t$ & $\tau=280$ $\times \de t$\\ \midrule
$\boldsymbol{\sigma}\mathbf{=0.25}$ \\    
MFM$_{\text{RBF}}$&1500 &1.462 $\pm$ 0.007 & 1.469 $\pm$ 0.005 & 1.470 $\pm$ 0.012 & 1.469 $\pm$ 0.008 & 1.469 $\pm$ 0.006 & 1.466 $\pm$ 0.008&\\
{MFM}$_{\text{LAND}}$&1500 &1.463 $\pm$ 0.007 & 1.469 $\pm$ 0.005 & 1.469 $\pm$ 0.012  &1.469 $\pm$ 0.008 & 1.469 $\pm$ 0.006 & 1.467 $\pm$ 0.008& \\
\textbf{Geometric$_{\text{RBF}}$ (our)}&500 &  
0.419 $\pm$ 0.052& 0.458 $\pm$ 0.063& 0.493 $\pm$ 0.031& 0.517 $\pm$ 0.022& 0.657 $\pm$ 0.040& 1.001 $\pm$ 0.077 &\\
\textbf{Geometric$_{\text{LAND}}$ (our)}&500 & 0.474 $\pm$ 0.034& $\mathbf{0.413 \pm 0.016}$& $\mathbf{0.514 \pm 0.068}$& $\mathbf{0.578 \pm 0.022 }$&0.687 $\pm$ 0.032& 0.993 $\pm$ 0.037&\\
\cmidrule{1-8}
$\boldsymbol{\sigma}\mathbf{=0.50}$\\
{MFM}$_{\text{RBF}}$&1500 & 1.516 $\pm$ 0.011 & 1.525 $\pm$ 0.006 & 1.538 $\pm$ 0.009 & 1.537 $\pm$ 0.017 & 1.528 $\pm$ 0.015 & 1.544 $\pm$ 0.019 &\\
{MFM}$_{\text{LAND}}$&1500 & 1.517 $\pm$ 0.011 & 1.526 $\pm$ 0.006 & 1.536 $\pm$ 0.009 & 1.537 $\pm$ 0.017 & 1.528 $\pm$ 0.015 & 1.545 $\pm$ 0.019&\\        
\textbf{Geometric$_{\text{RBF}}$ (our)}&500 & 0.653 $\pm$ 0.014& 0.690 $\pm$ 0.026& 0.694 $\pm$ 0.026& 0.761 $\pm$ 0.050& 0.798 $\pm$ 0.047& 0.933 $\pm$ 0.160 &\\
\textbf{Geometric$_{\text{LAND}}$ (our)}&500 &$\mathbf{0.462 \pm 0.019}$&$\mathbf{ 0.541 \pm 0.023}$& $\mathbf{0.621 \pm 0.012} $& $\mathbf{0.675 \pm 0.030}$& $\mathbf{0.750 \pm 0.038}$& $\mathbf{0.865 \pm 0.057}$&\\
\bottomrule
\end{tabular}}
\captionsetup{width=\textwidth}
\caption{Performance comparison in terms of weighted root mean square error (wRMSE) of proposed framework and Metric flow matching framework~\citep{kapusniak2024metric} employing the LAND~\citep{arvanitidis2019fast} and RBF~\citep{arvanitidis2017latent} metric approximation algorithm for different noise conditions $\sigma$ and inter-observation intervals $\tau$ for the Van der Pol system. }
\label{tab:LCmetric}
\end{table*}

\if False
\subsection{Robustness of geodesic estimation and geometric constraints}

To demonstrate the robustness of the geometric approximation of the invariant density, we set up a series of numerical experiments to systematically quantify \textbf{(i)} the \textbf{stability of the metric approximation}, \textbf{(ii)} the \textbf{stability} of the identified \textbf{geoedesics}, and \textbf{(iii)} the \textbf{convergence of the identified geodesics} to the \textbf{ground truth mean of the transition density between consecutive observations} when inter-observation interval $\tau$, total duration length $T$ of the time series vary for different values of the hyper-parameters of the metric approximation and the geodesic construction.

To assess how good the invariant geometry induced by the observations can be approximated by the metric approximation of Eq.~\ref{eq:metric_approx}, we characterise the stability of the approximated metric over different realisations of the observations, for different observation lengths $T$ and inter-observation intervals. 

To that end, for each considered system, we generate $\ell = 10$ independent trajectories

The parameter $\sigma_w$ involved in the metric approximation controls how fast the metric changes.

 a smoothly
changing metric that induces a Riemannian manifold


\textbf{(i) Across realizations (primary test).}
For each stochastic system and sampling configuration $(\Delta t, T)$, we generate multiple independent trajectories (different random seeds) from the ground-truth SDE. Using fixed metric hyperparameters $(\sigma_w^\star,\varepsilon^\star)$, we estimate a metric $G^{(s)}(x)$ from each realization $s$ and compare them on a common evaluation set $\{x_m\}$ using an SPD-aware discrepancy, e.g.\ the (log-)Euclidean or affine-invariant Riemannian distance:
\begin{equation}
D(G^{(s_1)},G^{(s_2)}) \;=\; \frac{1}{M} \sum_{m=1}^M 
d\!\left(G^{(s_1)}(x_m),\,G^{(s_2)}(x_m)\right).
\end{equation}
Small across-seed discrepancies indicate that the learned geometry is stable under finite-sample variability for the given $(\Delta t, T)$; large discrepancies signal that the observations are insufficient to reliably constrain the metric. This is the primary robustness criterion, as it directly answers whether repeating the same experiment would recover the same geometry.

\textbf{(ii) Across hyperparameters (secondary test).}
Conditioned on a representative dataset in a regime where across-seed variability is low, we probe sensitivity to the metric hyperparameters. For a grid of $(\sigma_w,\varepsilon)$ around $(\sigma_w^\star,\varepsilon^\star)$, we estimate $G_{\sigma_w,\varepsilon}(x)$ on the same data and compute $D(G_{\sigma_w,\varepsilon},G_{\sigma_w^\star,\varepsilon^\star})$. We observe a broad plateau where this discrepancy remains small, with degradation only for extreme under- or over-smoothing. This confirms that, once the sampling regime provides sufficient coverage, our results are not knife-edge with respect to $(\sigma_w,\varepsilon)$.

\textbf{Link to geodesics and dynamics.}
In the main text, we report how regimes with low across-seed metric discrepancy also exhibit (i) low variability of geometry-induced geodesics across realizations and (ii) small deviation between these geodesics and oracle paths derived from the known dynamics (Fokker--Planck effective drift or Onsager--Machlup paths). Thus, reliable geometric constraints emerge precisely when the metric is both stable across seeds and robust to moderate hyperparameter changes.

For each sampling sparsity $\tau$, we can quantify how long one must record to obtain reliable geometry-induced geodesics, and show that in this regime geodesics are both stable across resamples and consistent with the ground-truth transition density mean. The behaviour degrades in a controlled way for very short recordings or extreme hyperparameters


\fi

\subsection{Comparison of geodesic curves with Onsager--Machlup minimizers}\label{supp:om_comparison}

The \textbf{Onsager--Machlup} (OM) functional, we introduced in Section~\ref{app:OM}, provides a variational characterisation of the \textbf{most probable path} (MPP) between two system states over a prescribed time interval. In our framework, the geodesic curves, we employed to guide bridge augmentation, are meant to serve as coarse empirical approximations of the OM-minimisers, i.e. of the most probable paths between successive observations. A natural question, therefore, is how accurately these inferred geodesics capture the corresponding OM-minimisers. To address this question, we compare the identified geodesic curves to the ground-truth OM-minimisers, and compare their \textbf{distance} and \textbf{action gap} to the true OM-minimisers with that of a simple straight interpolant between consecutive observations. We perform these comparisons for different inter-observation intervals $\tau \in \{200,\,240,\,280\}\times \de t$, total trajectory lengths $T \in \{200,\,300,\,400,\,500\}$ time units, values of the hyper-parameter $\sigma_{\mathcal{M}}\in \{0.1,\,0.2,\,0.3,\,0.4,\,0.5,\, 1.0\}$ controlling the width of the kernel employed in the metric approximation, and for different levels of dynamical noise $\sigma \in \{0.1,\,0.25,\,0.50 \}$.

More precisely, among all sufficiently smooth curves $\mathbf{x}:[0,\tau]\to\mathbb{R}^d$ with boundary conditions
\begin{equation}
    \mathbf{x}(0)=\mathbf{x}_a,
    \qquad
    \mathbf{x}(\tau)=\mathbf{x}_b,
\end{equation}
the OM-minimiser is the path that minimises the \textbf{OM-action}
\begin{equation}
    \mathcal{S}_{\mathrm{OM}}[\mathbf{x}] = \int_0^\tau \left[ \frac{1}{4\,D}\,\|\dot{\mathbf{x}}(t)-\mathbf{f}(\mathbf{x}(t))\|_2^2
        + \frac{1}{2}\,\nabla\!\cdot \mathbf{f}(\mathbf{x}(t))
    \right]\de t.
    \label{app_eq:om_functional_continuous}
\end{equation}
Thereby, the OM-minimiser represents the most probable continuous-time path connecting two consecutive observations under the ground-truth stochastic dynamics for fixed inter-observation interval $\tau$.

In the framework presented in the main text, we employ geodesic curves on the empirical manifold as guidance for path augmentation. Since these curves are not, in general, identical to OM-minimisers, it is natural to quantify how well they approximate the corresponding most probable paths, and contrast its performance to a simple interpolant (straight curve) that joins the considered observations. To that end, for each pair of observations, we compare the identified geodesic curve $\boldsymbol{\gamma}_{\mathrm{geo}}$ with the OM-minimiser $\mathbf{x}_{\mathrm{OM}}^\star$ by measuring both a path-space distance and an action-based discrepancy, and compare these two discrepancies with the discrepancies of a straight interpolant to the OM-minimiser.

We characterise the path-space discrepancy between the two curves using their time-averaged $L^2$ path distance
\begin{equation}
    d_{\mathrm{path}}  =\left( \frac{1}{\tau}
        \int_0^\tau  \|\boldsymbol{\gamma}_{\mathrm{geo}}(t)-\mathbf{x}_{\mathrm{OM}}^\star(t)\|_2^2\,\de t \right)^{1/2}.
    \label{eq:path_distance_continuous}
\end{equation}
This quantifies how close the two curves remain in state space over the inter-observation interval.

Second, we compare the curves in terms of their OM action. Since the OM minimiser is defined as the minimiser of Eq.\ref{app_eq:om_functional_continuous}, the \textbf{action gap}, i.e. the difference between the action evaluated for the most probable path and for the geodesic curve
\begin{equation}
    \Delta \mathcal{S}_{\mathrm{geo}} =
    \mathcal{S}_{\mathrm{OM}}[\boldsymbol{\gamma}_{\mathrm{geo}}] -
    \mathcal{S}_{\mathrm{OM}}[\mathbf{x}_{\mathrm{OM}}^\star].
    \label{eq:action_gap_continuous}
\end{equation}
We expect $\Delta \mathcal{S}\geq 0$, and thus values close to zero indicate that the geodesic curve is nearly optimal in the OM sense.

Since the absolute value of the action can vary across settings, to compare how the discrepancy between the geodesic and the OM-minimiser changes for different values of $\sigma_{\mathcal{M}}$, total trajectory length $T$, and inter-observation interval $\tau$, we consider the
 \textbf{relative action gap}
\begin{equation}
    \Delta \mathcal{S}^{\mathrm{rel}}_{\mathrm{geo}}
    =
    \frac{
        \mathcal{S}_{\mathrm{OM}}[\boldsymbol{\gamma}_{\mathrm{geo}}]
        -
        \mathcal{S}_{\mathrm{OM}}[\mathbf{x}_{\mathrm{OM}}^\star]
    }{
        \bigl|\mathcal{S}_{\mathrm{OM}}[\mathbf{x}_{\mathrm{OM}}^\star]\bigr|+\varepsilon
    },
    \label{eq:relative_action_gap_continuous}
\end{equation}
where $\varepsilon= 10^{-12}$. This normalized quantity facilitates comparison across different parameter regimes.

We quantified the aforementioned discrepancies across different values of the metric-approximation parameter $\sigma_{\mathcal{M}}$, different total trajectory lengths, and different inter-observation intervals $\tau$. For each parameter setting we used $3$ independent runs, and within each run we evaluate the above quantities on $20$ different inter-observation intervals. Thus, for each configuration we obtain empirical summaries of the path distance Eq.~\ref{eq:path_distance_continuous}, the action gap Eq.~\ref{eq:action_gap_continuous}, and the relative action gap Eq.~\ref{eq:relative_action_gap_continuous} over a collection of independently generated bridge problems.

\paragraph{Discrete OM action.}
The OM minimiser is computed on a discretized time grid, that we denote here with
\begin{equation}
    0=t_0<t_1<\cdots<t_M=\tau.
\end{equation}
 $X_k \approx x(t_k)$ identifies a discrete path with fixed endpoints
\begin{equation}
    \mathbf{X}_0 = \mathbf{x}_a,
    \qquad
    \mathbf{X}_M = \mathbf{x}_b.
\end{equation}
Using midpoint discretisation, we approximate the continuous OM functional by
\begin{equation}
    \mathcal{S}_{\mathrm{OM}}^{(M)}(\mathbf{X}_{0:M})
    =
    \sum_{k=0}^{M-1}
    \left[
        \frac{1}{4D}
        \left\|
            \frac{\mathbf{X}_{k+1}-\mathbf{X}_k}{t_{k+1}-t_k}
            -
            f\!\left(\frac{\mathbf{X}_{k+1}+\mathbf{X}_k}{2}\right)
        \right\|^2
        +
        \frac{1}{2}
        \nabla\!\cdot f\!\left(\frac{\mathbf{X}_{k+1}+\mathbf{X}_k}{2}\right)
    \right]
    (t_{k+1}-t_k).
    \label{eq:om_functional_discrete}
\end{equation}
We define then discrete OM-minimiser as
\begin{equation}
    X_{0:M}^{\star}
    =
    \arg\min_{\mathbf{X}_1,\dots,\mathbf{X}_{M-1}}
    \mathcal{S}_{\mathrm{OM}}^{(M)}(\mathbf{X}_{0:M}),
    \label{eq:discrete_om_minimizer}
\end{equation}
and optimise only the interior points, while keeping the endpoints fixed.

\paragraph{Numerical computation of the OM-minimiser.}
The continuous OM problem is a two-point boundary-value variational problem: we want to identify the path connecting the prescribed endpoints that minimises the OM-action \eqref{app_eq:om_functional_continuous}. Rather than solving the Euler--Lagrange boundary-value equations directly, we solve the equivalent discrete optimisation problem \eqref{eq:discrete_om_minimizer}. 

To improve robustness, especially for longer intervals or smaller noise levels, where the action landscape can become stiff, we resorted to a multi-scale optimisation. More precisely, we start the optimisation on a coarse temporal grid obtained by subsampling the finest target grid, solve the discrete OM minimisation problem on that coarse grid, and then use the resulting path as the initialisation for the next finer discretisation. This refinement procedure is repeated until the target finest grid is reached. Thus, the optimisation proceeds through a sequence of nested discretisations (here we employed 3 nested discretisation). 

\paragraph{Computation of action gap and relative action gap.}
We denote the geodesic curve on the discretised grid by $X^{\mathrm{geo}}_{0:M}$. We evaluate both the geodesic path and the identified OM-minimising path using the \emph{same} discrete OM functional \eqref{eq:om_functional_discrete}. The discrete action gap is therefore
\begin{equation}
    \Delta \mathcal{S}^{(M)}
    =
    \mathcal{S}_{\mathrm{OM}}^{(M)}(\mathbf{X}^{\mathrm{geo}}_{0:M})
    -
    \mathcal{S}_{\mathrm{OM}}^{(M)}(\mathbf{X}^\star_{0:M}),
    \label{eq:discrete_action_gap}
\end{equation}
and the corresponding discrete relative action gap is
\begin{equation}
    \Delta \mathcal{S}^{(M)}_{\mathrm{rel}}
    =
    \frac{
        \mathcal{S}_{\mathrm{OM}}^{(M)}(\mathbf{X}^{\mathrm{geo}}_{0:M})
        -
        \mathcal{S}_{\mathrm{OM}}^{(M)}(\mathbf{X}^\star_{0:M})
    }{
        \bigl|\mathcal{S}_{\mathrm{OM}}^{(M)}(\mathbf{X}^\star_{0:M})\bigr|+\varepsilon
    }.
    \label{eq:discrete_relative_action_gap}
\end{equation}
Likewise, we compute the discrete path distance on the same grid as
\begin{equation}
    d_{\mathrm{path}}^{(M)}
    =
    \left(
        \frac{1}{M+1}
        \sum_{k=0}^{M}
        \|\mathbf{X}_k^{\mathrm{geo}}-X_k^\star\|^2
    \right)^{1/2}.
    \label{eq:discrete_path_distance}
\end{equation}

Together, these quantities allow us to assess not only whether the geodesic curves remain geometrically close to associated OM-minimisers, but also whether they are nearly optimal in the sense of minimising the OM-action. The former addresses similarity in state space, while the latter directly measures how well the geodesic approximates the most probable path under the assumed stochastic dynamics.

\begin{figure*}[ht!] 
  \begin{center}
  \begin{overpic}[width=0.98\textwidth]{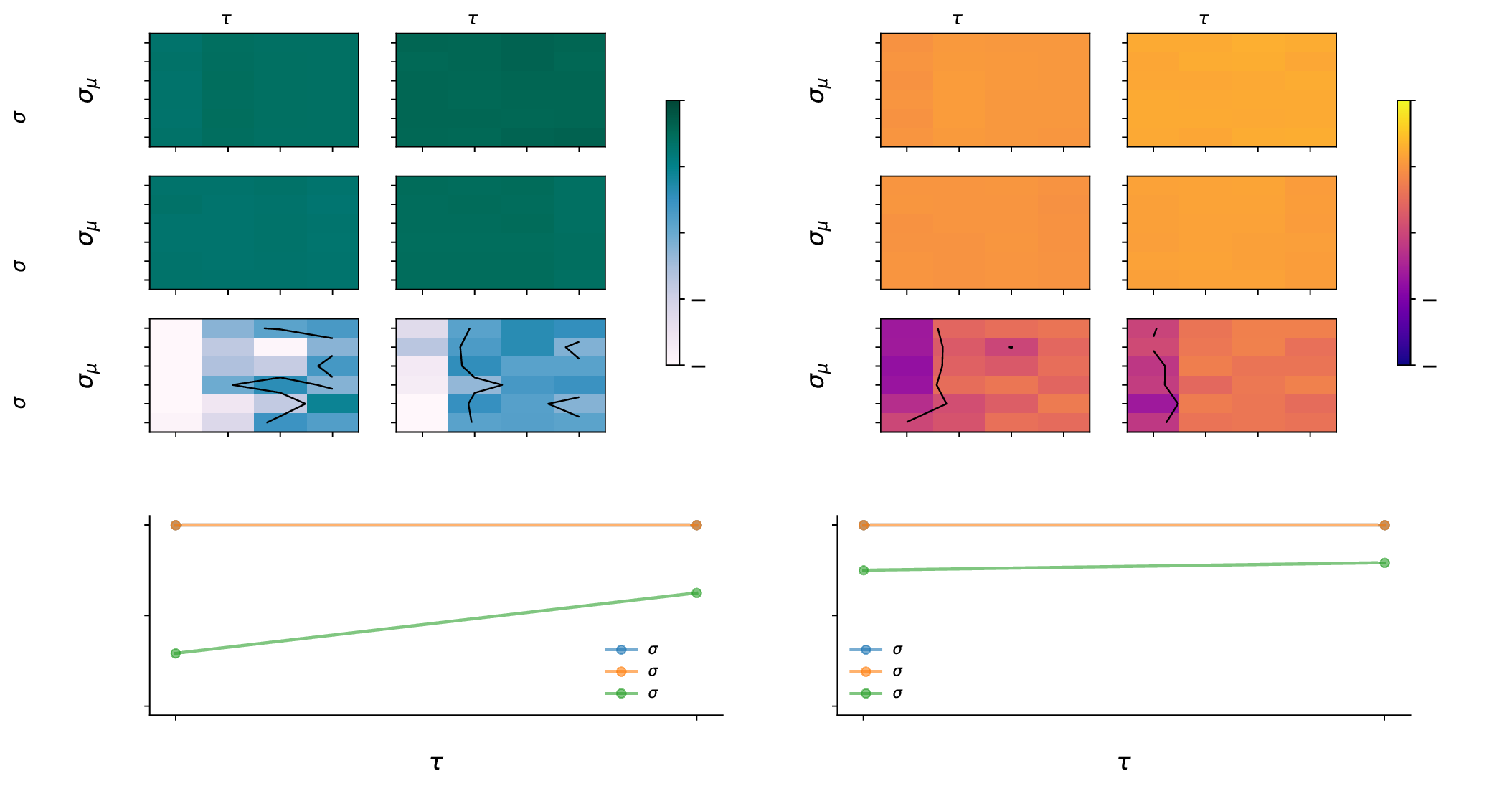}
  \end{overpic}
   \end{center}
  \captionsetup{width=0.95\textwidth}
  \caption{\textbf{Relative performance of geodesic curves in approximating the OM-minimisers against linear interpolants.}  ($\de t=0.01$) 
\textbf{a.)} improvement in OM action gap, quantified by $I_A=1-
\frac{\Delta A_{\mathrm{geo}}}{\Delta A_{\mathrm{lin}}},
\qquad \Delta A_{\mathrm{geo}}=
A[\gamma_{\mathrm{geo}}]-A[\gamma_{\mathrm{OM}}],\qquad
\Delta A_{\mathrm{lin}}=
A[\gamma_{\mathrm{lin}}]-A[\gamma_{\mathrm{OM}}].$ \textbf{b.)} Improvement in path distance to the OM-minimiser, quantified by $I_D
=
1-
\frac{d(\gamma_{\mathrm{geo}},\gamma_{\mathrm{OM}})}
{d(\gamma_{\mathrm{lin}},\gamma_{\mathrm{OM}})}$.
Positive values indicate parameter settings in which the geodesic curve provides an improved approximation to the OM-minimiser over the linear interpolant, $I=0$ indicates equal performance, and negative values indicate settings in which the linear interpolant is closer to the OM minimiser. 
Results computed 
across trajectory length $T$, metric bandwidth $\sigma_\mu$, dynamical noise $\sigma$, and inter-observation interval $\tau$, whether the geodesic reference curve is closer to the OM minimiser than a straight-line interpolant between the same endpoints. 
The black contour marks the transition $I=0$. 
\textbf{c.), d.)}  Summary plots indicating for each $\tau$ and $\sigma$, the fraction of tested $(T,\sigma_\mu)$ settings in which the geodesic outperforms the linear interpolant. 
\textbf{c.)}  Fraction of geodesics with smaller OM action gap, $\Delta A_{\mathrm{geo}} < \Delta A_{\mathrm{lin}}$.  \textbf{d.)}  Fraction of geodesics with smaller path distance to the OM minimiser, $d(\gamma_{\mathrm{geo}},\gamma_{\mathrm{OM}}) < d(\gamma_{\mathrm{lin}},\gamma_{\mathrm{OM}})$. 
Values closer to one indicate that the geodesic curves consistently provide a better approximation to the OM minimiser across metric bandwidths and trajectory lengths, while values closer to zero indicate that the straight-line interpolant provides a better approximation.
}
  \label{fig:omimprove} 
\end{figure*}

\begin{figure*}[ht!] 
  \begin{center}
  \begin{overpic}[width=0.98\textwidth]{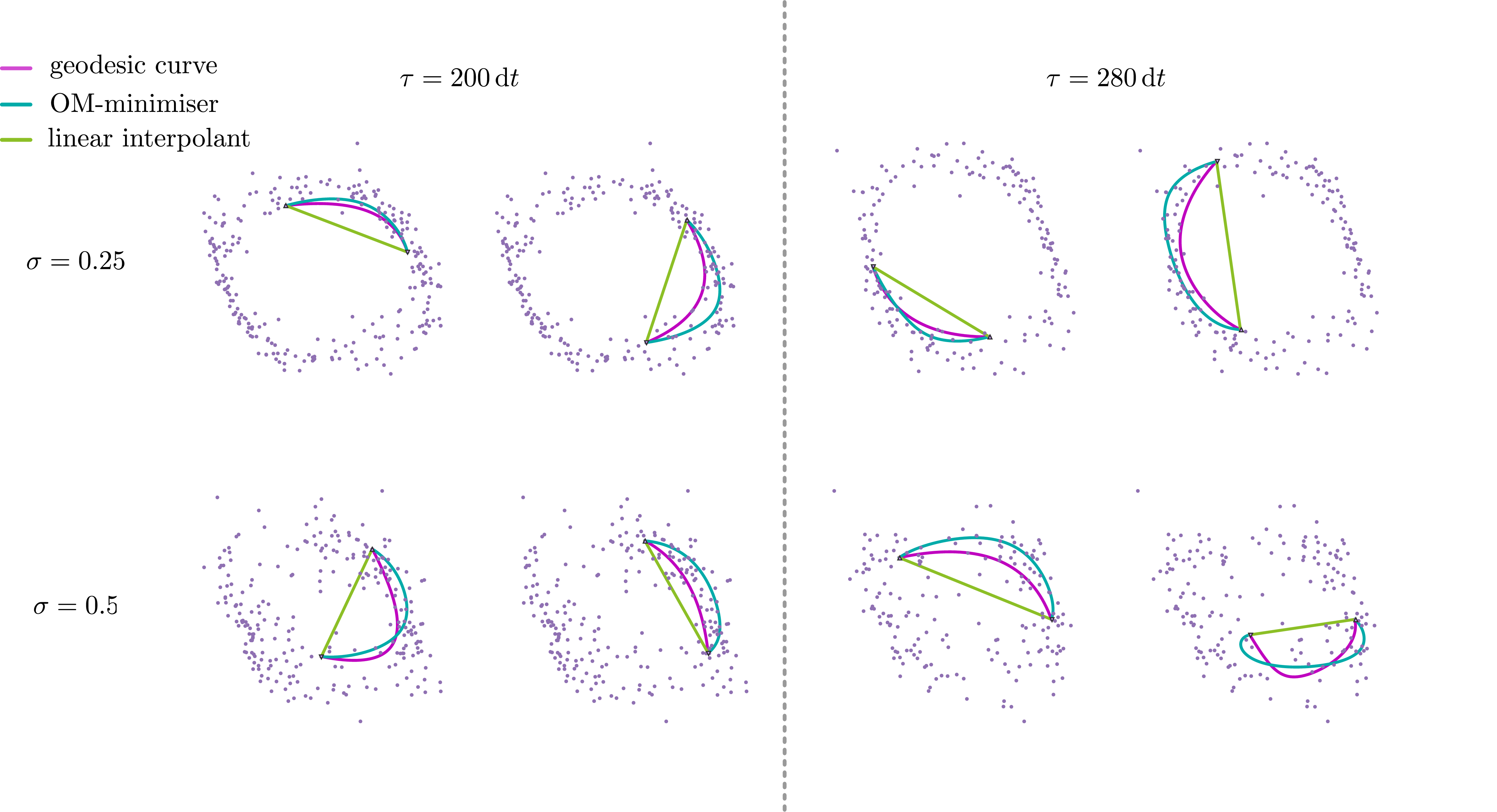}
  \end{overpic}
   \end{center}
  \captionsetup{width=0.95\textwidth}
  \caption{\textbf{Example Onsager-Machlup most probable paths, identified geodesics, and linear interpolants for the Van der Pol system.} Identified geodesics, OM-minimisers under the ground truth dynamics, and linear interpolants for two example inter-observation intervals for $\sigma=\{0.25, 0.50\}$ and inter-observation interval $\tau=200\,\de t$ and $\tau=280\,\de t$. ($\de t=0.01$)
}
  \label{fig:example-OM} 
\end{figure*}

Example OM-minimisers, geodesics curves and linear interpolants are illustrated in Fig.~\ref{fig:example-OM}.

\newpage

\newpage

 \section{Details on numerical experiments}\label{sec:details}
 We simulated a two dimensional Van der Pol oscillator
 with drift function 
 \begin{align}
        f_1(x,y) &= \mu(x - \frac{1}{3}x^3-y)\\
        f_2(x,y) &= \frac{1}{\mu}x,
 \end{align}
 starting from initial condition $ x0 = [1.81, -1.41]$
and under noise amplitudes $\sigma=\{0.25, 0.50, 0.75, 1.00\}$ for total duration of $T=\{500, 1000\}$ time units. The employed inter-observation intervals $\tau=\{80, 120, 160, 200, 240, 280,  320\} \times \de t$. The last inter-observation interval exceeds the half period of the oscillator and thus samples only a single state per period. This resulted in erroneous estimates. 
In this setting this indicates the upper limit of $\tau$ for which we can provide estimates. However for any inference method, if the observation process samples only one observation per period, identifying the underlying force field without additional assumptions is not possible with temporal methods. The discretisation time-step used for simulation of the ground truth dynamics, and path augmentation $\delta t=0.01$.
For sampling the controlled bridges we employed $N=100$ particles evolving the associated ordinary differential equation as described in~\citep{maoutsa2022deterministic}. The logarithmic gradient estimator used $M=40$ inducing points. The sparse Gaussian process for estimating the drift was based on a sparse kernel approximation of $S=300$ points. In the presented simulation we have employed a weighting parameter $\beta = 0.5$ (Eq.~\ref{apeq:free_energy2}). This provides a moderate pull towards the invariant density. The example in Figure~\ref{fig:res} was constructed with $\beta=1$ and provides a better approximation of the transition density, than $\beta=0.5$.

For the \textbf{out-of-equilibrium process} with harmonic trapping and circulation and a Gaussian repulsive obstacle in the centre we followed the description presented in~\cite{frishman2020learning} following the drift
\begin{equation}
f_{\mu}(\mathbf{x}) = -\Omega_{\mu\nu} x_\nu + \alpha e^{-x^2 / 2\sigma^2} x_\mu
\quad \text{with }  \quad
\Omega = \begin{pmatrix}
2 & 2 \\
-2 & 2
\end{pmatrix},
\end{equation}
for $\alpha = 10$ and simulated the stochastic system with noise amplitude $\sigma=0.5$ on a time grid of $\de t=0.01$ steps, observed at inter-observation intervals $\tau=\{150, 200, 250 \}\times \de t$ and for total duration $T= 1000$ time units.

For
the \textbf{Hopf system} we used the drift
\begin{align}
f_1({x_1, x_2}) &= z_2, \\
f_2({x_1, x_2}) &= -z_1 + \left(\mu - z_1^2\right) z_2,
\end{align}
with $\mu=0.35$ and integrated the system with noise amplitude $\sigma=0.15$ on a timegrid with $\de t = 0.01$ resolution, observed at $\tau=\{   200, 300, 400 \}\times \de t$ time intervals. This is the normal form of the Hopf bifurcation.

For the \textbf{Selkov glycolysis model}~\citep{sel1968self} we employed the drift
\begin{align}
f_1({x_1, x_2}) &= -x_1 + \alpha x_2 + x_1^2 x_2, \\
f_2({x_1, x_2}) &= 0.6 - \alpha x_2 - x_1^2 x_2,
\end{align}
with $a=0.06$ and noise amplitude $\sigma = 0.05$ for inter-observation intervals $\tau=\{  100, 200 \} \times \de t$ and simulation time grid of $dt=0.01$ spacing and for total duration $T= 1000$ time units.

This model is a minimal two-variable model of glycolytic oscillations, first introduced in~\citep{sel1968self}. It describes the autocatalytic feedback processes in the glycolysis pathway, focusing on how simple nonlinear interactions can give rise to oscillatory dynamics in concentrations of intermediates.
The first state variable $x_1$ represents the concentration of adenosine diphosphate, while
	$x_2$ represents the concentration of a glycolytic intermediate.

 For the \textbf{Repressilator dynamics}, we considered the system
introduced by Elowitz and Leibler.  
The repressilator consists of three genes arranged in a cyclic inhibitory loop, where each gene
represses the expression of the next one in the cycle.

The deterministic dynamics of the system are
\begin{align}
\dot{x}_1 &= \frac{\beta}{1 + \left( \frac{x_3}{k} \right)^n} - \gamma x_1, \\
\dot{x}_2 &= \frac{\beta}{1 + \left( \frac{x_1}{k} \right)^n} - \gamma x_2, \\
\dot{x}_3 &= \frac{\beta}{1 + \left( \frac{x_2}{k} \right)^n} - \gamma x_3,
\end{align}
where $\beta$ controls the maximal production rate, $\gamma$ is the degradation rate,
$k$ is the repression threshold, and $n$ is the Hill coefficient controlling the nonlinearity
of repression.

In our simulations we used $\beta=3.0$, $\gamma=1.0$, $k=1.0$, and $n=4.0$. We simulated the system with $\de t= 0.01$ for $T=1000$ time units, and observed the trajecotry every $\tau = 170 \times \de t$ time units. We used $M=120$ inducing points for the logarithmic gradient estimation, $N=400$ particles for sampling the flows, and employed noise amplitude $\sigma=0.2$

\subsection{On computation of geodesic curves}

For the computation of geodesic curves we followed the framework introduced in~\citep{arvanitidis2019fast}.  The geodesic equation relies on a non-parametric estimation of the
Riemannian metric, which is constructed using kernel-weighted local diagonal covariances, and has computational complexity $\mathcal{O}(ND)$, where  $D$ is
the dimensionality of the problem and $N$ denotes the
number of samples. The computational cost of solving the geodesic equation scales sublinearly with increasing
dimensionality.

\subsection{Details on baseline methods}\label{appsec:baselines}

We compared the performance of our method against a series of competing methods for inference of stochastic dynamics. In particular, we compared our method against methods specifically designed for inference of stochastic systems from single trajectories, and against systems that infer population dynamics.

We employed the following methods that assume single trajectories for drift inference:
\begin{enumerate}
\item Gaussian process regression without state estimation (\textbf{GP})
\item path augmentation with Ornstein-Uhlenbeck dynamics with Gaussian process inference (\textbf{OU})~\citep{batz2018approximate}
\item sparse variational inference with state estimation (\textbf{SVISE})~\citep{course2023state}
\item basis function approximation of Kramers-Moyal coefficients, i.e. the drift function (\textbf{KM-basis})~\citep{nabeel2025discovering}
\item latent SDE inference with amortized reparameterization with (\textbf{LatentSDE+GP-pre}) and without pre-training (\textbf{LatentSDE})~\citep{course2023amortized}.

\end{enumerate}

 We further compared our method with recent Schroedinger bridge generating frameworks that primary aim to infer population dynamics from snapshot data. In particular we considered the following frameworks:
 \begin{enumerate}
\item[I.] Metric Flow Matching \textbf{(MFM)} ~\citep{kapusniak2024metric}
\item[II.] Generalized Schr{\"o}dinger Bridge Matching \textbf{(GSBM)}  ~\citep{liu2023generalized}
\item[III.] Simulation-free Schr{\"o}dinger bridges via score and flow matching (\textbf{[SF]$^2$ M}) ~\citep{tong2023simulation}
 \end{enumerate}

 For these methods, we clustered the observations of each system into \emph{disjoint} subsets of adjacent points. We employed the k-Nearest neighbours algorithm~\citep{fix1985discriminatory, cover1967nearest} to construct the clusters as local neighbourhoods on the state space, 
 comprising each at most $64$ and minimum $20$ observations. 
We paired each cluster $\mathcal{J}_b$ with the set $\mathcal{J}^{\,+}_b \dot{=} \,\{\, \boldsymbol{\mathcal{O}}_{k+1} :\boldsymbol{\mathcal{O}}_k  \in \mathcal{J}_b \,\}\,$ of the next observation of each cluster member.  
We then considered each cluster pair $(\mathcal{J}_b, \mathcal{J}^+_b)$ as the initial and terminal condition for a Schr{\"o}dinger bridge problem, i.e.
\begin{align}
\pi^b_0 \,&\dot{=} \,\{\,\boldsymbol{\mathcal{O}}_k: \boldsymbol{\mathcal{O}}_k  \in \mathcal{J}_b \,\}\\
\pi^b_1 \,&\dot{=} \,\{\, \boldsymbol{\mathcal{O}}_{\ell} : \boldsymbol{\mathcal{O}}_\ell \in \mathcal{J}^+_b \,\}.
\end{align}

These serve as samples of the densities required as boundaries conditions for the Schr\"odinger bridges.

For the multi-marginal setting, starting from the cluster that contained the observation $\boldsymbol{\mathcal{O}}_1$ and subsequently created a sequence of cluster following the time ordering of the observations, i.e.
\begin{equation}
\pi^0_i \,{=} \,\{\, \boldsymbol{\mathcal{O}}_{k+i} : k \in \mathcal{J}_0 \,\}.
\end{equation}

We then employed a sequence of $50$ marginal densities $\{\pi^0_i\}^{\,49}_{\,i=0}$ as snapshot observations required by the framework.

\textbf{Metric Flow Matching.} For the Metric Flow Matching framework, we trained on observations resulting from total simulation length $T_{\text{MFM}}= 3\,T=1500$ (time units) to ensure sufficient data for each bridge. 
For each constructed bridge indexed by $b$, the flow network trained with the flow matching objective represents the velocity of the samples $\mathbf{u}_{b}(\mathbf{x},t)$ transferred within the normalised time $t\in[0,1]$ from the initial boundary condition to the terminal one. We approximate a time-independent local drift $\hat{ \mathbf{f}}_b(\mathbf{x})$ by rescaling the velocity field $\mathbf{u}_{b}(\mathbf{x},t)$ with the inter-observation interval $\tau$, i.e.,
\begin{equation}
\hat{ \mathbf{f}}_b(\mathbf{x}) \;=\; \frac{1}{\tau}\, \mathbf{u}_{b}(\mathbf{x},t) .
\end{equation}

To obtain a global drift estimate from the individual local estimates, we compute "responsibilities" or weights of each individual drift for each point $\mathbf{x}_m$ of a pre-defined two-dimensional evaluation grid that covers the state space region occupied by the observations. These weights indicate how relevant each bridge $b$ was for estimating the drift at each grid point $\mathbf{x}_m$. For each bridge, we compute support weights $\omega_b(\mathbf{x})$ on the grid employing kernel density estimation (KDE) over the bridge boundary condition samples. Then, for each grid point $\mathbf{x}_m$, we compute bridge responsibilities as
\begin{equation}
\rho_b(\mathbf{x}_m)\;=\;\frac{\omega_b(\mathbf{x}_m)}{\sum_{j=1}^B \omega_j(\mathbf{x}_m)}\,,\qquad \sum_{b=1}^B \rho_i(\mathbf{x}_m)=1.
\end{equation}
We estimate the global drift at each grid point by weighting the local estimated drifts with the corresponding bridge responsibility, i.e.,
\begin{equation}
\hat{ \mathbf{f}}(\mathbf{x}_m)\;=\;\sum_{b=1}^B \rho_b(\mathbf{x}_m)\,\hat{\mathbf{f}}_b(\mathbf{x}_m).
\end{equation}

\section{Computational complexity of the proposed method}

We provide here a breakdown of the computational cost of the main components of the proposed framework. We denote by $K$ the total number of sparse observations, $d$ the state-space dimension, $I_{\mathrm{EM}}$ the number of E.M. iterations, $L$ the number of temporal discretisation points used within each inter-observation interval, $N_p$ the number of particles used for bridge augmentation, and $M$ the number of inducing points used in the sparse Gaussian-process drift estimator. The total number of augmented state samples generated in one E.M. iteration is
\begin{equation}
N_{\mathrm{aug}} \approx (K-1)LN_p .
\end{equation}
The computational cost splits into four main components: metric estimation, geodesic construction, geometry-guided bridge augmentation, and drift estimation.

\paragraph{Metric estimation.}
We estimate the Riemannian metric once from the observed states. In the diagonal local-covariance approximation used in the experiments, evaluating the metric at $N_q$ query locations by direct summation over all observations costs
\begin{equation}
\mathcal{O}(N_q K d),
\end{equation}
with memory cost $\mathcal{O}(N_q d)$ for storing the diagonal metric values. This term is repeated only when we might want to re-estimate the metric. 

\paragraph{Geodesic construction.}
For each pair of consecutive observations $(O_k,O_{k+1})$, we compute a geodesic with respect to the learned metric. If each geodesic is represented by $L_\gamma$ discretisation points and the boundary-value solver requires $N_{\gamma}$ optimisation or solver iterations, the cost is approximately
\begin{equation}
\mathcal{O}\!\left((K-1)N_{\gamma}L_\gamma C_H(d)\right),
\end{equation}
where $C_H(d)$ denotes the cost of evaluating the metric and its derivatives along the curve. For the diagonal metric used here, this cost is substantially cheaper than for a full dense metric, for which metric inversion and Christoffel-symbol computations may introduce additional polynomial dependence on $d$. Geodesic construction is also performed independently for each inter-observation interval and is therefore parallelisable across the $K-1$ intervals. Since the geodesics are computed before the E.M. iterations and then reused, this component is smaller than the repeated bridge-sampling cost.

\paragraph{Geometry-guided bridge augmentation.}
The main computational bottleneck is the geometry-guided bridge augmentation step. At each E.M. iteration, and for each inter-observation interval, we solve a controlled bridge problem using a particle representation. The cost scales with the number of intervals, time-discretisation points, and particles as
\begin{equation}
\mathcal{O}\!\left(I_{\mathrm{EM}}(K-1)LN_p C_{\mathrm{step}}\right),
\end{equation}
where $C_{\mathrm{step}}$ is the cost of one particle propagation, weight update, and control evaluation. In our implementation, $C_{\mathrm{step}}$ includes evaluating the current drift estimate, applying the observation and geodesic potentials, and performing particle-filtering operations such as weighting and resampling. Thus, the bridge step scales as
\begin{equation}
\mathcal{O}\!\left(I_{\mathrm{EM}}(K-1)LN_p C_f(d)\right),
\end{equation}
where $C_f(d)$ is the cost of evaluating the current drift model at one state. This is the dominant computational component of the current implementation, especially as $d$, $L$, or $N_p$ increase. The memory requirement is
\begin{equation}
\mathcal{O}\!\left((K-1)LN_p d\right)
\end{equation}
if all augmented trajectories are stored. This could be reduced by streaming Monte Carlo estimates of the sufficient statistics needed for drift learning, however in the present implementation we store all the trajectories.

\paragraph{Sparse Gaussian-process drift estimation.}
Given the augmented trajectories, we estimate the drift using a sparse Gaussian-process approximation with $M$ inducing points. For $N_{\mathrm{aug}}$ augmented samples, the cost of sparse GP inference scales as
\begin{equation}
\mathcal{O}\!\left(d\,N_{\mathrm{aug}}M^2 + d\,M^3\right),
\end{equation}
since each dimensional component from the total $d$ components are estimated independently. The corresponding memory cost is approximately
\begin{equation}
\mathcal{O}\!\left(N_{\mathrm{aug}}M + M^2\right).
\end{equation}
 Since $M \ll N_{\mathrm{aug}}$, this sparse approximation avoids the cubic cost $\mathcal{O}(N_{\mathrm{aug}}^3)$ of full Gaussian-process regression.

\paragraph{Overall cost.}
Combining the above terms, the total runtime is approximately
\begin{equation}
\mathcal{O}(N_qKd)
+
\mathcal{O}\!\left((K-1)N_{\gamma}L_\gamma C_H(d)\right)
+
\mathcal{O}\!\left(I_{\mathrm{EM}}(K-1)LN_p C_{\mathrm{step}}\right)
+
\mathcal{O}\!\left(I_{\mathrm{EM}}d\,N_{\mathrm{aug}}M^2 + I_{\mathrm{EM}}d\,M^3\right).
\end{equation}
In practice, the bridge augmentation term dominates because it is repeated at every E.M. iteration and scales with the number of particles, the temporal resolution of the bridge, and the number of inter-observation intervals. By contrast, metric estimation and geodesic construction are performed once and are parallelisable across query locations or observation pairs. The $K-1$ controlled bridge problems are also independent conditional on the current drift estimate, so the augmentation step is naturally parallelisable across inter-observation intervals.

\begin{table}[ht]
\centering
\begin{tabular}{lll}
\toprule
Component & Leading cost & Main scaling variables \\
\midrule
Metric estimation
& $\mathcal{O}(N_q K d)$
& observations, query points, dimension \\

Geodesic construction
& $\mathcal{O}((K-1)N_{\gamma}L_\gamma C_H(d))$
& intervals, geodesic resolution, metric cost \\

Bridge augmentation
& $\mathcal{O}(I_{\mathrm{EM}}(K-1)LN_p C_{\mathrm{step}})$
& E.M. iterations, particles, bridge resolution \\

Sparse GP drift update
& $\mathcal{O}(I_{\mathrm{EM}}dN_{\mathrm{aug}}M^2 + I_{\mathrm{EM}}dM^3)$
& augmented samples, inducing points, dimension \\
\bottomrule
\end{tabular}
\caption{Runtime scaling of the main components of the geometry-guided augmentation framework. The dominant term in the current implementation is the particle-based bridge augmentation step.}
\label{tab:runtime_breakdown}
\end{table}
\vspace{-8pt}
In the present work we prioritised demonstrating the inference method in the sparse-observation setting and establishing the usefulness of geometry-guided augmentation for drift inference. The current implementation is therefore not intended to be a fully optimised scalable solver for high-dimensional systems. More scalable implementations could use amortised bridge samplers, adaptive particle allocation, reduced-rank or local approximations of the control, mini-batched drift updates, or GPU-exploiting implementation. Developing such improvements is an important direction for future work.

\section{Algorithmic details}

Here we provide the outline algorithm for each constituent component of our work.
Algorithm~\ref{algo:skeleton} provides the main skeleton of the framework. For the geometric approximation and the construction of the geodesics we defer the readers to~\cite{arvanitidis2019fast}. Algorithm~\ref{algo:control} outlines the solution of the control problem that implements the path augmentation. This part is an adapted version of the main algorithm proposed by~\cite{maoutsa2021deterministica}. Finally, Algorithm~\ref{algo:drift} describes the solution of the Gaussian process inference given the path augmentations (bridges) created for each augmentation pair. For the simulation of Fokker-Planck equation solutions we used the deterministic particle framework of~\cite{Maoutsa_2020}.

\renewcommand\thealgocf{A\arabic{algocf}}
\vspace{12pt}
\begin{algorithm}[H]\label{alg:outer}
\SetAlgoLined
\DontPrintSemicolon
\KwIn{
    $\boldsymbol{\mathcal{O}} = \{ (\mathbf{x}_k, t_k) \}_{k=1}^{K}$: observed states at timepoint $t_k$\\
}
\KwOut{
    $\hat{\mathbf{f}}$: posterior estimate of the drift function\\
    \hspace*{35pt}$\boldsymbol{B}^{(j)}$: augmented paths of latent states (optional)
}
\BlankLine

\tcp*{initialise $\hat{\mathbf{f}}$ with a coarse drift estimate }
Initialise drift estimate $\hat{\mathbf{f}}^{(0)}$ according to Eq.~\ref{eq:increments}\;

    \tcp*{Approximate Riemannian metric from observations (Eq.~\ref{appeq:metric})}
    ${H}_{dd} =$ ApproximateMetric$\left(\{\boldsymbol{\mathcal{O}_k}\}^K_{k=1}\right)$\;
    \tcp*{}

    \BlankLine
    \tcp*{Construct geodesics between $\boldsymbol{\mathcal{O}_k}$ and $\boldsymbol{\mathcal{O}_{k+1}}$ under the estimated metric as shortest paths}
    $\boldsymbol{\Gamma}^{(\ell)} =$ ConstructGeodesics$(\boldsymbol{\mathcal{O}}_k, \boldsymbol{\mathcal{O}}_{k+1} , H_{dd})$\;
    \tcp*{$\boldsymbol{\Gamma}^{k} = \{ \boldsymbol{\gamma}^{k}_{t^\prime} \}_{k=1}^{K}$ geodesic curves between selected observation pairs}

    \BlankLine
    \For{each iteration $j$}{
    \tcp*{ augment paths along geodesics using particle flow}
    $\boldsymbol{B}^{(j)} =$ AugmentPaths$(\{\boldsymbol{\mathcal{O}}\}^K_{k=1}, \boldsymbol{\Gamma}^{(j)}, \hat{\mathbf{f}}^{(j-1)})$\;
    \tcp*{uses the deterministic particle flow / bridge construction (Alg.~\ref{algo:control}) to sample augmented trajectories with $\hat{\mathbf{f}}^{(j-1)}$}

    \BlankLine
    \tcp*{ Gaussian process inference of the drift function}
    $\hat{\mathbf{f}}^{(\ell)} =$ GPDriftInference$(\{\boldsymbol{\mathcal{O}}\}^K_{k=1}, \boldsymbol{B}^{(j)})$\;
    \tcp*{update GP posterior over $\mathbf{f}$ using original and augmented data}
}

\BlankLine

\caption{{Skeleton of the proposed framework.}}
\label{algo:skeleton}
\end{algorithm}

\renewcommand\thealgocf{A\arabic{algocf}}
\vspace{12pt}
\begin{algorithm}[H]\label{alg1}
\SetAlgoLined
\DontPrintSemicolon
\KwIn{
    $N,M$: scalars, number of particles and number of inducing points\\
    \hspace*{35pt}$t_k,t_{k+1}, dt$: scalars, initial and final timepoints, and discretisation step\\
    \hspace*{35pt}$\boldsymbol{\mathcal{{O}}}_k, \boldsymbol{\mathcal{{O}}}_{k+1}$: $1 \times d$, $1 \times d$ initial and target state \\
    \hspace*{35pt}$\hat{\mathbf{f}}$: current drift estimate\\
    \hspace*{35pt}$\sigma$: noise amplitude\\
    \hspace*{35pt}$\boldsymbol{\gamma}_t$: geodesic curve in functional representation \\
 }

\KwOut{$F$: $d \times N\times (t_{k+1}-t_k)/dt$, sample representation of $q_t(x)$\\
 }
\BlankLine

$\ell = \frac{(t_{k+1}-t_{k})}{\de t}$   \tcp*{number of timesteps,}

\tcp*{, Forward filtering $\rho_t(x)$ (Eq.~\ref{eq:FPE2})         \,\, \;\;         }

$\epsilon = 10^{-3}$

$\mathbf{Z}_{ti=0} = \mathcal{N}(\boldsymbol{\mathcal{{O}}}_k, \epsilon \,\mathbf{I}_d)$   \tcp*{initialise particles' positions}

$\mathbf{Z}_{ti=1} = \mathbf{Z}_{0} + \de t \, \left( \hat{\mathbf{f}}(\mathbf{Z}_{0},t_0) - \frac{1}{2}\sigma^2 \frac{\mathbf{Z}_0-\boldsymbol{\mathcal{{O}}}_k  }{\epsilon } \right)$   \tcp*{1st step is with analytic score}
For $ti= 2:\ell$ \tcp*{deterministic propagation , }

 \hspace{36pt}$\mathbf{Z}_{ti+1} = \mathbf{Z}_{ti} + dt \left( \hat{\mathbf{f}}(\mathbf{Z}_{ti},t) -\frac{1}{2} \sigma^2 \nabla \log \rho(\mathbf{Z}_{ti};\mathbf{Z}_{ti})  \right) $   


\hspace{36pt} $W = \exp{\left(-U(\mathbf{Z}_{ti+1},t) \,dt \right) }$ \\
\hspace{36pt} $T^* =$ EnsembleTransformParticleFilter$(\mathbf{Z}_{ti+1},W)$\\
\hspace{36pt} $\mathbf{Z}_{ti+1} = \mathbf{Z}_{ti+1} \cdot T^*$\\

\tcp*{,Time-reversed propagation of flow $q_t(\mathbf{x})$             \hspace{155pt}  }
$\mathbf{B}_{ti=\ell} =  \mathcal{N}(\boldsymbol{\mathcal{{O}}}_{k+1}, \epsilon \,\mathbf{I}_d)$   \tcp*{initialise particles' positions }

\tcp*{1st step is stochastic }

$\mathbf{B}_{ti=\ell-1} = \mathbf{B}_{\ell} - \de t \left( \hat{\mathbf{f}}(\mathbf{B}_{\ell},t_1)    +\frac{1}{2} \sigma^2 \nabla \log \rho(\mathbf{B}_{\ell};\mathbf{Z}_{\ell})  - \frac{1}{2}\sigma^2 \frac{\mathbf{B}_\ell-\boldsymbol{\mathcal{{O}}}_{k+1}  }{\epsilon }\right)$   

For $ti= \ell-2:0$ \tcp*{deterministic propagation , }
\hspace{36pt} $\mathbf{B}_{ti-1} = \mathbf{B}_{ti} - \de t \left( \hat{\mathbf{f}}(\mathbf{B}_{ti},t) -\frac{1}{2} \sigma^2 \nabla \log \rho(\mathbf{B}_{ti}, \mathbf{Z}_{ti}) +\frac{1}{2} \sigma^2 \nabla \log q(\mathbf{B}_{ti},\mathbf{B}_{ti})  \right) $

\tcp*{,Compute control $u(x\,t)$ and controlled paths $\mathbf{F}_{0:T}$            \hspace{305pt}  }
For $ti= 1:\ell$

\hspace{36pt} $\mathbf{u}(\mathbf{x},ti) = \sigma^2 \nabla \log q(\mathbf{x};\mathbf{B}_{ti}) - \sigma^2 \nabla \log \rho(\mathbf{x};\mathbf{Z}_{ti}) $

\hspace{36pt} $\mathbf{F}_{ti+1} =\mathbf{F}_{ti} + \de t \,\left(\hat{\mathbf{f}}(\mathbf{F}_{ti},t) + \mathbf{u}(\mathbf{F}_{ti},ti)    - \frac{1}{2}\sigma^2 \frac{\mathbf{F}_0-\boldsymbol{\mathcal{{O}}}_k  }{\epsilon } \right) $

\caption{{Path augmentation algorithm employing Deterministic Particle Flow control} } \label{algo:control}
\end{algorithm}

With the notation $\nabla \log q(\mathbf{x};\mathbf{B}_{ti})$ we indicate the score function estimation in a functional form ($\mathbf{x}$) based on the density represented by the particles $\mathbf{B}_{ti}$, while $\nabla \log q(\mathbf{F}_{ti};\mathbf{B}_{ti})$ indicates the same score function evaluated at locations $\mathbf{F}_{ti}$.

\begin{algorithm}[H]\label{alg:gp_path_likelihood_1}
\SetAlgoLined
\DontPrintSemicolon
\KwIn{
    $\mathcal{Z} = \{\mathbf{z}_i\}_{i=1}^{S}$: inducing points for the sparse GP (Sp)\\
    \hspace*{35pt}$\{\mathbf{X}_j(t_\ell)\}_{j=1,\dots,N}^{\ell=1,\dots,T'}$: particle positions from the path measure $Q$ (BALL2)\\
    \hspace*{35pt}$\{\mathbf{g}(\mathbf{X}_j(t_\ell), t_\ell)\}$: effective drift evaluated along particles (gbALL2)\\
    \hspace*{35pt}$k^{\mathbf{f}}$: kernel with lengthscales $\ell_1,\ell_2,\ell_3$ (shared across dimensions)\\
    \hspace*{35pt}$g$: diffusion amplitude, $\sigma^2 = g^2$\\
    \hspace*{35pt}$\Delta t$: time step of the particle simulation\\
    \hspace*{35pt}$d$: state dimension, $N$: number of particles, $T'$: number of time steps
}
\KwOut{
    Approximations $I_1^{(i)}, I_2^{(i)}$ of the integrals over $A(\mathbf{x})$ and $B(\mathbf{x})$\\
}
\BlankLine

\tcp*{0. shorthand and initialisation}
Set $S \leftarrow |\mathcal{Z}|$ (number of inducing points)\;
Initialise $I_1 \in \mathbb{R}^{S \times S \times d}$ and $I_2 \in \mathbb{R}^{S \times d}$ to zero\;
Initialise $\Lambda \in \mathbb{R}^{S \times S \times d}$ and $\mathbf{d} \in \mathbb{R}^{S \times d}$ to zero\;

\BlankLine
\tcp*{1. compute kernel matrices on the inducing points}
Construct the inducing–inducing kernel matrix
\begin{equation}
\mathcal{K}_S = k^{\mathbf{f}}(\mathcal{Z},\mathcal{Z}) \in \mathbb{R}^{S \times S}
\end{equation}
and compute a regularised inverse
\begin{equation}
\mathcal{K}_S^{-1} = \big(\mathcal{K}_S + \varepsilon I\big)^{-1},\quad \varepsilon \approx 10^{-3}.
\end{equation}

Define the kernel map to inducing points
\begin{equation}
k^{\mathbf{f}}(\mathcal{Z}, \mathbf{x}) 
= \big( k^{\mathbf{f}}(\mathbf{z}_i, \mathbf{x}) \big)_{i=1}^S \in \mathbb{R}^S.
\end{equation}

\BlankLine
\tcp*{2. sample-based approximation of $A(\mathbf{x})$ and $B(\mathbf{x})$}
\For{$i = 1,\dots,d$}{  \tcp*{loop over state dimensions}
    \For{$\ell = 1,\dots,T'$}{ \tcp*{loop over time}
        Let $\mathbf{X}(t_\ell) \in \mathbb{R}^{d \times N}$ be the particle positions at time $t_\ell$\;
        For each particle position $\mathbf{X}_j(t_\ell)$, compute
        \begin{equation}
        \mathbf{k}_j = k^{\mathbf{f}}(\mathcal{Z}, \mathbf{X}_j(t_\ell)) \in \mathbb{R}^S.
        \end{equation}
        Stack them column-wise to obtain
        \begin{equation}
        K_\ell = \big[\mathbf{k}_1,\dots,\mathbf{k}_N\big] \in \mathbb{R}^{S \times N}.
        \end{equation}

        Let $g_i(\mathbf{X}_j(t_\ell), t_\ell)$ denote the $i$-th component of the effective drift at particle $j$ and time $t_\ell$\;

        \tcp*{accumulate Monte Carlo estimates of the integrals}
        Update
        \begin{equation}
        I_1^{(i)} \leftarrow I_1^{(i)} + K_\ell K_\ell^\top,
        \qquad
        I_2^{(i)} \leftarrow I_2^{(i)} + K_\ell \mathbf{g}_i(t_\ell),
        \end{equation}
        where $\mathbf{g}_i(t_\ell) = \big(g_i(\mathbf{X}_1(t_\ell),t_\ell),\dots,g_i(\mathbf{X}_N(t_\ell),t_\ell)\big)^\top$.
    }
    \tcp*{normalise by time and number of particles}
    \begin{equation}
    I_1^{(i)} \leftarrow \frac{\Delta t}{N} I_1^{(i)}, \qquad
    I_2^{(i)} \leftarrow \frac{\Delta t}{N} I_2^{(i)}.
    \end{equation}
}

\caption{{Gaussian process drift inference from an augmented path measure (part I)}}\label{algo:drift}
\end{algorithm}

In this algorithm Here $I_1^{(i)}$ approximates
$
\int k^{\mathbf{f}}(\mathcal{Z},\mathbf{x}) A(\mathbf{x}) k^{\mathbf{f}}(\mathbf{x},\mathcal{Z}) \,\de \mathbf{x}$, 
and $I_2^{(i)}$ approximates
$
\int k^{\mathbf{f}}(\mathcal{Z},\mathbf{x}) B_i(\mathbf{x}) \,\de \mathbf{x}.$

\begin{algorithm}[H]\label{alg:gp_path_likelihood_2}
\SetAlgoLined
\DontPrintSemicolon
\KwIn{
    Same inputs as Alg.~\ref{alg:gp_path_likelihood_1}\\
    \hspace*{35pt}$I_1^{(i)}, I_2^{(i)}$: Monte Carlo approximations from Alg.~\ref{alg:gp_path_likelihood_1}\\
    \hspace*{35pt}$\mathcal{K}_S, \mathcal{K}_S^{-1}$: inducing–inducing kernel matrix and its regularised inverse
}
\KwOut{
    Component-wise drift estimators $\hat{f}_i(\mathbf{x})$, $i=1,\dots,d$\\
    \hspace*{35pt}Expected negative log data likelihood $\mathcal{L}_{\text{path}}$ under $Q_f$
}
\BlankLine

\tcp*{3. compute $\Lambda$ and $\mathbf{d}$ for each component}
\For{$i = 1,\dots,d$}{
    \tcp*{match Eq.~\eqref{appeq:drift_measure} with sparse GP parametrisation}
    \begin{equation}
    \Lambda^{(i)} \leftarrow \frac{1}{\sigma^2} \,\mathcal{K}_S^{-1}\, I_1^{(i)}\, \mathcal{K}_S^{-1},
    \qquad
    \mathbf{d}^{(i)} \leftarrow \frac{1}{\sigma^2}\, \mathcal{K}_S^{-1}\, I_2^{(i)}.
    \end{equation}
}
This matches the definitions
\begin{equation}
\Lambda = \frac{1}{\sigma^2} \mathcal{K}_S^{-1} 
\left(   \int k^{\mathbf{f}}(\mathcal{Z},\mathbf{x}) A(\mathbf{x}) k^{\mathbf{f}}(\mathbf{x},\mathcal{Z}) \de \mathbf{x} \right) 
\mathcal{K}_S^{-1},
\quad
\mathbf{d} = \frac{1}{\sigma^2} \mathcal{K}_S^{-1} 
\left(   \int k^{\mathbf{f}}(\mathcal{Z},\mathbf{x}) B(\mathbf{x}) \de \mathbf{x} \right).
\end{equation}

\BlankLine
\tcp*{4. define the component-wise drift estimators}
For each component $i = 1,\dots,d$, define
\begin{equation}
\hat{f}_i(\mathbf{x})
= k^{\mathbf{f}}(\mathbf{x},\mathcal{Z})
\Big(I + \Lambda^{(i)} \mathcal{K}_S\Big)^{-1}\mathbf{d}^{(i)},
\end{equation}
so that the full drift estimate is
\begin{equation}
\hat{\mathbf{f}}_S(\mathbf{x}) = 
\big(\hat{f}_1(\mathbf{x}), \dots, \hat{f}_d(\mathbf{x})\big)^\top.
\end{equation}

\BlankLine
\tcp*{5. compute expected negative log data likelihood under $Q_f$}
Initialise accumulators $S_{\|f\|} \leftarrow 0$, $S_{\nabla\cdot f} \leftarrow 0$, $S_{f\cdot g} \leftarrow 0$\;

\For{$\ell = 1,\dots,T'$}{
    For all particle positions $\mathbf{X}_j(t_\ell)$, evaluate $\hat{\mathbf{f}}_S(\mathbf{X}_j(t_\ell))$\;
    Accumulate
    \begin{equation}
    S_{\|f\|} \leftarrow S_{\|f\|} + \sum_{j=1}^{N} \|\hat{\mathbf{f}}_S(\mathbf{X}_j(t_\ell))\|^2,
    \end{equation}
    \begin{equation}
    S_{f\cdot g} \leftarrow S_{f\cdot g} + \sum_{j=1}^{N} 
        \hat{\mathbf{f}}_S(\mathbf{X}_j(t_\ell))^\top 
        \mathbf{g}(\mathbf{X}_j(t_\ell), t_\ell),
    \end{equation}
    and compute the trace of the Jacobian 
    $\nabla \cdot \hat{\mathbf{f}}_S(\mathbf{X}_j(t_\ell))$ via automatic differentiation, accumulating it into $S_{\nabla\cdot f}$\;
}
Approximate the expected negative log data likelihood (up to constants) as
\begin{equation}
\mathcal{L}_{\text{path}}
= \frac{\Delta t}{N}
\left(
\tfrac{1}{2} S_{\|f\|} + S_{\nabla\cdot f} + S_{f\cdot g}
\right),
\end{equation}
which corresponds to evaluating the quadratic form in Eq.~\eqref{appeq:drift_measure} under the approximate posterior $Q_f$.

\caption{{Gaussian process drift inference from an augmented path measure (part II)}}
\end{algorithm}

\newpage

\end{document}

%% file: math_commands.tex
\usepackage{amsmath,amsfonts,bm}

\def\eqref#1{equation~\ref{#1}}

\def\1{\bm{1}}

\DeclareMathAlphabet{\mathsfit}{\encodingdefault}{\sfdefault}{m}{sl}
\SetMathAlphabet{\mathsfit}{bold}{\encodingdefault}{\sfdefault}{bx}{n}



%% file: geometry.tex
 Our method is based on the argument that
the invariant density\footnote{In the following the discussion concentrates around invariant measures. We point out here that the invariant density is the Radon-Nikodym derivative of the invariant measure with respect to some reference measure, often the Lebesgue measure if it exists~\citep{maharam1969invariant}.}
of the observed system imposes a low-dimensional structure on the state space, within which the observations are confined. We propose that this low-dimensional structure is well approximated by a Riemannian manifold ${\mathcal{M}_{\infty} \in \mathbb{R}^{m \leq d}}$ and that the observations $ \{\boldsymbol{\mathcal{O}}_k\}_{k=1}^{K}
$ offer a reliable discrete approximation to $\mathcal{M}_{\infty}$.

 We employ the notion of a "low-dimensional structure" as a concise way to refer to the fact that for many dissipative dynamical systems, the invariant measure has support on a subset of the state space with dimension
smaller than the ambient space dimension. This phenomenon arises due to the dissipative nature of these systems, which causes volume contraction in the state space, resulting in trajectories concentrating asymptotically on attractors of lower dimension than the state space dimension.
To provide further justification on this, in the following section, we start by building intuition from deterministic dynamical systems and then generalise to stochastic dynamics.

\subsection{Dimensionality of invariant measures induced by deterministic dynamics}
\label{app:low_dim_invariant}

We consider a dissipative deterministic dynamical system of the form
\begin{equation}
  \dot{\mathbf{x}}_t = \mathbf{f}(\mathbf{x}_t), \qquad \mathbf{x}_t \in \mathbb{R}^d,
\end{equation}
generating a semiflow $(\Phi^t)_{t \ge 0}$. 
Under standard assumptions, the dynamics
admit an invariant probability measure $\mu$ describing the distribution
of states along long-term typical trajectories. From an ergodic perspective, $\mu$ is the natural object
characterising the asymptotic behaviour of the system. For almost every initial condition in $\mu$, the
empirical measure
\begin{equation}
  \frac{1}{T}\int_0^T \delta_{\mathbf{x}_t}\,\de t
\end{equation}
converges (in the weak sense) to~$\mu$.

For dissipative systems, phase–space volumes contract along typical trajectories, so the
Lebesgue measure is not invariant under the dynamics, i.e. state space volume is not preserved when pushed forward through the flow~\citep{ruelle1979ergodic}. This implies that the system state does not explore the ambient space uniformly. Instead, trajectories
concentrate asymptotically on subsets of state space of vanishing Lebesgue measure.
In fact, this concentration phenomenon persists even in chaotic systems, where, although trajectories separate
exponentially along unstable directions, contraction along stable directions dominates
the evolution of infinitesimal volumes in the state space.

The resulting invariant measure $\mu$ typically has an \textbf{effective dimension} smaller than the ambient space dimension. To quantify this, we
require a notion of dimensionality that remains meaningful when the Lebesgue measure
vanishes. The Hausdorff dimension~\citep{ruelle1989chaotic,Young2002SRB,ott2002chaos} lends itself for such a purpose since it naturally extends from
sets to probability measures~\citep{young1982dimension}. More precisely, the Hausdorff dimension of an invariant
measure $\mu$ is defined as the smallest Hausdorff dimension among all measurable sets containing 
$\mu$
\begin{equation}
  \dim_H(\mu)
  =
  \inf\bigl\{\dim_H(\mathcal{A})\;:\;\mu(\mathcal{A})=1\bigr\}.
\end{equation}

A useful aspect of this formulation is its local interpretation. Under mild regularity
assumptions, $\dim_H(\mu)$ can be characterised by the scaling of probability
mass around typical points under $\mu$. If, for almost every $\mathbf{x}$,
\begin{equation}
  \mu\!\left(B_\varepsilon(\mathbf{x})\right)
  \sim \varepsilon^{\,d_\mu}
  \qquad \text{as }\varepsilon\to 0,
\end{equation}
then $d_\mu = \dim_H(\mu)$. Thus, this dimension reflects how probability mass
concentrates across scales.

In (smooth) deterministic dynamical systems, the
interplay between expansion and contraction along different directions governs this local scaling behaviour. This is well characterised by 
Lyapunov exponents that quantify the exponential deformation of infinitesimal neighbourhoods,
while the metric (Kolmogorov-Sinai) entropy $h_\mu$  quantifies the rate at which trajectories generate information.
Well known results in ergodic theory~\citep{ledrappier1985metric} show that
the Hausdorff dimension of an invariant measure can be expressed directly in terms of
these quantities, and \textbf{is strictly smaller than the ambient space dimension $d$ in dissipative systems} with non-trivial Lyapunov exponents, i.e. both positive and negative exponents.

More precisely, 
according to the Oseledets’ theorem~\citep{oseledets1968multiplicative}, the system has a Lyapunov spectrum $\lambda_1 \geq \cdots \geq \lambda_d$, and
dissipativity implies on average volume contraction, i.e.
\begin{equation}\label{app:dissi}
\sum_{i=1}^d \lambda_i < 0.
\end{equation}

\cite{ledrappier1985metric}
formulate an expression for the Hausdorff dimension of the invariant measure $\mu$ in terms of the Lyapunov exponents $\{\lambda_i\}^d_{i=1}$ and the Kolmogorov-Sinai entropy $h_\mu$~\citep{barany2017ledrappier}
\begin{equation}\label{appeq:youngDim}
\dim_H(\mu) = k + \frac{h_\mu - \sum_{i=1}^k \lambda_i}{|\lambda_{k+1}|},
\end{equation}
where $k$ is the largest integer for which $\sum_{i=1}^k \lambda_i \geq h_\mu$. This relation holds under standard smoothness and hyperbolicity
assumptions (for instance for $C^{1+\alpha}$ systems with
non-zero Lyapunov exponents almost everywhere). Intuitively, $k$ here quantifies the number of expanding dimensions needed to characterise the system's entropy. 

Since the sum of all Lyapunov exponents is negative (Eq.~\ref{app:dissi}), and the 
metric entropy is bounded by the sum of positive Lyapunov exponents~\citep{ruelle1978inequality}
\begin{equation}
 0 \leq h_\mu \leq \sum_{\lambda_i > 0} \lambda_i,
\end{equation}
the equality of Eq.~\ref{appeq:youngDim} implies
\begin{equation}
\dim_H(\mu) < d,
\end{equation}
indicating that
the invariant measure concentrates on a subset of the state space, whose
Hausdorff dimension is strictly smaller than the ambient space dimension $d$.

\subsection{Dimensionality of invariant measures induced by stochastic dynamics}
\label{app:stochastic_dimension}

We now consider stochastic dynamical systems of the form
\begin{equation}
\de \mathbf{X}_t = \mathbf{f}(\mathbf{X}_t)\de t + \sigma\,\de \mathbf{W}_t ,
\label{eq:sde_appendix}
\end{equation}
similar to the systems we discuss in the main text.
Under mild conditions on $\mathbf{f}$ and $\sigma$,
the corresponding Markov semigroup admits a unique invariant probability measure
$\mu_\sigma$, which coincides with the stationary solution of the associated
Fokker–Planck equation~\citep{risken1996fokker}.

The additive noise regularises the deterministic invariant measure, nevertheless its density concentrates exponentially around $\mathcal{A}$ as $\sigma \rightarrow 0$.
For non-degenerate noise $\sigma$, the H\"ormander condition ensures that $\mu_\sigma$ is
absolutely continuous with respect to the Lebesgue measure, and thus possesses a
\emph{smooth} invariant density~\citep{hormander1967hypoelliptic}. 
However, the invariant measure $\mu_\sigma$ of the stochastic system of Eq.~\ref{eq:sde_appendix} satisfies the following exponential concentration inequality around the deterministic attractor $\mathcal{A}$ \textbf{for sufficiently small noise amplitude $\sigma$}
\begin{equation}
\mu_\sigma\bigl(\{ \mathbf{x} \in \mathbb{R}^d : \mathrm{dist}(\mathbf{x}, \mathcal{A}) > \delta \} \bigr) \leq C(\delta) \exp\left(- \frac{c(\delta)}{\sigma^2} \right),
\end{equation}
for all $\delta > 0$, where $C(\delta), c(\delta) > 0$ denote $\delta$-dependent constants, that are nevertheless independent of noise amplitude $\sigma$ (see Theorem 4.2.1~\citep{kifer1988random}).
This exponential concentration indicates that, although $\mu_\sigma$ is absolutely continuous with respect to the Lebesgue measure for $\sigma > 0$, it becomes increasingly confined near $\mathcal{A}$ as $\sigma \to 0$. Thus the effective dimension of
$\mu_\sigma$ approaches that of the invariant measure of the deterministic system
$\mu_0$, while remaining bounded above by the ambient dimension $d$~\citep{kifer1988random, arnold2006random}.

In this sense, stochasticity does not destroy the low-dimensional structure induced by the deterministic dynamics, but
thickens the invariant measure around the
deterministic attractor. 
Thus even though $\mu_\sigma$ is smooth, its effective dimension can still be low-dimensional in the sense of mass being tightly concentrated near a lower-dimensional skeleton determined by the underlying deterministic dynamics.